\documentclass[preprint,12pt]{elsarticle}

\usepackage{amsmath, amsthm, amsfonts, amssymb}
\usepackage{upgreek}
\usepackage[font=small, labelsep=period]{caption}
\usepackage{longtable}
\usepackage{booktabs}
\usepackage{hyperref}
\usepackage[dvipsnames]{xcolor}
\usepackage{listings}
\allowdisplaybreaks
\theoremstyle{plain}
\newtheorem{thm}{Theorem}[section]
\newtheorem{prop}[thm]{Proposition}
\newtheorem{coro}[thm]{Corollary}
\newtheorem{lemm}[thm]{Lemma}
\newtheorem*{thms}{Main Theorem}

\theoremstyle{definition}
\newtheorem{deff}{Definition}
\newtheorem{examp}{Example}

\theoremstyle{remark}
\newtheorem{rema}{Remark}

\newcommand\twoscript[2]{\substack{{#1} \\ {#2}}}
\newcommand\threescript[3]{\substack{{#1} \\ {#2} \\ {#3}}}

\newcommand\rmi{\mathrm{i}}
\newcommand\rme{\mathrm{e}}
\newcommand\ob{\overline{\beta}}

\newcommand*\legendre[2]{\genfrac{(}{)}{}{}{#1}{#2}}
\newcommand*\tbtmat[4]{\left(\begin{smallmatrix}{#1} & {#2} \\ {#3} & {#4}\end{smallmatrix}\right)}
\newcommand*\tbtMat[4]{\begin{pmatrix}{#1} & {#2} \\ {#3} & {#4}\end{pmatrix}}
\newcommand*\abs[1]{\lvert#1\rvert}
\newcommand\etp[1]{\mathfrak{e}\left(#1\right)}
\newcommand\ord{\mathop{\mathrm{ord}}}
\newcommand\rad{\mathop{\mathrm{rad}}}
\newcommand\radE{\mathop{\mathrm{rad}_E}}
\newcommand\radO{\mathop{\mathrm{rad}_O}}
\newcommand\radp{\mathop{\mathrm{rad}'}}
\newcommand\irad{\mathop{\mathrm{irad}}}
\newcommand\iradp{\mathop{\mathrm{irad}'}}

\newcommand\numZ{\mathbb{Z}}
\newcommand\numQ{\mathbb{Q}}

\newcommand\numC{\mathbb{C}}
\newcommand\projQ{\mathbb{P}^1(\mathbb{Q})}

\newcommand\halfint{\frac{1}{2}\mathbb{Z}}
\newcommand*\numgeq[2]{\mathbb{#1}_{\geq #2}}

\newcommand\slZ{\mathrm{SL}_2(\mathbb{Z})}

\newcommand\pslZ{\mathrm{PSL}_2(\mathbb{Z})}
\newcommand\glpQ{\mathrm{GL}_2^{+}(\mathbb{Q})}
\newcommand\glpR{\mathrm{GL}_2^{+}(\mathbb{R})}
\newcommand\glptQ{\widetilde{\mathrm{GL}_2^{+}(\mathbb{Q})}}
\newcommand\glptR{\widetilde{\mathrm{GL}_2^{+}(\mathbb{R})}}

\newcommand\sltZ{\widetilde{\mathrm{SL}_2(\mathbb{Z})}}

\newcommand*\Gz[1]{\Gamma_0({#1})}
\newcommand*\Gzt[1]{\widetilde{\Gamma_0({#1})}}
\newcommand*\Gzp[1]{\overline{\Gamma_0({#1})}}

\newcommand\uhp{\mathfrak{H}}

\newcommand\MFormW[3]{M^!_{#1}(\Gzt{#2},#3)}
\newcommand\MForm[3]{M_{#1}(\Gzt{#2},#3)}
\newcommand\MFormt[2]{M_{#1}(\Gzt{#2})}
\newcommand\CForm[3]{S_{#1}(\Gzt{#2},#3)}
\newcommand\MFormWG[3]{M^!_{#1}(\widetilde{#2},#3)}
\newcommand\MFormG[3]{M_{#1}(\widetilde{#2},#3)}
\newcommand\CFormG[3]{S_{#1}(\widetilde{#2},#3)}

\journal{Journal of Number Theory}

\begin{document}

\begin{frontmatter}



\title{Double coset operators and eta-quotients}


\author[label1]{Hai-Gang Zhou}
\ead{haigangz@tongji.edu.cn}
\ead[https://orcid.org/0000-0002-4711-1301]{https://orcid.org/0000-0002-4711-1301}

\author[label1]{Xiao-Jie Zhu\corref{cor}}
\ead{zhuxiaojiemath@outlook.com}
\ead[https://orcid.org/0000-0002-6733-0755]{https://orcid.org/0000-0002-6733-0755}

\cortext[cor]{Corresponding author}

\affiliation[label1]{organization={School of Mathematical Sciences, Tongji University},
            addressline={1239 Siping Road}, 
            city={Shanghai},
            postcode={200092}, 
            country={P.R. China}}

\begin{abstract}
We study a type of generalized double coset operators which may change the characters of modular forms. For any pair of characters $v_1$ and $v_2$, we describe explicitly those operators mapping modular forms of character $v_1$ to those of $v_2$. We give three applications, concerned with eta-quotients. For the first application, we give many pairs of eta-quotients of small weights and levels, such that there are operators maps one eta-quotient to another. We also find out these operators. For the second application, we apply the operators to eta-powers whose exponents are positive integers not greater than $24$. This results in recursive formulas of the coefficients of these functions, generalizing Newman's theorem. For the third application, we describe a criterion and an algorithm of whether and how an eta-power of arbitrary integral exponent can be expressed as a linear combination of certain eta-quotients.
\end{abstract}




\end{frontmatter}

\tableofcontents

\section{Introduction}
\label{sec:Intro}
\subsection{Hecke operators and Double coset operators}
\label{subsec:Operators}
The \textit{discriminant function} $\Delta(z)$ can be defined by the infinite product
\begin{equation}
\Delta(z) = q\prod_{n\geq1}(1-q^n)^{24},
\end{equation}
where $q=\rme^{2\uppi\rmi z}$ with $z\in\numC$ such that $\Im(z)>0$. The \textit{Ramanujan $\tau$-function} is defined to be the Fourier coefficients of the $\Delta$-function, that is,
\begin{equation}
\Delta(z)=\sum_{n\geq1}\tau(n)q^n.
\end{equation}
Ramanujan conjectured some properties of $\tau(n)$ in \cite{Ram00}, some of which were proved by Mordell in 1917. Hecke developed a general theory of a certain kind of linear operators on modular forms in \cite{Hec37a} and \cite{Hec37b}, which could settle some of Ramanujan's conjectures on $\tau(n)$ and a wide range of similar problems. It is a generalization of such operators which is the main subject of this paper.

Roughly speaking, a Hecke operator maps a modular form $f$ to a linear combination of functions of the form $f\left(\frac{az+b}{d}\right)$. The resulted function is still a modular form with the same transformation properties as $f$. So Hecke operators may be viewed as operators from certain spaces of modular forms to themselves. In 1957, Wohlfahrt \cite{Woh57} defined a type of Hecke operators, but with the property of changing multiplier systems of modular forms. Wohlfahrt applied these operators to certain theta-functions defined in \cite{Sch39}, and obtained formulas concerning representations of integers by quadratic forms of odd number of variables. In \cite{vL57} the author applied this type of Hecke operators to eta-quotients and, for instance, gave independent proofs of many of Newman’s relations on Fourier coefficients of eta-powers (\cite{New55}, \cite{New56} and \cite{New58}). He also gave new relations and studied linear dependency of functions occurring in the definition of generalized Hecke operators. For the definitions of Dedekind eta-function, eta-quotients and eta-powers, see \eqref{eq:defEta}, the paragraph containing \eqref{eq:etaQuoChar}, and the paragraph containing \eqref{eq:Prn} respectively.

A general theory of Hecke operators acting on modular forms of half-integral weights was developed by Shimura \cite{Shi73}. He explored various aspects of these operators which have integral-weights parallels developed by Hecke. For instance, Shimura studied functional equations of L-functions associated with modular forms, Euler products and eigenfunctions. More importantly, Shimura found a deep relation between modular forms of half-integral weights and those of integral weights. Notice that Shimura's operators do not change the multiplier system (more precisely, change very little). In fact, modular forms of half-integral weights in the sense of Shimura are just those with some special multiplier systems---the multiplier system of the classical theta series ($\theta(\tau)$ in Example \ref{examp:onedimExamp2}) times a Dirichlet character (see the paragraph before Theorem \ref{thm:charCompatibleEtaDirichlet}).

The theory of Hecke operators (developed by Hecke for integral weights and by Shimura for half-integral weights) play important roles in arithmetical properties of Fourier coefficients of modular forms, and also in structures of spaces of modular forms. (e.g., these operators lead to natural decomposition of spaces of modular forms into eigenspaces of all the operators.) For example, Ono \cite{Ono11} obtained the action of Hecke operators on generating functions of the partition function. More precisely, Ono derived a closed formula of half-integral Hecke operators $T_{l^2}$ acting on $\eta^{-1}(24z)$. More recently, Carney, Etropolski and Pitman \cite{CEP12} obtained closed formulas of half-integral Hecke operators $T_{l^2}$ acting on eta-powers $\eta^r$ with the variable $z$ replaced by some $Nz$ of suitable integer $N$. Such closed formulas are useful in studying congruence properties of Fourier coefficients of eta-powers. It should be mentioned that, Str\"{o}mberg \cite{Str07,Str08} has succeeded in applying Hecke-like operators based on the ideas of Wohlfahrt to Maass Waveforms with eta multiplier.

By contrast, those operators that change multiplier systems rarely have applications.

The central objects of this paper are the generalized double coset operators. Notice that Hecke operators, not only on spaces of modular forms and modular functions, but also on various kinds of automorphic forms (such as Jacobi forms, see \cite{EZ85} and \cite{Ajo15}, and Siegel modular forms, see \cite{HW02}), have two relevant definitions: one using one-side cosets of modular groups, and the other using double cosets. The latter definition is usually more general than the former. The generalized Hecke operators of Wohlfahrt are defined using one-side cosets. It may therefore be interesting to see how we can develop
a similar theory using double coset operators\footnote{The approach by Wohlfahrt is essentially the same as the double coset operators for square-free integers, as is pointed out to us by the referee.}. This paper contains such an attempt (Section \ref{sec:Operators}).

Among the properties of such generalized double coset operators we proved, the most important one is Theorem \ref{thm:charCompatibleEtaDirichlet}, which is given below as `Main Theorem', using notation introduced in Sections \ref{sec:Modular forms} and \ref{sec:Operators}.
\begin{thms}
Let $v_1$ and $v_2$ be complex linear characters on the modular group $\Gz{N}$ or its double cover, each of which is the product of the character of an eta-quotient and a Dirichlet character. Let $l$ be a positive integer. Then we can define an operator $T_l$ that maps modular forms of character $v_1$ to those of character $v_2$ if and only if $l$ satisfies certain explicit congruence conditions, i.e. \eqref{eq:charCom1}, \eqref{eq:charCom2}, \eqref{eq:charCom3} and \eqref{eq:charCom4}.
\end{thms}

\subsection{Dedekind Eta-function and Eta-quotients}
The Dedekind eta-function\footnote{It is unfortunate that a usual variable of modular forms is $\tau$, which is used widely in literature. But this contradicts the symbol for Ramanujan $\tau$-function. Since we never use Ramanujan $\tau$-function in the remaining of the text, we shall use the symbol $\tau$ to denote the variable of modular forms, instead of $z$.} $\eta(\tau)$ is defined as
\begin{equation*}
\eta(\tau)=q^{1/24}\prod_{n\geq 1}(1-q^n),
\end{equation*}
where $q=\rme^{2\uppi\rmi \tau}$ with $\tau\in\numC$ such that $\Im(\tau)>0$. Note that $\Delta(\tau)=\eta^{24}(\tau)$. It was Dedekind who first (in 1877) introduced this function, proved modularity properties, and gave an explicit formula of its character. See \cite[Section 3]{Apo90} for all of these. This function is useful, since many modular forms or functions can be expressed using it, and many combinatorial generating functions can be constructed from it. More precisely, many modular forms can be expressed as a linear combination of some eta-quotients --- products of some $\eta^{r_j}(j\tau)$ (see the paragraph containing \eqref{eq:etaQuoChar}). For example, let $r_k(n)$ be the number of solutions of the equation $n=x_1^2+x_2^2+\dots+x_k^2$ with $x_1,\,x_2,\dots,x_k\in\numZ$. Then we have
\begin{equation}
\sum_{n\geq 0}r_k(n)q^n=\frac{\eta^{5k}(2\tau)}{\eta^{2k}(\tau)\eta^{2k}(4\tau)}.
\end{equation}
The only difficulty for a proof is the case $k=1$. We have sketched a proof for this case in Example \ref{examp:onedimExamp2}. For another example (as Ono observed), let $E_4(\tau)$ and $E_6(\tau)$ be the standard Eisenstein series for the full modular group, normalised such that the constant term is $1$. Then
{\small
\begin{gather}
E_4(\tau)=\frac{\eta^{16}(\tau)}{\eta^8(2\tau)}+2^8\cdot\frac{\eta^{16}(2\tau)}{\eta^8(\tau)},\\
E_6(\tau)=\frac{\eta^{24}(\tau)}{\eta^{12}(2\tau)}-2^5\cdot3\cdot5\cdot\eta^{12}(2\tau)-2^9\cdot3\cdot11\cdot\frac{\eta^{12}(2\tau)\eta^8(4\tau)}{\eta^8(\tau)}+2^{13}\frac{\eta^{24}(4\tau)}{\eta^{12}(2\tau)}.
\end{gather}
}%
See \cite[Theorem 1.67]{Ono04}. This implies that the algebra of modular forms on the full modular group with trivial character is contained in the algebra generated by $\eta(\tau)$, $\eta(2\tau)$, and $\eta(4\tau)$. Other important results on the classification of spaces of modular forms which are generated by eta-quotients have been obtained in \cite{RW15}, \cite{All20} and \cite{Bha17}.

Eta-quotients have connections with many other areas of mathematics. For the connection with theta functions, see \cite[p. 30]{BvdGHZ08} and \cite{LO13}. For that with the partition function, see \cite{Ono00}. For that with elliptic curves, see \cite{MO97}.

It is natural to consider Hecke operators acting on eta-quotients. In this direction, van Lint \cite{vL57} has proved that the eta-power $\eta^r(\tau)$ is a Hecke eigenform for any Hecke operators in the sense of Wohlfahrt, for $r=1,\,2,\dots,24$. Martin has succeeded in figuring out all integral-weight holomorphic eta-quotients that are simultaneous eigenforms for all possible Hecke operators in \cite{Mar96}. But the eta-quotients Martin considered are those with characters induced by Dirichlet characters, and the Hecke operators Martin used are in the sense of Hecke himself.

In this paper, after developing the theory of generalized double coset operators in Section 3, we will give three applications concerning eta-quotients. Before describing these applications, we emphasize that the operators $T_l$ (see Definition \ref{deff:doubleCosetOperators} and the paragraph after the proof of Lemma \ref{lemm:matrixSetDecomp}) we used may act on any eta-quotients, not only those with characters induced by Dirichlet characters. This distinguishes our treatment from previous ones.

For \emph{the first application}, we obtain formulas of the form
\begin{equation}
\label{eq:app1Rough}
T_l \prod_{j}\eta^{r_j}(j\tau) = c\cdot\prod_{j'}\eta^{r_{j'}}(j'\tau),
\end{equation}
where $r_j$ and $r_{j'}$ ($j,\,j'\geq1$) are integers and only finitely many of them are nonzero. For the precise version, see Theorem \ref{thm:TlIntegralWeight} and Theorem \ref{thm:TlHalfintegralWeight}. For concrete examples, see Example \ref{examp:onedimExamp1} and Example \ref{examp:onedimExamp2}. We obtain 568 pairs of integral-weight eta-quotients, and 441 pairs of half-integral-weight ones for which \eqref{eq:app1Rough} holds for some $l$, and consequently for infinitely many $l$. In fact, we obtain more, but listing all of them would lead to a too long paper.

For \emph{the second application}, we generalize some recursive relations concerning Fourier coefficients of eta-powers $\eta^r(\tau)$ with $r=1,\,2,\dots,24$. Such relations were first discovered and p{}roved by Newman in \cite{New55}, \cite{New56} and \cite{New58}. The thesis \cite{vL57} improved and reproved these relations using (Wohlfahrt's) Hecke operators. Their results could be rewritten as $T_l\eta^r=c_r\eta^r$ with $c_r\in\numC$, for $l$ being some prime power. What we proved are formulas of this form for $l$ being general positive integers, not only prime powers, and we write out the explicit expression. See Theorem \ref{thm:mainThmEtar} for details.

For \emph{the third application}, we give a criterion for whether a certain kind of modular form can be expressed as linear combinations of some eta-quotients, and we give an algorithm to compute the coefficients. The modular forms under study are
\begin{equation*}
F_{r, p^{\beta}}(\tau)=q^{pr/24}\sum_{n \in\numZ}P_r\left(p^{\beta}n+\frac{r(p^{\beta+1}-1)}{24}\right)q^n,
\end{equation*}
where $p$ is a prime and $r$ is an integer satisfying $24 \mid r(p^2-1)$, and $\beta$ is a positive odd number. The quantity $P_r(n)$ occurring in these functions is the $n$-th Fourier coefficient of the function $q^{-r/24}\eta^r(\tau)$. Similarly, for even $\beta$, we also have a class of functions
\begin{equation*}
F_{r, p^{\beta}}(\tau)=q^{r/24}\sum_{n \in\numZ}P_r\left(p^{\beta}n+\frac{r(p^{\beta}-1)}{24}\right)q^n.
\end{equation*}
These functions arise naturally when one investigates congruence properties, vanishing properties and recursive relations of the sequences $P_r(n)$. See, for example, \cite[Section 5]{Ono04} for the case $r=-1$. If one can express $F_{r, p^{\beta}}(\tau)$ as a linear combination of some eta-quotients on the modular group $\Gz{p}$, then by passing to $\numZ/p\numZ$ for example one can find congruence relations modulo $p$. Recently, Du, Liu and Zhao \cite[Theorem 1.1]{DLZ19} give an existence result that $F_{r, p^{\beta}}(\tau)$ can be expressed as such a linear combination when $p=2,\,3,\,5,\,7,\,13$. Our result has the features that there may be some linear expansion for other $p$ (arbitrary large) and that the coefficients are given explicitly by some recursive relations. See Theorem \ref{thm:mainThmLinearCombina}.

\subsection{Outline of the Paper and Notational Convention}
The structure of this paper is as follows. In Section \ref{sec:Modular forms}, we recall basic facts on modular forms and eta-quotients that will be used later. In particular, we emphasize the necessity to use double covers of modular groups, and we give a proof of the valence formula suitable for our purpose. In Section \ref{sec:Operators}, we develop a theory of generalized double coset operators. We begin with some group-theoretical preliminaries. We then give the definition and basic properties. The section ends with a criterion (Theorem \ref{thm:charCompatibleEtaDirichlet}) on the existence of operators between two arbitrarily chosen characters. Section \ref{sec:one-dimensional spaces} gives the first application on finding those eta-quotients and operators such that the result of applying an operator on an eta-quotient is still an eta-quotient, up to a constant factor. The main results are Theorem \ref{thm:TlIntegralWeight} and Theorem \ref{thm:TlHalfintegralWeight}. Section \ref{sec:eta-powers} contains the second application, which is a generalization of Newman's Theorem. We give general recursive relations on Fourier coefficients of eta-powers $\eta^r$ with $r\in\{0,\,1,\,2,\dots,24\}$. The result is summarized in Theorem \ref{thm:mainThmEtar}. In Section \ref{sec:Express}, we give the third application, expressing a certain kind of modular form as linear combinations of eta-quotients of prime levels. The final theorem is a criterion (a sufficient condition) on whether such expressions exist, and also an algorithm on how to compute the coefficients. See Theorem \ref{thm:mainThmLinearCombina}. Section \ref{sec:Miscellaneous} contains some miscellaneous observations that do not fit into any previous sections. All tables, that are parts of the Theorem \ref{thm:TlIntegralWeight} and Theorem \ref{thm:TlHalfintegralWeight}, are arranged in \ref{apx:tables}.

We collect some notational conventions here. We shall always use $\tau$ to denote the variable of a modular form in next sections, so $\tau\in\uhp$, the complex numbers with positive imaginary part. By the symbols $\etp{\tau}$ and $q$, we mean $\rme^{2\uppi\rmi \tau}$. Set $q^n=\etp{n\tau}$ for $n \in \numQ$. For certain groups $G$, the symbol $\widetilde{G}$ denotes its double cover, and $\overline{G}$ denotes its image under a projection. Section \ref{sec:Modular forms} will explain these in more detail. The notation $G\alpha G$ denotes a double coset, see the paragraph before Lemma \ref{lemm:charCoset}. By $G_\tau$, we mean the isotropy group, or the stablizers of $\tau$, under the action of $G$ on $\uhp$, or some other set. By $\bigsqcup_jA_j$, we mean the disjoint union of the sets $A_j$. Some specific groups, such as $\pslZ$, $\slZ$, $\glpQ$, $\glpR$ and $\Gz{N}$, are defined in Section \ref{sec:Modular forms}. The symbol $\legendre{d}{c}$ denotes the Kronecker-Jacobi symbol, which we define in the paragraph after \eqref{eq:etaChar}.
The symbol $f\vert_kA$ and $v_1\vert_\alpha v_2$ are defined in \eqref{eq:slashAction} and Lemma \ref{lemm:charCoset} respectively.

A final remark. Some results are obtained with the aid of a computer algebra system---SageMath \cite{Sage}. For instance, all tables are generated, and many computations are carried out, by SageMath programs. The Sage code can be obtained from the repository \cite{Zhu23}. We will explain the usage of the code in \ref{sec:Usage of code}. In particular, we will explain how to use them to verify most theorems, lemmas, propositions, identities, how to produce all tables and how to produce identities automatically.

\section{Modular forms and Eta-quotients}
\label{sec:Modular forms}
We recall some basic notions and facts on modular forms and eta-quotients in this section.

The \textit{full modular group}, $\slZ$, is the group of $2\times2$ integral matrices of determinant one, which is a subgroup of $\glpR$, the group of $2\times2$ real matrices with positive determinants. Let $N$ be a positive integer, then $\Gz{N}$ denotes the matrices in $\slZ$ whose left-bottom entries are divisible by $N$. These groups serve as tools to state modular transformation laws and the action of Hecke operators. To deal with half-integral weights, we need the so-called metaplectic cover of $\glpR$, denoted by $\glptR$, consisting of ordered pairs of the form $(A, \varepsilon)$, where $A \in \glpR$, and $\varepsilon=\pm 1$. The composition law of $\glptR$ is given by
\begin{equation}
\label{eq:compositionCover}
\left(\tbtmat{a_1}{b_1}{c_1}{d_1}, \varepsilon_1\right)\left(\tbtmat{a_2}{b_2}{c_2}{d_2}, \varepsilon_2\right)=\left(\tbtmat{a_1}{b_1}{c_1}{d_1}\tbtmat{a_2}{b_2}{c_2}{d_2}, \varepsilon_1\varepsilon_2\delta\right),
\end{equation}
where $\delta \in \{\pm 1\}$ is determined by
\begin{equation}
\delta=\frac{\sqrt{c_1(a_2\tau+b_2)/(c_2\tau+d_2)+d_1}\sqrt{c_2\tau+d_2}}{\sqrt{c_1(a_2\tau+b_2)+d_1(c_2\tau+d_2)}}
\end{equation}
for any $\tau$ in the upper half plane. For the square root involved, we choose the principal branch, i.e., $-\uppi/2<\arg\sqrt{z}\leq\uppi/2$. The quantity $\delta$ is in fact a cocycle, only depending on the matrices $A=\tbtmat{a_1}{b_1}{c_1}{d_1}$, $B=\tbtmat{a_2}{b_2}{c_2}{d_2}$ and is usually denoted by $\sigma(A,B)$. See \cite[Section 4]{Str13} for details on the metaplectic group and cocycles. If $G$ if a subgroup of $\glpR$, then $\widetilde{G}$ denotes the preimage of $G$ under the natural projection that sends $(A, \varepsilon)$ to $A$. If $A \in \glpR$, then $\widetilde{A}$ denotes $(A, 1)$. Three specific matrices, $\tbtmat{1}{1}{0}{1}$, $\tbtmat{0}{-1}{1}{0}$, and $\tbtmat{1}{0}{0}{1}$, are denoted by $T$, $S$, and $I$ respectively. It is well-known that $T$ and $S$ generates $\slZ$. Similarly, $\widetilde{T}$ and $\widetilde{S}$ generates $\sltZ$. Note that $\widetilde{A}\cdot\widetilde{B}=\widetilde{AB}\cdot(I,\sigma(A,B))$, which is a direct consequence of \eqref{eq:compositionCover}. The groups of symmetries that we consider are all subgroups of $\slZ$. Therefore only matrices in $\glpQ$, the group of $2\times2$ rational matrices with positive determinants, are actually used, since non-rational matrices would lead to non-commensurability problems.

The modular group $\slZ$, and more generally, the group $\glptR$ acts from the left on the upper half plane $\uhp=\{\tau \in \numC: \Im \tau > 0\}$ by $\left(\tbtmat{a}{b}{c}{d},\varepsilon\right)\tau = \frac{a\tau+b}{c\tau+d}$. So the group should acts on functions from the right, which is given by the so-called \textit{slash operators}. Fix an integer or half integer $k \in \halfint$, we define
\begin{equation}
\label{eq:slashAction}
f \vert_k \left(\tbtmat{a}{b}{c}{d}, \varepsilon\right)(\tau)=\varepsilon^{-2k}(ad-bc)^{k/2}(c\tau+d)^{-k}f\left(\frac{a\tau+b}{c\tau+d}\right).
\end{equation}
The group law of $\glptR$ ensures that this is really a right group action. (If $k$ is not integer, one can similarly define the ``action'' of $\glpR$ on functions, but this does not obey the axioms of group actions.) We restrict the functions considered to meromorphic or continuous ones from upper half plane $\uhp$ to complex numbers $\numC$. Note that, if $A \in \glpR$, then $f\vert_k A$ means $f\vert_k \widetilde{A}$.

Let $f$ be a holomorphic function on $\uhp$, let $N$ be a positive integer and $k$ be an integer or half-integer, and let $v$ be a linear character of $\Gzt{N}$, that is to say, $v\colon \Gzt{N} \rightarrow \numC^\times$ is a group homomorphism. We say that $f$ is a \textit{weakly holomorphic modular form} of weight $k$ for the group $\Gzt{N}$ with multiplier system (or character) $v$, if the following conditions hold:
\begin{enumerate}
\item $f \vert_k \gamma = v(\gamma)f$ for each $\gamma \in \Gzt{N}$ (The modular transformation law).
\item For any $\gamma \in \sltZ$, the function $f\vert_k \gamma$ has a Fourier expansion, i.e., $f\vert_k \gamma(\tau)=\sum_{n}a_n q^n$, where the summation is over a lower-bounded set of rational numbers with bounded denominators, and $a_n \in \numC$. Those series are required to converge at any $\tau \in \uhp$.
\end{enumerate}
If additionally, $f\vert_k \gamma(\tau)=\sum_{m \geq 0}a_m q^m$ for any $\gamma \in \sltZ$, then we say $f$ is a \textit{modular form}. If $f\vert_k \gamma(\tau)=\sum_{m > 0}a_m q^m$ for any $\gamma \in \sltZ$, we say it is a \textit{cusp form}. Note that $m$ may be rational. The weakly holomorphic modular forms, modular forms, cusp forms of weight $k$ for the group $\Gzt{N}$ with multiplier system $v$, are denoted by $\MFormW{k}{N}{v}$, $\MForm{k}{N}{v}$, and $\CForm{k}{N}{v}$, respectively. Moreover, if $f$ is just meromorphic on $\uhp$, and satisfies Condition 1, we say $f$ transforms like a modular form. If in addition, for any $\gamma \in \sltZ$, $f\vert_k \gamma$ has a Fourier expansion with finitely many terms of negative power of $q$ which converges for $\Im \tau$ greater than some non-negative number, then we say $f$ is a \textit{meromorphic modular form}. The integer $N$ is called the \textit{level}. One can also define such kinds of forms for arbitrary subgroup of $\sltZ$.

For instance, the elliptic invariant $j=E_4^3/\Delta$ and the generating function of the partition function $\eta^{-1}$ are weakly holomorphic modular forms of weight $0$ and $-1/2$ respectively. On the other hand, the function $E_6^{-1}$ is a meromorphic modular form of weight $-6$, which has simple poles precisely at points equivalent to $\tau=\rmi$ modulo $\slZ$.

We now state the \textit{valence formula}, an essential tool for studying meromorphic modular forms. For this purpose, we require some more notation. Since this formula has a geometric nature, we should use the projective special linear group $\pslZ$, defined as the quotient group of $\slZ$ by $\{\pm I\}$, instead of $\slZ$ itself, or its metaplectic cover. For a subgroup $G$ in $\slZ$, $\overline{G}$ denotes the image of $G$ under the natural projection of $\slZ$ onto $\pslZ$. Let $f$ be a meromorphic function on $\uhp$. The notation $\ord_\tau(f)$ denotes the exponent of the first non-vanishing term of the Laurent expansion around $\tau$, as usual in complex analysis. If $f$ is a meromorphic modular form of weight $k$, and $A \in \slZ$, then $\ord_{\rmi\infty}(f\vert_k A)$ denotes\footnote{This definition may differ from some authors', who define the order of a form $f$ at a cusp as the order defined in our paper times the width of the cusp.} the exponent of the first non-vanishing term of the Fourier expansion of $f\vert_k A$. For any group $H$ and some subgroup $G$ in $H$, the cardinality of the orbit space of $G$ acting on $H$ by left translations is denoted by $[H\colon G]$ or $\abs{G\backslash H}$. Suppose that $H$ acts on a set $X$ from the left or right. The stablizer of an element $x \in X$ under $H$ is denoted by $H_x$, which is a subgroup of $H$. The orbit space $\{Hx\colon x \in X\}$, or $\{xH\colon x \in X\}$ is denoted by $H\backslash X$ or $X/H$, according to whether $H$ acts from the left or the right.
\begin{thm}
\label{thm:valence}
Let $G$ be a finite index subgroup of $\slZ$, and $v\colon \widetilde{G} \rightarrow \numC^{\times}$ a character. Let $k$ be an integer or half-integer. Suppose $f$ is a nonzero meromorphic modular form of weight $k$ for the group $\widetilde{G}$ with multiplier $v$. Then we have
\begin{equation}
\label{eq:valence}
\sum_{\overline{\gamma} \in \overline{G}\backslash\pslZ}\ord\nolimits_{\rmi\infty}(f\vert_k \gamma) + \sum_{\tau \in \overline{G}\backslash \uhp}\frac{1}{\abs{\overline{G}_\tau}}\ord\nolimits_\tau(f)=\frac{1}{12}[\pslZ : \overline{G}]k
\end{equation}
\end{thm}
This is a classical, well-known formula. However, we prefer to give a sketch of proof, since we have not seen this particular form in literature.
\begin{proof}
Set $m=[\slZ : G]$. Partition the group $\slZ$ into disjoint union of orbits: $\slZ=\bigsqcup_{1 \leq i \leq m} G\cdot M_i$. It can be proved that $\sltZ=\bigsqcup_{1 \leq i \leq m} \widetilde{G}\cdot \widetilde{M_i}$, still a disjoint union. Set $g=\prod_{i}f\vert_k(\widetilde{M_i})$. Then $g$ is a nonzero meromorphic modular form of weight $mk$ for the group $\sltZ$ with some character on $\sltZ$. Since the linear characters of $\sltZ$ form a cyclic group of order $24$ (see, for example, \cite{Boy13}),  the function $g^{24}$ is a nonzero meromorphic modular form of weight $24mk$ for  the group $\sltZ$ with trivial character. Applying the usual valence formula (see, for example,  \cite[Chap. VII, Sec. 3]{Ser73}) to $g^{24}$ gives us
\begin{equation}
\label{eq:valenceProof}
\sum_{1\leq i \leq m}\ord\nolimits_{\rmi\infty}(f\vert_k \widetilde{M_i}) + \sum_{\tau \in \sltZ\backslash \uhp}\frac{1}{e_\tau}\sum_{1\leq i\leq m}\ord\nolimits_\tau(f\vert_k\widetilde{M_i})=\frac{1}{12}mk,
\end{equation}
where $e_\tau=3$, if $\tau$ is equivalent to $\frac{-1+\sqrt 3\rmi}{2}$, $e_\tau=2$, if $\tau$ is equivalent to $\rmi$ modulo $\sltZ$, and $e_\tau=1$ otherwise. Using $\ord_\tau(f\vert_k \widetilde{M_i})=\ord_{\widetilde{M_i}\tau}(f)$, and $e_\tau=e_{\widetilde{M_i}\tau}$, we obtain
\begin{equation*}
\sum_{\tau \in \sltZ\backslash \uhp}\frac{1}{e_\tau}\sum_{1\leq i\leq m}\ord\nolimits_\tau(f\vert_k\widetilde{M_i}) = \sum_\tau\sum_i\frac{\ord_{\widetilde{M_i}\tau}(f)}{e_{\widetilde{M_i}\tau}},
\end{equation*}
which equals
\begin{equation*}
\sum_{\tau \in \sltZ\backslash \uhp}\sum_{\tau_1 \in \widetilde{G}\backslash\sltZ\tau}\frac{\ord_{\tau_1}(f)[\sltZ_{\tau_1}: \widetilde{G}_{\tau_1}]}{e_{\tau_1}}=\sum_{\tau \in \widetilde{G}\backslash\uhp}\frac{[\sltZ_{\tau}: \widetilde{G}_{\tau}]}{e_\tau}\ord\nolimits_\tau(f).
\end{equation*}
Combining this and \eqref{eq:valenceProof} gives
\begin{equation}
\label{eq:valenceProof2}
\sum_{\gamma \in \widetilde{G}\backslash\sltZ}\ord\nolimits_{\rmi\infty}(f\vert_k \gamma) + \sum_{\tau \in \widetilde{G}\backslash \uhp}\frac{1}{e_\tau}[\sltZ_{\tau}: \widetilde{G}_{\tau}]\ord\nolimits_\tau(f)=\frac{1}{12}[\sltZ : \widetilde{G}]k.
\end{equation}
Note that if $-I \in G$, then $[\sltZ : \widetilde{G}]=[\pslZ : \overline{G}]$, $[\sltZ_\tau : \widetilde{G}_\tau]=[\pslZ_\tau : \overline{G}_\tau]$, and $e_\tau=\lvert\pslZ_\tau\rvert$. Combining these identities and \eqref{eq:valenceProof2} gives \eqref{eq:valence}. On the other hand, if $-I \notin G$, then $[\sltZ : \widetilde{G}]=2[\pslZ : \overline{G}]$, $[\sltZ_\tau : \widetilde{G}_\tau]=2[\pslZ_\tau : \overline{G}_\tau]$, and $e_\tau=\lvert\pslZ_\tau\rvert$. Again by \eqref{eq:valenceProof2} we deduce \eqref{eq:valence}, since
\begin{equation*}
\sum_{\gamma \in \widetilde{G}\backslash\sltZ}\ord\nolimits_{\rmi\infty}(f\vert_k \gamma) = 2\sum_{\overline{\gamma} \in \overline{G}\backslash\pslZ}\ord\nolimits_{\rmi\infty}(f\vert_k \gamma).
\end{equation*}
\end{proof}

Let $\projQ$ be the usual one-dimensional projective space over $\numQ$, whose elements can be represented by a rational number $q/p$, or $\rmi\infty=1/0$. The level $N$ modular group $\Gz{N}$ (also $\Gzt{N}$ and $\Gzp{N}$) acts on $\projQ$ on the left by $\tbtmat{a}{b}{c}{d}\frac{q}{p}=\frac{aq+b}{cp+d}$. An orbit of this action is called a \textit{cusp}. By abuse of language, we also call a representative in an orbit a cusp. Let $\overline{s}$ be a cusp with representative $s$. The \textit{width} of $\overline{s}$ is defined as $w_{\overline{s}} = [\pslZ_s : \Gzp{N}_s]$. If $f$ is a non-zero meromorphic modular form of weight $k \in \halfint$, then define $\ord_{\overline{s}}(f)=\ord_s(f)=\ord_{\rmi\infty}(f\vert_k\gamma)$, where $\gamma \in \sltZ$ is any element such that $\gamma(\rmi\infty)=s$. All these considerations remain true if we replace $\Gz{N}$ with arbitrary subgroup $G$ of $\slZ$ of finite index. Using the notion of cusps, \eqref{eq:valence} can be reformulated as
\begin{equation}
\label{eq:valenceCusp}
\sum_{\overline{s} \in \overline{G}\backslash\projQ}w_{\overline{s}}\ord\nolimits_{\overline{s}}(f) + \sum_{\tau \in \overline{G}\backslash \uhp}\frac{1}{\abs{\overline{G}_\tau}}\ord\nolimits_\tau(f)=\frac{1}{12}[\pslZ : \overline{G}]k.
\end{equation}

A typical example of a modular form of general character (not necessarily induced by some Dirichlet character) is the \textit{Dedekind eta function}, defined by
\begin{equation}
\label{eq:defEta}
\eta(\tau)=q^{1/24}\prod_{n \geq 1}(1-q^n).
\end{equation}
This is a modular form of weight $1/2$ for the full modular group, with a character which is denoted by $v_{\eta}$. An explicit formula for $v_\eta$ is given by
\begin{equation}
\label{eq:etaChar}
v_\eta\left(\tbtmat{a}{b}{c}{d},\varepsilon\right)=\begin{cases}
\varepsilon\cdot\legendre{d}{\abs{c}}\etp{\frac{1}{24}\left((a+d-3)c-bd(c^2-1)\right)}   & \text{if }2 \nmid c, \\
\varepsilon\cdot\legendre{c}{d}\etp{\frac{1}{24}\left((a-2d)c-bd(c^2-1)+3d-3\right)}   & \text{if }2 \mid c.
\end{cases}
\end{equation}
The symbol $\legendre{\cdot}{\cdot}$ denotes the Kronecker-Jacobi symbol, an extension of Jacobi symbol, defined as follows:
\begin{enumerate}
\item $\legendre{m}{p}$ is the usual Legendre symbol if $p$ is an odd prime.
\item $\legendre{m}{2}$ equals $0$ if $2 \mid m$, and equals $(-1)^{(m^2-1)/8}$ if $2 \nmid m$.
\item $\legendre{m}{-1}$ equals $1$ if $m \geq 0$, and equals $-1$ otherwise.
\item $\legendre{m}{1}=1$ by convention.
\item $\legendre{m}{n}$ is defined to make it a complete multiplicative function of $n \in \numZ-\{0\}$.
\item $\legendre{m}{0}=0$ if $m \neq \pm1$, and $\legendre{\pm1}{0}=1$.
\end{enumerate}
The proof of \eqref{eq:etaChar}  is quite tedious. It is given in a number of places. See, for instance, \cite{Kno70}.

An \textit{eta-quotient}, which is a central object in this paper, is a function of the form $\prod_{n \in \numgeq{Z}{1}}\eta^{r_n}(n\tau)$, where $r_n$'s are integers and only finitely many of them are non-zero. Any eta-quotient can be written as $\eta^{\mathbf{r}}(\tau)=\prod_{n \mid N}\eta^{r_n}(n\tau)$ for some positive integer $N$ and $\mathbf{r}=(r_1, r_2, \dots)$. Define $v_\mathbf{r}$ as follows:
\begin{equation}
\label{eq:etaQuoChar}
v_\mathbf{r}\colon \Gzt{N} \rightarrow \numC^{\times} \quad \left(\tbtmat{a}{b}{c}{d}, \varepsilon\right) \mapsto \prod_{n \mid N}v_\eta^{r_n}\left(\tbtmat{a}{bn}{c/n}{d}, \varepsilon\right).
\end{equation}
Then $\eta^{\mathbf{r}}$ is a weakly holomorphic modular form of weight $\frac{1}{2}\sum{r_i}$ for the group $\Gzt{N}$ with character $v_\mathbf{r}$. If an eta-quotient is also a modular form, we call it a \textit{holomorphic eta-quotient}. We shall need the order of an eta-quotient at any cusp $d/c$, where $c$ and $d$ are coprime integers with $c>0$. The formula reads
\begin{equation}
\label{eq:etaOrder}
\ord\nolimits_{d/c}\left(\prod_{n \mid N}\eta^{r_n}(n\tau)\right)=\frac{1}{24}\sum_{n \mid N}\frac{\gcd(n,c)^2}{n}r_n.
\end{equation}
For a proof, see \cite[Section 5.9]{CoS17}.

By an \textit{eta-power}, we mean an eta-quotient $\eta^r(\tau)$ with $r \in \numZ$. Define the sequence $P_r(n)$ as
\begin{equation}
\label{eq:Prn}
\eta^r(\tau)=q^{r/24}\sum_{n \geq 0}P_r(n)q^n.
\end{equation}
For $n \in \numQ$ not a non-negative integer, set $P_r(n)=0$.

\section{Generalized Double Coset Operators}
\label{sec:Operators}
For modular forms with general characters, we present a new kind of double coset operators, under which the character of a form may change, but the level and weight remain unchanged.

Some group-theoretical notations. Let $H$ be an arbitrary group, let $G$ be a subgroup and $\alpha \in H$. Let $v$ be a linear complex character on $G$. Define $G^\alpha = \alpha^{-1}G\alpha$, and $v^\alpha$ to be the character on $G^\alpha$ that maps $\alpha^{-1}g\alpha$ to $v(g)$. The symbol $G\alpha G$ denotes the double coset $\{g_1\alpha g_2 \colon g_1, g_2 \in G\}$.
\begin{lemm}
\label{lemm:charCoset}
Suppose $v_1$ and $v_2$ are two linear complex characters on $G$, satisfying that $v_1^\alpha=v_2$ on $G\cap G^\alpha$. Then the map
\begin{equation}
G\alpha G \rightarrow \numC^{\times} \qquad g_1\alpha g_2 \mapsto v_1(g_1)v_2(g_2)
\end{equation}
is well-defined. This map is denoted by $v_1\vert_\alpha v_2$, and we say $v_1$ and $v_2$ are $\alpha$-compatible in this case.
\end{lemm}
\begin{proof}
Suppose that $g_1\alpha g_2=g_1'\alpha g_2'$ with $g_1, g_2, g_1', g_2' \in G$. The aim is to show $v_1(g_1)v_2(g_2)=v_1(g_1')v_2(g_2')$. We have
\begin{equation*}
\alpha^{-1}g_1'^{-1}g_1\alpha=g_2'g_2^{-1} \in G \cap G^\alpha.
\end{equation*}
It follows from $v_1^\alpha=v_2$ on $G\cap G^\alpha$ that
\begin{equation*}
v_2(g_2'g_2^{-1})=v_2(\alpha^{-1}g_1'^{-1}g_1\alpha)=v_1^\alpha(\alpha^{-1}g_1'^{-1}g_1\alpha)=v_1(g_1'^{-1}g_1).
\end{equation*}
So $v_2(g_2')v_2(g_2)^{-1}=v_1(g_1')^{-1}v_1(g_1)$, and the required equation follows.
\end{proof}

It is appropriate to introduce the generalized double coset operators acting on modular forms at this place.
\begin{deff}
\label{deff:doubleCosetOperators}
Let $G$ be a subgroup of $\slZ$, and $v_1, v_2\colon \widetilde{G} \rightarrow \numC^{\times}$ are linear characters. Let $\alpha \in \glptQ$ such that $v_1$ and $v_2$ are $\alpha$-compatible. Let $k \in \halfint$. We define a linear operator denoted by $T_\alpha$ as follows: For any meromorphic modular form $f$ of weight $k$ for the group $\widetilde{G}$ with character $v_1$, set
\begin{equation}
\label{eq:doubleCosetOperators}
T_\alpha f=\sum_{h \in R} (v_1\vert_\alpha v_2)^{-1}(h)f\vert_k h,
\end{equation}
where $R$ is any set of complete representatives of the orbit space $\widetilde{G}\backslash \widetilde{G}\alpha \widetilde{G}$ (In another words, the decomposition $\widetilde{G}\alpha \widetilde{G}=\bigcup_{h \in R}\widetilde{G}h$ is a disjoint union).
\end{deff}
Note that the operator $T_\alpha$ depends not only on $\alpha$, but also on the group $\widetilde{G}$, the characters $v_1$ and $v_2$, and the weight $k$. We do not insert these information into $T_\alpha$ to avoid heavy notations. The reader should infer these from context.

\begin{rema}
Wohlfahrt has developed a theory of generalized Hecke operators in \cite{Woh57}. See also \cite[Section 1]{vL57}, in which some more operators have been defined to study the Fourier coefficients of eta-powers. These operators are of course similar to those given in Definition \ref{deff:doubleCosetOperators}. Our treatment has the advantage that we consider the double covers of modular groups to avoid multiplier systems that are not group characters. Moreover, we can give a practical and general criterion on whether an operator can be defined from one character to another character, as stated and proved in Theorem \ref{thm:charCompatibleEtaDirichlet}.
\end{rema}

We list some basic properties of these operators.
\begin{prop}
\label{prop:basicPropTalpha}
Use the notations and assumptions of Definition \ref{deff:doubleCosetOperators}.
\begin{enumerate}
\item $T_\alpha f$ is independent of the choice of $R$.
\item The sum in \eqref{eq:doubleCosetOperators} is a finite sum. Equivalently, the set $R$ is finite.
\item If we have a disjoint union $\widetilde{G}\alpha \widetilde{G}=\bigsqcup_{i}\widetilde{G}\alpha g_i$, then $T_\alpha f=\sum_i v_2^{-1}(g_i)f\vert_k \alpha g_i$.
\item The function $T_\alpha f$ is still a meromorphic modular form of the same weight $k$, for the same group $\widetilde{G}$, but with the possibly different character $v_2$.
\item The operator $T_\alpha$ preserves regularity conditions, in other words, maps $\MFormWG{k}{G}{v_1}$, $\MFormG{k}{G}{v_1}$, and $\CFormG{k}{G}{v_1}$ to $\MFormWG{k}{G}{v_2}$, $\MFormG{k}{G}{v_2}$, and $\CFormG{k}{G}{v_2}$, respectively.
\end{enumerate}
\end{prop}
\begin{proof}
Part 3 is an immediate consequence of the definition. Any set of representatives can be written in this form.

To prove Part 1, suppose that $R'=\{t_i'\alpha g_i'\}$ is another set of representatives. Let $R=\{t_i\alpha g_i\}$.Then we have $\widetilde{G}\alpha \widetilde{G}=\bigsqcup_{i}\widetilde{G}\alpha g_i = \bigsqcup_{i}\widetilde{G}\alpha g_i'$. Without loss of generality, assume $\widetilde{G}\alpha g_i=\widetilde{G}\alpha g_i'$ for any $i$. It follows that there exists a unique $h_i \in \widetilde{G}$ such that $t_i\alpha g_i=h_it_i'\alpha g_i'$, for any $i$. Thus, we have
\begin{align*}
\sum_i (v_1\vert_\alpha v_2)^{-1}(t_i\alpha g_i)f\vert_k t_i\alpha g_i &= \sum_i (v_1\vert_\alpha v_2)^{-1}(h_it_i'\alpha g_i')f\vert_k h_it_i'\alpha g_i'\\
 &= \sum_i v_2(g_i')^{-1}v_1(h_it_i')^{-1}(f\vert_k h_it_i')\vert_k \alpha g_i' \\
  &= \sum_i (v_1\vert_\alpha v_2)^{-1}(t_i'\alpha g_i')f\vert_k t_i'\alpha g_i'.
\end{align*}
This means the value of $T_\alpha f$ does not depend on $R$.

Part 2 follows from the fact that $R_0$ is a set of complete representatives of $(\widetilde{G}\cap \widetilde{G}^\alpha)\backslash\widetilde{G}$ if and only if $\alpha R_0=\{\alpha h\colon h\in R_0\}$ is a set of  complete representatives of $\widetilde{G}\backslash\widetilde{G}\alpha\widetilde{G}$, and that $[\widetilde{G} : \widetilde{G}\cap \widetilde{G}^\alpha] < +\infty$, since $\alpha \in \glptQ$. For a detailed proof of the former assertion, see \cite[Lemma 2.7.1]{Mi06}, while the latter assertion says nothing but that elements of $\glptQ$ are all in the commensurator of $\widetilde{G}$, for which the reader may consult \cite[Section 10.5]{CoS17}.

We now proceed to prove Part 4. First of all, $T_\alpha f$ is obvious meromorphic on $\uhp$. Secondly, we shall prove it satisfies the modular transformation laws. Let $\gamma \in \widetilde{G}$ be arbitrary. For any $h \in R$, there is a unique $r_h \in R$ and $g_h \in \widetilde{G}$ such that $h\gamma=g_h r_h$. Since $\{h\gamma\colon h\in R\}$ is again a complete set of representatives of $\widetilde{G}\backslash \widetilde{G}\alpha \widetilde{G}$, the map $h \mapsto r_h$ is a bijection on $R$. Hence
\begin{align*}
(T_\alpha f)\vert_k \gamma &= \sum_{h \in R} (v_1\vert_\alpha v_2)^{-1}(h)f\vert_k h\gamma \\
&= \sum_{h \in R} (v_1\vert_\alpha v_2)^{-1}(g_h r_h \gamma^{-1})f\vert_k g_h r_h \\
&= \sum_{h \in R} (v_1\vert_\alpha v_2)^{-1}(g_h r_h \gamma^{-1})v_1(g_h)f\vert_k r_h \\
&= v_2(\gamma)\sum_{h \in R} (v_1\vert_\alpha v_2)^{-1}(r_h)f\vert_k r_h \\
&= v_2(\gamma)T_\alpha f.
\end{align*}
Finally, we shall prove $T_\alpha f$ is meromorphic at any cusp. For this, let $\gamma \in \sltZ$ be arbitrary. It suffices to prove that the Fourier expansion of $f\vert_k h\gamma$ has finitely many terms of negative power of $q$, for any $h \in R$. This can be shown by writing
\begin{equation}
\label{eq:decomphgamma}
h\gamma = (A, 1)(\tbtmat{a}{b}{0}{d}, 1)(I, \varepsilon)
\end{equation}
for some $A \in \slZ$, $a, b, d \in \numQ$ and $\varepsilon=\pm 1$ (automatically $a/d > 0$), and using the fact that the Fourier expansion of $f\vert_k(A, 1)$ has finitely many terms of negative power of $q$.

Now proceed to prove Part 5. It is obviously that $T_\alpha$ maps $\MFormWG{k}{G}{v_1}$ to $\MFormWG{k}{G}{v_2}$, since Part 4 has been established. To prove $T_\alpha$ maps $\MFormG{k}{G}{v_1}$ to $\MFormG{k}{G}{v_2}$, use the method in the proof of Part 4, i.e., \eqref{eq:decomphgamma} and the Fourier expansion of $f\vert_k(A, 1)$. A similar argument shows $T_\alpha$ maps $\CFormG{k}{G}{v_1}$ to $\CFormG{k}{G}{v_2}$.
\end{proof}

We are mainly interested in the case $G=\Gz{N}$. To compute the action of operators in Definition \ref{deff:doubleCosetOperators} explicitly, the following classical formulas are useful:
\begin{lemm}
\label{lemm:matrixSetDecomp}
Let $N$ and $l$ be positive integers. Let $\widetilde{M}_2(\numZ; l, N)$ be the set consisting of pairs $\left(\tbtmat{a}{b}{c}{d}, \varepsilon\right)$, where $\tbtmat{a}{b}{c}{d}$ is an integral matrix of determinant $l$ such that $\gcd (N, a)=1$ and $N \mid c$, and $\varepsilon=\pm 1$. Denoted by $\widetilde{M}^\star_2(\numZ; l, N)$ the set of $(A, \varepsilon) \in \widetilde{M}_2(\numZ; l, N)$ where the greatest common divisor of the entries of $A$ is $1$. Let $m>0$ satisfying $m^2 \mid l$ and $\gcd (N, m)=1$. Then we have
\begin{gather}
\widetilde{M}_2(\numZ; l, N)=\bigsqcup_{\twoscript{a,\,d>0,\, 0\leq b < d}{\gcd(N,a)=1,\, ad=l}}\Gzt{N}\widetilde{\tbtmat{a}{b}{0}{d}}, \label{eq:gamma0Ndecomp}\\
\widetilde{M}_2(\numZ; l, N)=\bigsqcup_{\twoscript{a>0,\,a^2 \mid l}{\gcd (N,a)=1}}\Gzt{N}\widetilde{\tbtmat{a}{0}{0}{l/a}}\Gzt{N}, \label{eq:gamma0NdecompDouble}\\
\widetilde{M}^\star_2(\numZ; l, N)=\Gzt{N}\widetilde{\tbtmat{1}{0}{0}{l}}\Gzt{N}\label{eq:gamma0NdecompStar}, \\
\Gzt{N}\widetilde{\tbtmat{m}{0}{0}{l/m}}\Gzt{N}=\bigsqcup_{\threescript{a,\,d>0,\, 0\leq b < d}{\gcd(N,a)=1,\, ad=l}{\gcd(a, b, d)=m}}\Gzt{N}\widetilde{\tbtmat{a}{b}{0}{d}}. \label{eq:doubleDecompLeft}
\end{gather}
\end{lemm}
\begin{proof}
These are just the lifts of the representatives from $\glpQ$ which are well known (see e.g. \cite[Section 6.5]{CoS17}).
\end{proof}

Let $G=\Gz{N}$ in Definition \ref{deff:doubleCosetOperators}. For any $\alpha \in \glptQ$, there exists an $\alpha' \in \widetilde{M}_2(\numZ; l, N')$ for some positive integer $N'$ and $l$, such that the operators $T_\alpha$ and $T_{\alpha'}$ are the same up to a constant factor. Among all operators $T_\alpha$ with $\alpha \in \glptQ$, we only studying those with $\alpha \in \widetilde{M}_2(\numZ; l, N)$. The reason for such restriction is that for these $T_\alpha$'s we have nice explicit descriptions, thanks to Lemma \ref{lemm:matrixSetDecomp}. In fact, we can further restrict ourselves to the operators $T_\alpha$ with $\alpha=\widetilde{\tbtmat{1}{0}{0}{l}}$, since by \eqref{eq:doubleDecompLeft} we have
\begin{align*}
\Gzt{N}\widetilde{\tbtmat{1}{0}{0}{l/m^2}}\Gzt{N}&=\bigsqcup_{\threescript{a,\,d>0,\, 0\leq b < d}{\gcd(N,a)=1,\, ad=l/m^2}{\gcd(a, b, d)=1}}\Gzt{N}\widetilde{\tbtmat{a}{b}{0}{d}},\\
\Gzt{N}\widetilde{\tbtmat{m}{0}{0}{l/m}}\Gzt{N}&=\bigsqcup_{\threescript{a',\,d'>0,\, 0\leq b' < d'}{\gcd(N,a'm)=1,\, a'd'=l/m^2}{\gcd(a', b', d')=1}}\Gzt{N}\widetilde{\tbtmat{a'm}{b'm}{0}{d'm}},
\end{align*}
and hence $T_\alpha$ with $\alpha=\widetilde{\tbtmat{m}{0}{0}{l/m}}$ is, up to a constant factor, the same as $T_\alpha$ with $\alpha=\widetilde{\tbtmat{1}{0}{0}{l/m^2}}$, and \eqref{eq:gamma0NdecompDouble} tells us $T_\alpha$ for any $\alpha \in \widetilde{M}_2(\numZ; l, N)$ is essentially the same as some $T_\alpha$ with $\alpha=\widetilde{\tbtmat{1}{0}{0}{l}}$ for some positive integer $l$. By abuse of language, let $T_l$ denote the operator $T_\alpha$ with $\alpha=\widetilde{\tbtmat{1}{0}{0}{l}}$.

We can express $T_l f$ by its Fourier series expansion. If both $v_1$ and $v_2$ are trivial, or induced by some Dirichlet character, the operator $T_l$ is just the usual Hecke operator, whose Fourier series expressions  are well-known. Otherwise, it is difficult to deduce a simple formula for general $v_1$ and $v_2$. We present the following formula under some special circumstances. For a positive integer $N$, the symbol $\rad(N)$ denotes the product of distinct prime factors of $N$.
\begin{prop}
\label{prop:radlcaseFourier}
Let $l$ and $N$ be positive integers such that $\rad(l) \mid \rad(N)$, and let $k \in \halfint$. Let $v_1$ and $v_2$ be complex linear characters on $\Gzt{N}$ which are $\widetilde{\tbtmat{1}{0}{0}{l}}$-compatible. Let $f$ be a meromorphic modular form of weight $k$, for the group $\Gzt{N}$, and with character $v_1$. Then 
\begin{equation}
\label{eq:radlcaseFourier1}
T_l f(\tau)=l^{-k/2}\sum_{0\leq b < l}v_2(\widetilde{T})^{-b}f\left(\frac{\tau+b}{l}\right).
\end{equation}
Moreover, if $f(\tau)=\sum_n a_n q^n$, where the summation range is a lower-bounded set of rationals with bounded denominators, then
\begin{equation}
\label{eq:radlcaseFourier2}
T_l f(\tau)=l^{1-k/2}\sum_{\etp{n/l}=v_2(\widetilde{T})}a_nq^{n/l}.
\end{equation}
\end{prop}
\begin{proof}
Use the representatives in \eqref{eq:doubleDecompLeft} with $m=1$. The assumption $\rad(l) \mid \rad(N)$ implies that $a$ can only be $1$ in this set of representatives. Hence
\begin{equation*}
T_l f(\tau)=\sum_{0 \leq b < l}(v_1\vert_\alpha v_2)^{-1}\widetilde{\tbtmat{1}{b}{0}{l}}f\vert_k\widetilde{\tbtmat{1}{b}{0}{l}}(\tau),
\end{equation*}
where $\alpha=\widetilde{\tbtmat{1}{0}{0}{l}}$. We have $v_1\vert_\alpha v_2\widetilde{\tbtmat{1}{b}{0}{l}}=v_2(\widetilde{T})^b$, and $f\vert_k\widetilde{\tbtmat{1}{b}{0}{l}}(\tau)=l^{-k/2}f\left(\frac{\tau+b}{l}\right)$, by Lemma \ref{lemm:charCoset} and \eqref{eq:slashAction} respectively. Then \eqref{eq:radlcaseFourier1} follows immediately. Finally, a straightforward calculation shows that
\begin{align*}
T_l f(\tau)&=l^{-k/2}\sum_n\sum_{0 \leq b < l} \left(\etp{n/l}v_2(\widetilde{T})^{-1}\right)^b a_nq^{n/l} \\
&=l^{-k/2}\left(\sum_{\etp{n/l}=v_2(\widetilde{T})}la_nq^{n/l}+\sum_{\etp{n/l}\neq v_2(\widetilde{T})}\frac{1-\etp{n}v_2(\widetilde{T})^{-l}}{1-\etp{n/l}v_2(\widetilde{T})^{-1}}a_nq^{n/l}\right).
\end{align*}
It remains to show that the second sum in the right-hand side vanishes. By compatibility, we have $v_1(\widetilde{T})=v_2(\widetilde{T})^l$, since $\widetilde{T}^l\in\Gzt{N}\cap\Gzt{N}^\alpha$ and $v_1^\alpha(\widetilde{T}^l)=v_2(\widetilde{T}^l)$. Applying the operator $\vert_k\widetilde{T}$ to $f$ we find that if $a_n\neq 0$ then $\etp{n}=v_1(\widetilde{T})$. Therefore
\begin{equation*}
1-\etp{n}v_2(\widetilde{T})^{-l}=1-\etp{n}v_1(\widetilde{T})^{-1}=0.
\end{equation*}
\end{proof}
\begin{rema}
If we have not the condition $\rad(l) \mid \rad(N)$, then by \eqref{eq:doubleCosetOperators} and \eqref{eq:doubleDecompLeft},
\begin{multline}
\label{eq:Tlgeneral}
T_l f(\tau)=l^{-k/2}\sum_{0\leq b < l}v_2(\widetilde{T})^{-b}f\left(\frac{\tau+b}{l}\right) \\
+l^{-k/2}\cdot\sum_{\twoscript{1<a \mid l}{\gcd(N,a)=1}}a^k\sum_{\twoscript{0 \leq b < l/a}{gcd(a, b, l/a)=1}}\left(v_1\vert_{\alpha}v_2\right)^{-1}\widetilde{\tbtMat{a}{b}{0}{l/a}}f\left(\frac{a\tau+b}{l/a}\right),
\end{multline}
where $\alpha=\widetilde{\tbtmat{1}{0}{0}{l}}$. The only obstacle to giving a more explicit expression for arbitrary $v_1$, $v_2$, $l$ and $N$ is the factor $\left(v_1\vert_{\alpha}v_2\right)^{-1}\widetilde{\tbtmat{a}{b}{0}{l/a}}$. We now give an algorithm to calculate this quantity. Fix some $1<a \mid l$, and $0 \leq b < l/a$ such that $\gcd(a, b, l/a)=1$. Put $d=l/a$. We choose some $x,\,y,\,z \in \numZ$ as follows:
\begin{enumerate}
\item Choose an $x\in \numZ$ such that $\gcd(Nd, -Nb+ax)=1$. (This can be done since $\gcd(Nd, Nb, a)=1$.)
\item Then choose $y,\,z \in \numZ$ such that $(-Nb+ax)y + Ndz=1$.
\end{enumerate}
It follows by a straightforward calculation, where the cocycle conditions should be taken care of (see \cite[Theorem 4.1]{Str13}), that
\begin{equation}
\label{eq:matrixab0dDecom}
\widetilde{\tbtMat{a}{b}{0}{d}}=\widetilde{\tbtMat{-Nb+ax}{z}{-Nd}{y}}\widetilde{\tbtMat{1}{0}{0}{l}}\widetilde{\tbtMat{ay}{by-dz}{N}{x}}.
\end{equation}
Consequently,
\begin{equation}
\label{eq:v1v2charGeneral}
\left(v_1\vert_{\alpha}v_2\right)^{-1}\widetilde{\tbtMat{a}{b}{0}{d}}=v_1^{-1}\widetilde{\tbtMat{-Nb+ax}{z}{-Nd}{y}}v_2^{-1}\widetilde{\tbtMat{ay}{by-dz}{N}{x}},
\end{equation}
since $\widetilde{\tbtmat{-Nb+ax}{z}{-Nd}{y}}$ and $\widetilde{\tbtmat{ay}{by-dz}{N}{x}}$ are in $\Gzt{N}$.
\end{rema}

Also, under the same assumption $\rad(l) \mid \rad(N)$, we need an estimate of the lower-bound of the order of $T_l f$ at cusps, which plays an important role in the application in Section \ref{sec:Express}.
\begin{prop}
\label{prop:orderTl}
Use notations and assumptions in Proposition \ref{prop:radlcaseFourier}. Let $a,\, c$ be co-prime integers with $c>0$. Then
\begin{equation*}
\ord\nolimits_{\frac{a}{c}}(T_l f) \geq \min_{\lambda=0, 1,\dots, l-1}\frac{\gcd(a+c\lambda,cl)^2}{l}\ord\nolimits_{\frac{a+c\lambda}{cl}}(f).
\end{equation*}
\end{prop}
\begin{proof}
Let $b,\, d$ be any integers such that $ad-bc=1$. For $\lambda=0, 1, \dots, l-1$, put $g =\gcd(a+c\lambda, cl)$. Let $x,\, y$ be any integers such that $(a+c\lambda)g^{-1}x+clg^{-1}y=1$.  Then we have
\begin{equation*}
\widetilde{\begin{pmatrix}1 & \lambda \\ 0 & l\end{pmatrix}}\widetilde{\begin{pmatrix}a & b \\ c & d\end{pmatrix}}=\widetilde{\begin{pmatrix}(a+c\lambda)g^{-1} & -y \\ clg^{-1} & x\end{pmatrix}}\widetilde{\begin{pmatrix}g & (b+d\lambda)x+dly \\ 0 & lg^{-1}\end{pmatrix}},
\end{equation*}
where the four matrices occurring are all integral. If the dependency on $\lambda$ should be indicated, we write $x_\lambda$, $y_\lambda$, $g_\lambda$ instead of $x$, $y$, $g$. Using this identity and an equivalent form of \eqref{eq:radlcaseFourier1}, we obtain
\begin{align*}
\label{eq:proofOrder}
T_l f\vert_k\widetilde{\tbtmat{a}{b}{c}{d}}&=\sum_{0\leq \lambda < l}v_2(\widetilde{T})^{-\lambda}f\vert_k\widetilde{\tbtmat{1}{\lambda}{0}{l}}\widetilde{\tbtmat{a}{b}{c}{d}} \\
&=\sum_{0\leq \lambda < l}v_2(\widetilde{T})^{-\lambda}f\vert_k\widetilde{\tbtmat{(a+c\lambda)g_\lambda^{-1}}{-y_\lambda}{clg_\lambda^{-1}}{x_\lambda}}\widetilde{\tbtmat{g_\lambda}{(b+d\lambda)x_\lambda+dly_\lambda}{0}{lg_\lambda^{-1}}}.
\end{align*}
Since $\widetilde{\tbtmat{(a+c\lambda)g_\lambda^{-1}}{-y_\lambda}{clg_\lambda^{-1}}{x_\lambda}} \in \sltZ$, $\ord_{\frac{a}{c}}(T_l f)=\ord_{\rmi\infty}\left(T_l f\vert_k\widetilde{\tbtmat{a}{b}{c}{d}}\right)$, and 
\begin{equation*}
\ord\nolimits_{\frac{a+c\lambda}{cl}}(f)=\ord\nolimits_{\rmi\infty}\left( f\vert_k\widetilde{\tbtmat{(a+c\lambda)g_\lambda^{-1}}{-y_\lambda}{clg_\lambda^{-1}}{x_\lambda}}\right),
\end{equation*}
the required inequality follows.
\end{proof}

We now turn to concrete examples of operators acting on spaces with characters given by eta-quotients. The key point is to determine when two characters given by eta-quotients are $\widetilde{\tbtmat{1}{0}{0}{l}}$-compatible. The following fact allows us to reduce the problem on verifying formulas for elements in a double cover $\Gzt{N}$ to the one on verifying formulas for elements in a subset of the matrix group $\Gz{Nl}$.

\begin{lemm}
\label{lemm:charCompatible}
Let $l,\, N$ be positive integers. Let $v_1$ and $v_2$ be complex characters on $\Gzt{N}$. Then $v_1$ and $v_2$ are $\widetilde{\tbtmat{1}{0}{0}{l}}$-compatible if and only if $v_1\widetilde{\tbtmat{a}{b}{c}{d}}=v_2\widetilde{\tbtmat{a}{b\cdot l}{c/l}{d}}$ for any $\tbtmat{a}{b}{c}{d}=\tbtmat{1}{1}{0}{1},\, \tbtmat{-1}{0}{0}{-1}$, and matrices in $\Gz{Nl}$ with $c>0$.
\end{lemm}
\begin{proof}
Set $\alpha=\tbtmat{1}{0}{0}{l}$. By Lemma \ref{lemm:charCoset}, $v_1$ and $v_2$ are $\widetilde\alpha$-compatible if and only if $v_1^{\widetilde\alpha}=v_2$ on $\Gzt{N} \cap \Gzt{N}^{\widetilde\alpha}=\{(\tbtmat{a}{b}{c}{d}, \pm 1) \in \sltZ \colon N\mid c,\, l \mid b\}$. This is equivalent to saying that
\begin{equation}
\label{eq:charCompatible1}
v_1\widetilde{\tbtmat{a}{b/l}{c\cdot l}{d}}=v_2\widetilde{\tbtmat{a}{b}{c}{d}}
\end{equation}
for any $\widetilde{\tbtmat{a}{b}{c}{d}}$ with $\tbtmat{a}{b}{c}{d}\in \Gz{N} \cap \Gz{N}^\alpha$, since another lift of $\tbtmat{a}{b}{c}{d}$ satisfies
\begin{equation*}
(\tbtmat{a}{b}{c}{d}, -1)=\widetilde{\tbtmat{a}{b}{c}{d}}\cdot(I, -1)=\widetilde{\tbtmat{a}{b}{c}{d}}\cdot\widetilde{\tbtmat{-1}{0}{0}{-1}}^2,
\end{equation*}
and $-I \in \Gz{N} \cap \Gz{N}^\alpha$. By the group isomorphism $\Gz{N} \cap \Gz{N}^\alpha\rightarrow \Gz{Nl}$ that sends $\tbtmat{a}{b}{c}{d}$ to $\tbtmat{a}{b/l}{c\cdot l}{d}$, \eqref{eq:charCompatible1} is equivalent to
\begin{equation}
\label{eq:charCompatible2}
v_1\widetilde{\tbtmat{a}{b}{c}{d}}=v_2\widetilde{\tbtmat{a}{bl}{c/l}{d}}
\end{equation}
for any $\tbtmat{a}{b}{c}{d} \in \Gz{Nl}$. This now holds in the case $c=0, a<0$ by the assumption and a careful analysis of the involved cocycles. It also holds in the case $c<0$ since the cocycle $\sigma\left(\tbtmat{-a}{-b}{-c}{-d},-I\right)=\sigma\left(\tbtmat{-a}{-bl}{-c/l}{-d},-I\right)=-1$ and hence
\begin{align}
\notag v_1\widetilde{\tbtmat{a}{b}{c}{d}} &= v_1\left(\tbtmat{-a}{-b}{-c}{-d},1\right)v_1(-I,-1)\\
&= v_1\left(\tbtmat{-a}{-b}{-c}{-d},1\right)v_1(-I,1)^3\\
&= v_2\left(\tbtmat{-a}{-bl}{-c/l}{-d},1\right)v_2(-I,1)^3\\
\notag &= v_2\widetilde{\tbtmat{a}{bl}{c/l}{d}}.
\end{align}
In the remaining case $c=0,a>0$, \eqref{eq:charCompatible2} follows easily from the assumption.
\end{proof}

\begin{rema}
\label{rema:charCompatible}
For a specific group $\Gz{Nl}$, if we know a particular set $X$ of generators, then $v_1$ and $v_2$ are $\widetilde{\tbtmat{1}{0}{0}{l}}$-compatible if and only if $v_1\widetilde{\tbtmat{a}{b}{c}{d}}=v_2\widetilde{\tbtmat{a}{b\cdot l}{c/l}{d}}$ for any $\tbtmat{a}{b}{c}{d} \in X$, since $\{\widetilde{\gamma}\colon \gamma \in X\}$ is a set of generators of $\Gzt{Nl}$. To prove this, let $\left(\gamma,\varepsilon\right)\in\Gzt{Nl}$ be arbitrary. We need to prove $\left(\gamma,\varepsilon\right)$ belongs to the group generated by $\{\widetilde{\gamma}\colon \gamma \in X\}$. Since $X$ generates $\Gz{Nl}$, $\gamma$ can be written as a finite product $\prod_{i}\gamma_i$ of elements $\gamma_i$ in $X\cup X^{-1}$. Therefore, $\prod_{i}\widetilde{\gamma_i}=\left(\gamma,\varepsilon'\right)$ for some $\varepsilon'\in\{\pm 1\}$. If $\varepsilon'=\varepsilon$, then the assertion follows. Now assume that $\varepsilon'=-\varepsilon$. Write $-I=\prod_{j}g_j$ with $g_j\in X\cup X^{-1}$. Then $\prod_{j}\widetilde{g_j}=(-I,\varepsilon'')$ with $\varepsilon''\in\{\pm 1\}$, and consequently $\left(\prod_{j}\widetilde{g_j}\right)^2=(I,-1)$. Hence $\left(\gamma,\varepsilon\right)=\prod_{i}\widetilde{\gamma_i}\left(\prod_{j}\widetilde{g_j}\right)^2$, from which the desired assertion follows.
\end{rema}

We need some facts about the action of \emph{Fricke involutions} on eta-quotients. For a positive integer $N$ and an integral sequence $\mathbf{r}=(r_1,r_2,\dots)$ with the property $r_n\neq 0 \implies n \mid N$, let $F_N=\tbtmat{0}{-1}{N}{0}$, and let $F_N\mathbf{r}$ be the sequence whose $n$-th term is $r_{N/n}$ if $n \mid N$, and $0$ otherwise. Set $k=\sum r_n/2$.
\begin{lemm}
We have $\widetilde{F_N}^{-1}\Gzt{N}\widetilde{F_N}=\Gzt{N}$, and
\begin{equation*}
\eta^{\mathbf{r}}\vert_k\widetilde{F_N}=N^{-\frac{k}{2}}\etp{-\frac{k}{4}}\left(\prod_{n \mid N}n^{r_{N/n}/2}\right)\eta^{F_N\mathbf{r}}.
\end{equation*}
In particular, $v_{\mathbf{r}}^{\widetilde{F_N}}=v_{F_N\mathbf{r}}$ on $\Gzt{N}$.
\end{lemm}
\begin{proof}
The first identity is the lift of the well-known identity $F_N^{-1}\Gz{N}F_N=\Gz{N}$ and the last one follows from the second one. To prove the second one, note that
\begin{align*}
\eta^{\mathbf{r}}\vert_k\widetilde{F_N}(\tau)&=\prod_{n\mid N}\eta(n\tau)^{r_n}\vert_{r_n/2}\widetilde{F_N}\\
&=\prod_{n\mid N}(\sqrt{N}\tau)^{-r_n/2}\eta\left(-\frac{1}{(N/n)\tau}\right)^{r_n}\\
&=\prod_{n\mid N}(\sqrt{N}\tau)^{-r_n/2}\left(\frac{N}{n}\tau\right)^{r_n/2}\etp{-\frac{r_n}{8}}\eta\left(\frac{N}{n}\tau\right)^{r_n}\\
&=N^{-\frac{k}{2}}\etp{-\frac{k}{4}}\left(\prod_{n \mid N}n^{r_{N/n}/2}\right)\eta^{F_N\mathbf{r}}.
\end{align*}
\end{proof}

Based on Lemma \ref{lemm:charCompatible}, we are able to give a very nice criterion of $\widetilde{\tbtmat{1}{0}{0}{l}}$-compatibility. This criterion deals with group characters each of which is the character of an eta-quotient times that induced by some Dirichlet character. We say that a complex linear character $v$ on $\Gzt{N}$ is induced by a Dirichlet character $\chi$ modulo $N$, if $v\left(\tbtmat{a}{b}{c}{d}, \varepsilon\right)=\chi(d)$ for any $\tbtmat{a}{b}{c}{d}\in \Gz{N}$ and $\varepsilon = \pm 1$. By slight abuse of language, we also write $\chi\left(\tbtmat{a}{b}{c}{d}, \varepsilon\right)=\chi(d)$. Note that $\chi^{\widetilde{F_N}}=\chi^{-1}$.
\begin{lemm}
\label{lemm:adjointCompatible}
Let $l$ and $N$ be positive integers, and let $\mathbf{r}=(r_1, r_2,\dots)$, $\mathbf{r'}=(r_1',r_2',\dots)$ be integral sequences such that $r_n \neq 0$ or $r_n' \neq 0$ implies $n \mid N$ for any positive integer $n$. Let $\chi$ and $\chi'$ be Dirichlet characters modulo $N$. Then $\chi\cdot v_\mathbf{r}$ and $\chi' \cdot v_\mathbf{r'}$ are $\widetilde{\tbtmat{1}{0}{0}{l}}$-compatible if and only if $\chi'^{-1}\cdot v_{F_N\mathbf{r'}}$ and $\chi^{-1} \cdot v_{F_N\mathbf{r}}$ are $\widetilde{\tbtmat{1}{0}{0}{l}}$-compatible.
\end{lemm}
\begin{proof}
Set $\alpha=\widetilde{\tbtmat{1}{0}{0}{l}}$. The desired assertion can be proved by the following sequence of logical equivalences:
\begin{align*}
&\left(\chi\cdot v_\mathbf{r}\right)^\alpha=\chi' \cdot v_\mathbf{r'} \text{ on } \Gzt{N}\cap\Gzt{N}^\alpha\\
\Longleftrightarrow &\chi\cdot v_\mathbf{r}=\left(\chi'\cdot v_\mathbf{r'}\right)^{\alpha^{-1}} \text{ on } \Gzt{N}^{\alpha^{-1}}\cap\Gzt{N}\\
\Longleftrightarrow &\left(\chi\cdot v_\mathbf{r}\right)^{\widetilde{F_N}}=\left(\chi' \cdot v_\mathbf{r'}\right)^{\alpha^{-1}\widetilde{F_N}} \text{ on } \Gzt{N}^{\alpha^{-1}\widetilde{F_N}}\cap\Gzt{N}^{\widetilde{F_N}}\\
\Longleftrightarrow &\chi^{-1}\cdot v_{F_N\mathbf{r}}=\left(\chi' \cdot v_\mathbf{r'}\right)^{l^{-1}\widetilde{F_N}\alpha} \text{ on } \Gzt{N}^{l^{-1}\widetilde{F_N}\alpha}\cap\Gzt{N}\\
\Longleftrightarrow &\left(\chi'^{-1} \cdot v_{F_N\mathbf{r'}}\right)^{\alpha}=\chi^{-1}\cdot v_{F_N\mathbf{r}} \text{ on } \Gzt{N}\cap\Gzt{N}^{\alpha}.
\end{align*}
\end{proof}

Now we state the Main Theorem.

\begin{thm}
\label{thm:charCompatibleEtaDirichlet}
Let $l$ and $N$ be positive integers, and let $\mathbf{r}=(r_1, r_2,\dots)$, $\mathbf{r'}=(r_1',r_2',\dots)$ be integral sequences such that $r_n \neq 0$ or $r_n' \neq 0$ implies $n \mid N$ for any positive integer $n$, and that $\sum_{n \mid N}r_n=\sum_{n \mid N}r_n'$. Let $\chi$ and $\chi'$ be Dirichlet characters modulo $N$. Let $k=\sum r_n/2$ and $\delta=l^{2\abs{k}}\prod_{2 \nmid r_n-r_n'}n$. Then $\chi\cdot v_\mathbf{r}$ and $\chi' \cdot v_\mathbf{r'}$ are $\widetilde{\tbtmat{1}{0}{0}{l}}$-compatible if and only if the following four conditions hold:
\begin{gather} 
l\sum_{n \mid N}\frac{N}{n}r_n \equiv \sum_{n \mid N}\frac{N}{n}r_n' \pmod{24} \label{eq:charCom1}\\
\sum_{n \mid N}nr_n \equiv l\sum_{n \mid N}nr_n' \pmod{24} \label{eq:charCom2}\\
2 \nmid Nl \Rightarrow \delta \equiv 1 \bmod 4 \label{eq:charCom3}\\
\chi(d)\chi^{\prime-1}(d)=\legendre{\delta}{d} \text{ if } \gcd(Nl, d)=1,\, d \in \numZ_{\neq 0}. \label{eq:charCom4}
\end{gather}
In particular, $v_\mathbf{r}$ and $v_\mathbf{r'}$ are $\widetilde{\tbtmat{1}{0}{0}{l}}$-compatible if and only if \eqref{eq:charCom1} and \eqref{eq:charCom2} hold, and $\delta$ is a perfect square.
\end{thm}
\begin{proof}
Assume the four given conditions hold. By Lemma \ref{lemm:charCompatible} we shall prove
\begin{equation}
\label{eq:requireToProofCompatible}
\chi\cdot v_\mathbf{r}\widetilde{\tbtmat{a}{b}{c}{d}}=\chi'\cdot v_\mathbf{r'}\widetilde{\tbtmat{a}{b\cdot l}{c/l}{d}}
\end{equation}
for any $\tbtmat{a}{b}{c}{d}=\tbtmat{1}{1}{0}{1},\, \tbtmat{-1}{0}{0}{-1}$, and matrices in $\Gz{Nl}$ with $c>0$. For $\tbtmat{a}{b}{c}{d}=\tbtmat{1}{1}{0}{1}$, the required equation is equivalent to $\etp{\frac{1}{24}\sum_{n \mid N}nr_n}=\etp{\frac{1}{24}l\sum_{n \mid N}nr_n'}$, using \eqref{eq:etaChar} and \eqref{eq:etaQuoChar}. This is true, since we have \eqref{eq:charCom2} by assumption. For $\tbtmat{a}{b}{c}{d}=\tbtmat{-1}{0}{0}{-1}$, the required equation \eqref{eq:requireToProofCompatible} is equivalent to $\chi(-1)\etp{-\frac{1}{4}\sum_{n \mid N}r_n}=\chi'(-1)\etp{-\frac{1}{4}\sum_{n \mid N}r_n'}$, again using \eqref{eq:etaChar} and \eqref{eq:etaQuoChar}. This is also true, since we have  $\sum_{n \mid N}r_n=\sum_{n \mid N}r_n'$, and $\chi(-1)=\legendre{\delta}{-1}\chi'(-1)=\chi'(-1)$ by \eqref{eq:charCom4} with $d=-1$.

It remains to prove that \eqref{eq:requireToProofCompatible} holds for $\tbtmat{a}{b}{c}{d} \in \Gz{Nl}$ with $c>0$. As a prerequisite, we first evaluate the quotient $v_\eta^{r_n'}\widetilde{\tbtmat{a}{bnl}{c/(nl)}{d}}/ v_\eta^{r_n}\widetilde{\tbtmat{a}{bn}{c/n}{d}}$, which is denoted by $Q$ in the remaining. There are two cases: $2 \nmid \frac{c}{n}$, and $2 \mid \frac{c}{n}$. Assume that $2 \nmid \frac{c}{n}$. Hence $2 \nmid \frac{c}{nl}$. Using the first line of \eqref{eq:etaChar}, we have
\begin{multline}
\label{eq:Qeq1}
Q=\legendre{d}{c}^{r_n'-r_n}\legendre{d}{n}^{r_n'-r_n}\legendre{d}{l}^{r_n'} \\
\cdot\etp{\frac{1}{24}\left((a+d-3)\frac{c}{Nl}\frac{N}{n}(r_n'-lr_n)-bd\frac{c^2}{Nl}\frac{N}{n}(r_n'-lr_n)+bdn(lr_n'-r_n)\right)}.
\end{multline}
For the case $2 \mid \frac{c}{n}$, we split the evaluation into two subcases: $2 \mid \frac{c}{nl}$ and $2 \nmid \frac{c}{nl}$. If $2 \mid \frac{c}{nl}$, we use the second line of \eqref{eq:etaChar}, and if $2 \nmid \frac{c}{nl}$, we use the first line of \eqref{eq:etaChar} for the numerator of $Q$, and the second line for the denominator. In both subcases, we obtain
\begin{multline}
\label{eq:Qeq2}
Q=\legendre{c}{d}^{r_n'-r_n}\legendre{n}{d}^{r_n'-r_n}\legendre{l}{d}^{r_n'}\\
\cdot\etp{\frac{1}{24}\left((a-2d)\frac{c}{Nl}\frac{N}{n}(r_n'-lr_n)-bd\frac{c^2}{Nl}\frac{N}{n}(r_n'-lr_n)+bdn(lr_n'-r_n)\right)}\\
\cdot\etp{\frac{1}{8}(d-1)(r_n'-r_n)}.
\end{multline}
Multiplying \eqref{eq:Qeq1} or \eqref{eq:Qeq2} for all positive divisors $n$ of $N$ and using \eqref{eq:charCom1} and \eqref{eq:charCom2} give us
\begin{multline}
\label{eq:Qeq3}
v_{\mathbf{r'}}\widetilde{\tbtmat{a}{b\cdot l}{c/l}{d}}/v_{\mathbf{r}}\widetilde{\tbtmat{a}{b}{c}{d}}= \\
\prod_{2 \nmid \frac{c}{n}}\legendre{d}{c}^{r_n'-r_n}\legendre{d}{n}^{r_n'-r_n}\legendre{d}{l}^{r_n'}\etp{\frac{1}{8}(d-1)\frac{c}{Nl}\frac{N}{n}(r_n'-lr_n)} \\
\cdot\prod_{2 \mid \frac{c}{n}}\legendre{c}{d}^{r_n'-r_n}\legendre{n}{d}^{r_n'-r_n}\legendre{l}{d}^{r_n'}\etp{\frac{1}{8}(d-1)(r_n'-r_n)}.
\end{multline}
We shall prove that the right-hand side of the above equation equals $\legendre{\delta}{d}$ for $d \neq 0$, from which and \eqref{eq:charCom4} the required equation \eqref{eq:requireToProofCompatible} for $\tbtmat{a}{b}{c}{d} \in \Gz{Nl}$ with $c>0$ follows.

We prove by cases, to compute the right-hand side of \eqref{eq:Qeq3}, which is denoted by $R$ in the remaining. There are three cases:
\begin{equation}
\label{eq:Rcases}
2 \mid l; \qquad 2 \nmid l,\,2 \nmid c; \qquad 2 \nmid l,\, 2 \mid c.
\end{equation}
If $2 \mid l$, then there does not exist $n \mid N$ such that $2 \nmid \frac{c}{n}$. It follows that, using $\sum_n {r_n}=\sum_n {r_n'}$,
\begin{equation}
\label{eq:R2midl}
R=\prod_{n \mid N}\legendre{n}{d}^{r_n'-r_n}\legendre{l}{d}^{r_n'}=\legendre{\delta}{d},
\end{equation}
as required. If $2 \nmid l$, and $2 \nmid c$, then there exists no $n \mid N$ such that $2 \mid \frac{c}{n}$. In this case, using \eqref{eq:charCom1} and the assumption $\sum_n {r_n}=\sum_n {r_n'}$, we have
\begin{equation}
\label{eq:R2nmidl2nmidc}
R=\prod_{n \mid N}\legendre{d}{n}^{r_n'-r_n}\legendre{d}{l}^{r_n'}=\legendre{d}{\delta}.
\end{equation}
Only at this place, we use \eqref{eq:charCom3}. Since $2 \nmid l$, and $2 \nmid c$, we have $2 \nmid Nl$. Consequently $\delta \equiv 1 \bmod 4$. Hence we have $\legendre{d}{\delta}=\legendre{\delta}{d}$ for any $d \in \numZ$. It follows that $R=\legendre{\delta}{d}$, as required. The last case is $2 \nmid l$ and  $2 \mid c$. Note that in this case $2 \nmid d$. We shall use the following version of the Jacobi reciprocity law: if $m_1$ and $m_2$ are co-prime odd numbers, then
\begin{equation}
\legendre{m_1}{m_2}\legendre{m_2}{m_1}=\varepsilon(m_1, m_2)(-1)^{\frac{m_1-1}{2}\frac{m_2-1}{2}},
\end{equation}
where $\varepsilon(m_1, m_2)$ is equal to $-1$ if $m_1$ and $m_2$ are both negative, and to $1$ otherwise. Apply this law to $\legendre{d}{c/n}$ and $\legendre{d}{l}$, and we get
\begin{align*}
\prod_{2 \nmid \frac{c}{n}}&\legendre{d}{c}^{r_n'-r_n}\legendre{d}{n}^{r_n'-r_n}\legendre{d}{l}^{r_n'} = \prod_{2 \nmid \frac{c}{n}}\legendre{d}{c/n}^{r_n'-r_n}\legendre{d}{l}^{r_n'} \\
&= \prod_{2 \nmid \frac{c}{n}}\legendre{c/n}{d}^{r_n'-r_n}\legendre{l}{d}^{r_n'}(-1)^{\frac{d-1}{2}\frac{c/n-1}{2}(r_n'-r_n)}(-1)^{\frac{d-1}{2}\frac{l-1}{2}r_n'} \\
&= \prod_{2 \nmid \frac{c}{n}}\legendre{c}{d}^{r_n'-r_n}\legendre{n}{d}^{r_n'-r_n}\legendre{l}{d}^{r_n'}(-1)^{\frac{d-1}{2}\frac{c/n-1}{2}(r_n'-r_n)}(-1)^{\frac{d-1}{2}\frac{l-1}{2}r_n'}\\
&= \prod_{2 \nmid \frac{c}{n}}\legendre{c}{d}^{r_n'-r_n}\legendre{n}{d}^{r_n'-r_n}\legendre{l}{d}^{r_n'}\etp{-\frac{1}{8}(d-1)\left(l\left(\frac{c}{nl}-1\right)r_n'-\left(\frac{c}{n}-1\right)r_n\right)}.
\end{align*}
Inserting this into \eqref{eq:Qeq3} gives
\begin{multline*}
R = \legendre{\delta}{d}\prod_{2 \nmid \frac{c}{n}}\etp{\frac{1}{8}(d-1)\frac{c}{Nl}\frac{N}{n}(r_n'-lr_n)}\\
\cdot\etp{-\frac{1}{8}(d-1)\left(l\left(\frac{c}{nl}-1\right)r_n'-\left(\frac{c}{n}-1\right)r_n\right)}\prod_{2 \mid \frac{c}{n}}\etp{\frac{1}{8}(d-1)(r_n'-r_n)}.
\end{multline*}
Since $8 \mid (d-1)(l-1)(\frac{c}{nl}-1)$ for $2 \nmid \frac{c}{n}$, we have
\begin{align}
R = \legendre{\delta}{d}&\prod_{2 \nmid \frac{c}{n}}\etp{\frac{1}{8}(d-1)\left(\left(\frac{c}{nl}-1\right)r_n'-\left(\frac{c}{n}-1\right)r_n\right)}\etp{\frac{1}{8}(d-1)(r_n'-r_n)} \notag\\
\cdot&\etp{-\frac{1}{8}(d-1)\left(l\left(\frac{c}{nl}-1\right)r_n'-\left(\frac{c}{n}-1\right)r_n\right)}\prod_{2 \mid \frac{c}{n}}\etp{\frac{1}{8}(d-1)(r_n'-r_n)} \notag\\
&=\legendre{\delta}{d}\prod_{2 \nmid \frac{c}{n}}\etp{\frac{1}{8}(d-1)(r_n'-r_n)}\prod_{2 \mid \frac{c}{n}}\etp{\frac{1}{8}(d-1)(r_n'-r_n)}\notag\\
&=\legendre{\delta}{d}.\label{eq:R2nmidl2midc}
\end{align}
This concludes the proof of the three cases in \eqref{eq:Rcases}, hence the `if' direction.

Conversely, assume that one of the four conditions does not hold. If \eqref{eq:charCom2} does not hold, then $\chi\cdot v_\mathbf{r}\left(\widetilde{T}\right)\neq\chi'\cdot v_\mathbf{r'}\left(\widetilde{T}^l\right)$ by \eqref{eq:etaChar} and \eqref{eq:etaQuoChar}, and hence the characters are not compatible by Lemma \ref{lemm:charCompatible}. If \eqref{eq:charCom1} does not hold, then
\begin{equation*}
\chi'^{-1}\cdot v_{F_N\mathbf{r'}}\left(\widetilde{T}\right)\neq\chi^{-1} \cdot v_{F_N\mathbf{r}}\left(\widetilde{T}^l\right),
\end{equation*}
and hence $\chi'^{-1}\cdot v_{F_N\mathbf{r'}}$ and $\chi^{-1} \cdot v_{F_N\mathbf{r}}$ are not compatible again by Lemma \ref{lemm:charCompatible}. It follows from this and Lemma \ref{lemm:adjointCompatible} that $\chi\cdot v_\mathbf{r}$ and $\chi' \cdot v_\mathbf{r'}$ are not compatible. If \eqref{eq:charCom4} does not hold for $d=-1$, then $\chi\cdot v_\mathbf{r}\left(\widetilde{-I}\right)\neq\chi'\cdot v_\mathbf{r'}\left(\widetilde{-I}\right)$ and non-compatibility follows. Now suppose \eqref{eq:charCom1}, \eqref{eq:charCom2} and \eqref{eq:charCom4} with $d=-1$ all hold. Then \eqref{eq:Qeq3} holds since its validity only relies on \eqref{eq:charCom1} and \eqref{eq:charCom2}. Assume that $\eqref{eq:charCom3}$ does not hold. Then $2 \nmid Nl$ but $\delta \not\equiv 1 \bmod 4$. Since $2 \nmid Nl$, $\delta$ is odd, so $\delta\equiv 3 \bmod 4$. Now for $\tbtmat{a}{b}{c}{d}=\tbtmat{Nl-1}{-1}{Nl}{-1}$ the equation \eqref{eq:R2nmidl2nmidc} holds. The non-compatibility thus follows from this, $\legendre{\delta}{-1}\neq\legendre{-1}{\delta}$ and \eqref{eq:charCom4} with $d=-1$. Finally, assume that \eqref{eq:charCom3} holds but \eqref{eq:charCom4} does not hold for some $d_0 \in \numZ_{\neq 0}$ with $\gcd(Nl,d_0)=1$. If $2 \mid l$, then the non-compatibility follows from \eqref{eq:R2midl} with $\tbtmat{a}{b}{c}{d}=\tbtmat{a_0}{b_0}{Nl}{d_0}\in\Gz{Nl}$. Otherwise if $2 \nmid l$, then the non-compatibility follows from \eqref{eq:R2nmidl2nmidc}, \eqref{eq:charCom3} in the case $2 \nmid N$ and from \eqref{eq:R2nmidl2midc} in the case $2 \mid N$ with $\tbtmat{a}{b}{c}{d}=\tbtmat{a_0}{b_0}{Nl}{d_0}\in\Gz{Nl}$. Thereby we have shown that if one of the four given conditions does not hold, then the given characters are not $\widetilde{\tbtmat{1}{0}{0}{l}}$-compatible.

The conditions on compatibility of $v_\mathbf{r}$ and $v_\mathbf{r'}$ follow immediately from the general case.
\end{proof}

The full modular group case, that is, the case $N=1$, is worthwhile to state independently. Note that the group of complex linear characters on $\sltZ$ is a cyclic group of order $24$ with a single generator $v_\eta$. Hence the following corollary covers all characters on $\sltZ$.
\begin{coro}
\label{coro:compatibleFullmodularGroup}
Let $r$ be an integer, and $l$ be a positive integer.
\begin{enumerate}
\item If $r$ is even, then the characters $v_\eta^r$ and $v_\eta^r$ are $\widetilde{\tbtmat{1}{0}{0}{l}}$-compatible if and only if $24 \mid r(l-1)$.
\item If $r$ is odd, then the characters $v_\eta^r$ and $v_\eta^r$ are $\widetilde{\tbtmat{1}{0}{0}{l}}$-compatible if and only if $24 \mid r(l-1)$, and $l$ is a square.
\end{enumerate}
\end{coro}
\begin{proof}
This follows immediately from Theorem \ref{thm:charCompatibleEtaDirichlet}.
\end{proof}

\section{Operators Acting on One-dimensional Spaces}
\label{sec:one-dimensional spaces}
As a first application, we try to find out some one-dimensional spaces of modular forms, and $T_l$'s acting on them. For two such spaces $M_1$ and $M_2$ and an operator $T_l$ from $M_1$ to $M_2$, there must exist some $c \in \numC$ such that $T_l f_1=cf_2$ for any $f_1 \in M_1$ and $f_2 \in M_2$. The aim of this section is to find such identities among certain eta-quotients.

\begin{lemm}
\label{lemm:dimone}
Let $N$ be a positive integer, $k$ be a positive integer or half-integer and $v$ be a complex linear character on $\Gzt{N}$. If $kN\prod_{p \mid N} \left(1+p^{-1}\right)<12$ where $p$ denotes primes, then we have $\dim_\numC M_k(\Gzt{N}, v) \leq 1$. Moreover, if $kN\prod_{p \mid N} \left(1+p^{-1}\right)=12$, then $\dim_\numC S_k(\Gzt{N}, v) \leq 1$.
\end{lemm}
\begin{proof}
Let $f,\, g \in M_k(\Gzt{N}, v)$ be nonzero forms. We shall prove they are $\numC$-linearly dependent, from which the estimate on $\dim_\numC M_k(\Gzt{N}, v)$ follows. Since they are nonzero, their sets of zeros are discrete in $\uhp$. Since $\Gzp{N}_{\tau} \subseteq \pslZ_{\tau}$, the set of $\tau\in\uhp$ with $\abs{\Gzp{N}_{\tau}} > 1$ is discrete. Therefore we can find a $\tau_0\in\uhp$ such that $\abs{\Gzp{N}_{\tau_0}}=1$, $f(\tau_0) \neq 0$ and $g(\tau_0) \neq 0$. Set $c_1=g(\tau_0)$ and $c_2=-f(\tau_0)$, and consider $h=c_1f+c_2g$. If $h \neq 0$ then by Theorem \ref{thm:valence} we have $1 \leq \ord_{\tau_0}(h) \leq \frac{1}{12}[\pslZ : \Gzp{N}]k$, which contradicts the assumption since $[\pslZ : \Gzp{N}]=[\slZ : \Gz{N}]=N\prod_p \left(1+p^{-1}\right)$ (see, for instance, \cite[Section 6.2]{CoS17}). Hence $h=0$, that is, $f$ and $g$ are dependent. The latter assertion can be proved in a similar manner with the inequalities on $\ord_{\tau_0}(h)$ replaced by
\begin{equation*}
 1 \leq \ord\nolimits_{\tau_0}(h) < \sum_{\overline{\gamma} \in \overline{G}\backslash\pslZ}\ord\nolimits_{\rmi\infty}(h\vert_k \gamma) + \ord\nolimits_{\tau_0}(h) \leq \frac{1}{12}[\pslZ : \Gzp{N}]k.
\end{equation*}
\end{proof}
We list all pairs $(N, k)$ satisfying $kN\prod_{p \mid N} \left(1+p^{-1}\right)<12$ in the following table.
\begin{table}[ht]
\centering
\caption{The level $N$ and weight $k$ such that $[\pslZ : \Gzp{N}]k<12$ \label{table:NkDimone}}
\begin{tabular}{llll}
\toprule
$N$ & $k$ & N & k \\
\midrule
$1$ & $1/2,\,1, \dots, 23/2$ & $8$ & $1/2$ \\
$2$ & $1/2,\,1, \dots, 7/2$ & $9$ & $1/2$ \\
$3$ & $1/2,\,1, \dots, 5/2$ & $10$ & $1/2$ \\
$4$ & $1/2,\, 1,\, 3/2$ & $11$ & $1/2$ \\
$5$ & $1/2,\, 1,\, 3/2$ & $13$ & $1/2$ \\
$6$ & $1/2$ & $17$ & $1/2$ \\
$7$ & $1/2,\, 1$ & $19$ & $1/2$ \\
\bottomrule
\end{tabular}
\end{table}

For each $(N, k)$ in this table, using \eqref{eq:etaOrder}, we find out all holomorphic eta-quotients on group $\Gzt{N}$, of weight $k$ and with some character. The results are listed in Table \ref{table:EtaDimonePrime} for $N$ being a prime or $1$, and in Table \ref{table:EtaDimoneNonPrime} for $N=4, 6, 8, 9, 10$. These tables are generated by SageMath programs. See \ref{sec:Usage of code}.
\begin{table}[ht]
\centering
\caption{For each $(N, k)$ in Table \ref{table:NkDimone} with $N$ being a prime or $1$ the set of all holomorphic eta-quotients on $\Gzt{N}$ of weight $k$, where $\eta^r_n$ denoted $\eta^r(n\tau)$ \label{table:EtaDimonePrime}}
\begin{tabular}{ll}
\toprule
$(N, k)$ & Holomorphic eta-quotients \\
\midrule
$(1, *)$ & $\eta^{1}_1,\,\eta^2_1,\dots,\eta^{23}_1$ \\
$(2,*)$ & $\eta^{-2k+i}_1\eta^{4k-i}_2$ for $i=0, 1,\dots, 6k$ \\
$(3, *)$ & $\eta^{-[k]+i}_1\eta^{2k+[k]-i}_3$ for $i=0, 1,\dots, 2k+2[k]$ \\
$(5,*)$ & $\eta^{i}_1\eta^{2k-i}_5$ for $i=0, 1,\dots, 2k$ \\
$(7,*)$ & $\eta^{i}_1\eta^{2k-i}_7$ for $i=0, 1,\dots, 2k$ \\
$(11,1/2)$ & $\eta^{i}_1\eta^{1-i}_{11}$ for $i=0, 1$ \\
$(13,1/2)$ & $\eta^{i}_1\eta^{1-i}_{13}$ for $i=0, 1$ \\
$(17,1/2)$ & $\eta^{i}_1\eta^{1-i}_{17}$ for $i=0, 1$ \\
$(19,1/2)$ & $\eta^{i}_1\eta^{1-i}_{19}$ for $i=0, 1$ \\
\bottomrule
\end{tabular}
\end{table}
\begin{table}[p]
\centering
\caption{For each $(N, k)$ in Table \ref{table:NkDimone} with $N=4, 6, 8, 9, 10$ the set of all holomorphic eta-quotients on $\Gzt{N}$ of weight $k$, where $\eta^r_n$ denoted $\eta^r(n\tau)$ \label{table:EtaDimoneNonPrime}}
\begin{tabular}{llllll}
\toprule
$(N, k)$ & \multicolumn{5}{l}{Holomorphic eta-quotients} \\
\midrule
$(4, 1/2)$ & $\eta^{-2}_1\eta^{5}_2\eta^{-2}_4$& $\eta^{-1}_1\eta^{2}_2\eta^{0}_4$& $\eta^{-1}_1\eta^{3}_2\eta^{-1}_4$& $\eta^{0}_1\eta^{-1}_2\eta^{2}_4$& $\eta^{0}_1\eta^{0}_2\eta^{1}_4$\\
& $\eta^{0}_1\eta^{1}_2\eta^{0}_4$& $\eta^{0}_1\eta^{2}_2\eta^{-1}_4$& $\eta^{1}_1\eta^{-1}_2\eta^{1}_4$& $\eta^{1}_1\eta^{0}_2\eta^{0}_4$& $\eta^{2}_1\eta^{-1}_2\eta^{0}_4$\\
$(4, 1)$ & $\eta^{-4}_1\eta^{10}_2\eta^{-4}_4$& $\eta^{-3}_1\eta^{7}_2\eta^{-2}_4$& $\eta^{-3}_1\eta^{8}_2\eta^{-3}_4$& $\eta^{-2}_1\eta^{4}_2\eta^{0}_4$& $\eta^{-2}_1\eta^{5}_2\eta^{-1}_4$\\
& $\eta^{-2}_1\eta^{6}_2\eta^{-2}_4$& $\eta^{-2}_1\eta^{7}_2\eta^{-3}_4$& $\eta^{-1}_1\eta^{1}_2\eta^{2}_4$& $\eta^{-1}_1\eta^{2}_2\eta^{1}_4$& $\eta^{-1}_1\eta^{3}_2\eta^{0}_4$\\
& $\eta^{-1}_1\eta^{4}_2\eta^{-1}_4$& $\eta^{-1}_1\eta^{5}_2\eta^{-2}_4$& $\eta^{0}_1\eta^{-2}_2\eta^{4}_4$& $\eta^{0}_1\eta^{-1}_2\eta^{3}_4$& $\eta^{0}_1\eta^{0}_2\eta^{2}_4$\\
& $\eta^{0}_1\eta^{1}_2\eta^{1}_4$& $\eta^{0}_1\eta^{2}_2\eta^{0}_4$& $\eta^{0}_1\eta^{3}_2\eta^{-1}_4$& $\eta^{0}_1\eta^{4}_2\eta^{-2}_4$& $\eta^{1}_1\eta^{-2}_2\eta^{3}_4$\\
& $\eta^{1}_1\eta^{-1}_2\eta^{2}_4$& $\eta^{1}_1\eta^{0}_2\eta^{1}_4$& $\eta^{1}_1\eta^{1}_2\eta^{0}_4$& $\eta^{1}_1\eta^{2}_2\eta^{-1}_4$& $\eta^{2}_1\eta^{-2}_2\eta^{2}_4$\\
& $\eta^{2}_1\eta^{-1}_2\eta^{1}_4$& $\eta^{2}_1\eta^{0}_2\eta^{0}_4$& $\eta^{2}_1\eta^{1}_2\eta^{-1}_4$& $\eta^{3}_1\eta^{-2}_2\eta^{1}_4$& $\eta^{3}_1\eta^{-1}_2\eta^{0}_4$\\
& $\eta^{4}_1\eta^{-2}_2\eta^{0}_4$&&&&\\
$(4, 3/2)$ & $\eta^{-6}_1\eta^{15}_2\eta^{-6}_4$& $\eta^{-5}_1\eta^{12}_2\eta^{-4}_4$& $\eta^{-5}_1\eta^{13}_2\eta^{-5}_4$& $\eta^{-4}_1\eta^{9}_2\eta^{-2}_4$& $\eta^{-4}_1\eta^{10}_2\eta^{-3}_4$\\
& $\eta^{-4}_1\eta^{11}_2\eta^{-4}_4$& $\eta^{-4}_1\eta^{12}_2\eta^{-5}_4$& $\eta^{-3}_1\eta^{6}_2\eta^{0}_4$& $\eta^{-3}_1\eta^{7}_2\eta^{-1}_4$& $\eta^{-3}_1\eta^{8}_2\eta^{-2}_4$\\
& $\eta^{-3}_1\eta^{9}_2\eta^{-3}_4$& $\eta^{-3}_1\eta^{10}_2\eta^{-4}_4$& $\eta^{-2}_1\eta^{3}_2\eta^{2}_4$& $\eta^{-2}_1\eta^{4}_2\eta^{1}_4$& $\eta^{-2}_1\eta^{5}_2\eta^{0}_4$\\
& $\eta^{-2}_1\eta^{6}_2\eta^{-1}_4$& $\eta^{-2}_1\eta^{7}_2\eta^{-2}_4$& $\eta^{-2}_1\eta^{8}_2\eta^{-3}_4$& $\eta^{-2}_1\eta^{9}_2\eta^{-4}_4$& $\eta^{-1}_1\eta^{0}_2\eta^{4}_4$\\
& $\eta^{-1}_1\eta^{1}_2\eta^{3}_4$& $\eta^{-1}_1\eta^{2}_2\eta^{2}_4$& $\eta^{-1}_1\eta^{3}_2\eta^{1}_4$& $\eta^{-1}_1\eta^{4}_2\eta^{0}_4$& $\eta^{-1}_1\eta^{5}_2\eta^{-1}_4$\\
& $\eta^{-1}_1\eta^{6}_2\eta^{-2}_4$& $\eta^{-1}_1\eta^{7}_2\eta^{-3}_4$& $\eta^{0}_1\eta^{-3}_2\eta^{6}_4$& $\eta^{0}_1\eta^{-2}_2\eta^{5}_4$& $\eta^{0}_1\eta^{-1}_2\eta^{4}_4$\\
& $\eta^{0}_1\eta^{0}_2\eta^{3}_4$& $\eta^{0}_1\eta^{1}_2\eta^{2}_4$& $\eta^{0}_1\eta^{2}_2\eta^{1}_4$& $\eta^{0}_1\eta^{3}_2\eta^{0}_4$& $\eta^{0}_1\eta^{4}_2\eta^{-1}_4$\\
& $\eta^{0}_1\eta^{5}_2\eta^{-2}_4$& $\eta^{0}_1\eta^{6}_2\eta^{-3}_4$& $\eta^{1}_1\eta^{-3}_2\eta^{5}_4$& $\eta^{1}_1\eta^{-2}_2\eta^{4}_4$& $\eta^{1}_1\eta^{-1}_2\eta^{3}_4$\\
& $\eta^{1}_1\eta^{0}_2\eta^{2}_4$& $\eta^{1}_1\eta^{1}_2\eta^{1}_4$& $\eta^{1}_1\eta^{2}_2\eta^{0}_4$& $\eta^{1}_1\eta^{3}_2\eta^{-1}_4$& $\eta^{1}_1\eta^{4}_2\eta^{-2}_4$\\
& $\eta^{2}_1\eta^{-3}_2\eta^{4}_4$& $\eta^{2}_1\eta^{-2}_2\eta^{3}_4$& $\eta^{2}_1\eta^{-1}_2\eta^{2}_4$& $\eta^{2}_1\eta^{0}_2\eta^{1}_4$& $\eta^{2}_1\eta^{1}_2\eta^{0}_4$\\
& $\eta^{2}_1\eta^{2}_2\eta^{-1}_4$& $\eta^{2}_1\eta^{3}_2\eta^{-2}_4$& $\eta^{3}_1\eta^{-3}_2\eta^{3}_4$& $\eta^{3}_1\eta^{-2}_2\eta^{2}_4$& $\eta^{3}_1\eta^{-1}_2\eta^{1}_4$\\
& $\eta^{3}_1\eta^{0}_2\eta^{0}_4$& $\eta^{3}_1\eta^{1}_2\eta^{-1}_4$& $\eta^{4}_1\eta^{-3}_2\eta^{2}_4$& $\eta^{4}_1\eta^{-2}_2\eta^{1}_4$& $\eta^{4}_1\eta^{-1}_2\eta^{0}_4$\\
& $\eta^{4}_1\eta^{0}_2\eta^{-1}_4$& $\eta^{5}_1\eta^{-3}_2\eta^{1}_4$& $\eta^{5}_1\eta^{-2}_2\eta^{0}_4$& $\eta^{6}_1\eta^{-3}_2\eta^{0}_4$&\\
$(6, 1/2)$ & $\eta^{-1}_1\eta^{1}_2\eta^{2}_3\eta^{-1}_6$& $\eta^{-1}_1\eta^{2}_2\eta^{0}_3\eta^{0}_6$& $\eta^{-1}_1\eta^{2}_2\eta^{1}_3\eta^{-1}_6$& $\eta^{0}_1\eta^{0}_2\eta^{-1}_3\eta^{2}_6$& $\eta^{0}_1\eta^{0}_2\eta^{0}_3\eta^{1}_6$\\
& $\eta^{0}_1\eta^{0}_2\eta^{1}_3\eta^{0}_6$& $\eta^{0}_1\eta^{0}_2\eta^{2}_3\eta^{-1}_6$& $\eta^{0}_1\eta^{1}_2\eta^{0}_3\eta^{0}_6$& $\eta^{1}_1\eta^{-1}_2\eta^{-1}_3\eta^{2}_6$& $\eta^{1}_1\eta^{0}_2\eta^{0}_3\eta^{0}_6$\\
& $\eta^{2}_1\eta^{-1}_2\eta^{-1}_3\eta^{1}_6$& $\eta^{2}_1\eta^{-1}_2\eta^{0}_3\eta^{0}_6$&&&\\
$(8, 1/2)$ & $\eta^{-2}_1\eta^{5}_2\eta^{-2}_4\eta^{0}_8$& $\eta^{-1}_1\eta^{2}_2\eta^{0}_4\eta^{0}_8$& $\eta^{-1}_1\eta^{3}_2\eta^{-1}_4\eta^{0}_8$& $\eta^{0}_1\eta^{-2}_2\eta^{5}_4\eta^{-2}_8$& $\eta^{0}_1\eta^{-1}_2\eta^{2}_4\eta^{0}_8$\\
& $\eta^{0}_1\eta^{-1}_2\eta^{3}_4\eta^{-1}_8$& $\eta^{0}_1\eta^{0}_2\eta^{-1}_4\eta^{2}_8$& $\eta^{0}_1\eta^{0}_2\eta^{0}_4\eta^{1}_8$& $\eta^{0}_1\eta^{0}_2\eta^{1}_4\eta^{0}_8$& $\eta^{0}_1\eta^{0}_2\eta^{2}_4\eta^{-1}_8$\\
& $\eta^{0}_1\eta^{1}_2\eta^{-1}_4\eta^{1}_8$& $\eta^{0}_1\eta^{1}_2\eta^{0}_4\eta^{0}_8$& $\eta^{0}_1\eta^{2}_2\eta^{-1}_4\eta^{0}_8$& $\eta^{1}_1\eta^{-1}_2\eta^{1}_4\eta^{0}_8$& $\eta^{1}_1\eta^{0}_2\eta^{0}_4\eta^{0}_8$\\
& $\eta^{2}_1\eta^{-1}_2\eta^{0}_4\eta^{0}_8$&&&&\\
$(9, 1/2)$ & $\eta^{0}_1\eta^{0}_3\eta^{1}_9$& $\eta^{0}_1\eta^{1}_3\eta^{0}_9$& $\eta^{1}_1\eta^{0}_3\eta^{0}_9$&&\\
$(10, 1/2)$ & $\eta^{-1}_1\eta^{2}_2\eta^{0}_5\eta^{0}_{10}$& $\eta^{0}_1\eta^{0}_2\eta^{-1}_5\eta^{2}_{10}$& $\eta^{0}_1\eta^{0}_2\eta^{0}_5\eta^{1}_{10}$& $\eta^{0}_1\eta^{0}_2\eta^{1}_5\eta^{0}_{10}$& $\eta^{0}_1\eta^{0}_2\eta^{2}_5\eta^{-1}_{10}$\\
& $\eta^{0}_1\eta^{1}_2\eta^{0}_5\eta^{0}_{10}$& $\eta^{1}_1\eta^{0}_2\eta^{0}_5\eta^{0}_{10}$& $\eta^{2}_1\eta^{-1}_2\eta^{0}_5\eta^{0}_{10}$&&\\
\bottomrule
\end{tabular}
\end{table}

For any two eta-quotients $\eta^{\mathbf{r}}$ and $\eta^{\mathbf{r'}}$ listed in Table \ref{table:EtaDimonePrime} and Table \ref{table:EtaDimoneNonPrime} on the same group $\Gzt{N}$ of the same weight $k$, we can use Theorem \ref{thm:charCompatibleEtaDirichlet} to find positive integers $l$ such that $v_{\mathbf{r}}$ and $v_{\mathbf{r'}}$ are $\widetilde{\tbtmat{1}{0}{0}{l}}$-compatible. Such $l$ should satisfy \eqref{eq:charCom1}, \eqref{eq:charCom2}, and that $\delta=l^{2k}\prod_{2 \nmid r_n-r_n'}n$ is a perfect square. Once such an $l$ is found, we have $T_l \eta^{\mathbf{r}} = c\eta^{\mathbf{r'}}$ for some $c \in \numC$, as we aim to seek. All these can be done by a computer algebra system.

\subsection{The case of integral weights}
In this subsection, $k$ always denotes positive integers. This is an easier case, since $l^{2k}$ is already a perfect square. Hence $l^{2k}\prod_{2 \nmid r_n-r_n'}n$ is a perfect square if and only if $\prod_{2 \nmid r_n-r_n'}n$ is, and the latter condition is independent of $l$. Thus we need only find solutions $l$ among integers from $1$ to $24$, and all solutions are of the form $l_0+24j$ for $j \in \numgeq{Z}{0}$ and $l_0$ being a solution in  $\{1, 2, \dots, 24\}$. All tuples $(\eta^{\mathbf{r}}, \eta^{\mathbf{r'}}, l)$ have been found by calling a SageMath function \lstinline{findl(N, k)}, which are listed in Table \ref{table:Lintegral} in \ref{apx:tables}. For the usage of SageMath code, see \ref{sec:Usage of code}.

We summarize:
\begin{thm}
\label{thm:TlIntegralWeight}
Let $(N, k)$ be a pair listed in Table \ref{table:NkDimone} with $k$ an integer. Let $\eta^{\mathbf{r}}$ and $\eta^{\mathbf{r'}}$ be holomorphic eta-quotients on $\Gzt{N}$ of weight $k$. Then $\eta^{\mathbf{r}}$ and $\eta^{\mathbf{r'}}$ must be functions listed in Table \ref{table:EtaDimonePrime} or in Table \ref{table:EtaDimoneNonPrime}, of the corresponding $(N, k)$. Moreover, let $l$ be a positive integer. Then $\mathbf{r}$, $\mathbf{r'}$ and $l$ satisfy \eqref{eq:charCom1}, \eqref{eq:charCom2}, and that $\prod_{2 \nmid r_n'-r_n}n$ is a square, if and only if $(\eta^\mathbf{r}, \eta^\mathbf{r'}, l)$ is an entry in Table \ref{table:Lintegral}. For such a tuple $(\eta^\mathbf{r}, \eta^\mathbf{r'}, l)$, there exists some $c \in \numC$, such that $T_l \eta^\mathbf{r}=c\eta^\mathbf{r'}$.
\end{thm}

We illustrate this theorem by an example, in which the formula $T_l \eta^\mathbf{r}=c\eta^\mathbf{r'}$ is evaluated more explicitly.
\begin{examp}
\label{examp:onedimExamp1}
Look at the entry $(\eta^{-3}_1\eta^7_2, \eta^1_1\eta^3_2)$ in Table \ref{table:Lintegral} for which $N=2$ and $k=2$. Set $\mathbf{r}=(-3, 7)$, $\mathbf{r'}=(1, 3)$ and $\alpha=\widetilde{\tbtmat{1}{0}{0}{l}}$. The numbers $l$ such that $v_\mathbf{r}$ and $v_\mathbf{r'}$ are $\alpha$-compatible are precisely, by Theorem \ref{thm:charCompatibleEtaDirichlet}, $l=5+24j$, with $j$ being non-negative integers. For simplicity, we only consider the case $l$ is a prime: $l=5$, $29$, $53$, $101$, $149$, $173$, $197$, $269$, $293$, $317$, $389,\dots$. Note that
\begin{equation*}
\widetilde{\tbtMat{l}{0}{0}{1}}=\widetilde{\tbtMat{l}{-\frac{l-1}{2}}{-2}{1}}\widetilde{\tbtMat{1}{0}{0}{l}}\widetilde{\tbtMat{l}{\frac{l-1}{2}}{2}{1}}.
\end{equation*}
It follows that
\begin{equation*}
\left(v_\mathbf{r}\vert_\alpha v_{\mathbf{r'}}\right)^{-1}\widetilde{\tbtMat{l}{0}{0}{1}}=v_\mathbf{r}^{-1}\widetilde{\tbtMat{l}{-\frac{l-1}{2}}{-2}{1}}v_\mathbf{r'}^{-1}\widetilde{\tbtMat{l}{\frac{l-1}{2}}{2}{1}}=-1,
\end{equation*}
by \eqref{eq:etaChar} and \eqref{eq:etaQuoChar}. Hence, using \eqref{eq:doubleCosetOperators} and \eqref{eq:doubleDecompLeft} with $m=1$, we have
\begin{multline*}
T_l \eta^{-3}(\tau)\eta^7(2\tau)=-l\eta^{-3}(l\tau)\eta^7(2l\tau) \\
+\frac{1}{l}\sum_{\lambda=0}^{l-1}\etp{-\frac{7\lambda}{24}}\eta^{-3}\left(\frac{\tau+\lambda}{l}\right)\eta^7\left(\frac{2\tau+2\lambda}{l}\right).
\end{multline*}
Thus, theorem \ref{thm:TlIntegralWeight} tells us that, for any prime $l$ congruent to $5$ modulo $24$  there exists a constant $c_l \in \numC$, depending on $l$, such that
\begin{multline}
\label{eq:identityEtaQuotExam}
-l\eta^{-3}(l\tau)\eta^7(2l\tau)+ \frac{1}{l}\sum_{\lambda=0}^{l-1}\etp{-\frac{7\lambda}{24}}\eta^{-3}\left(\frac{\tau+\lambda}{l}\right)\eta^7\left(\frac{2\tau+2\lambda}{l}\right) \\
=c_l\eta(\tau)\eta^3(2\tau).
\end{multline}
To determine $c_l$, express this identity in form of Fourier series. Suppose that
\begin{equation*}
\eta^{-3}(\tau)\eta^7(2\tau)=q^{11/24}\sum_{n \in \numgeq{Z}{0}}C(n)q^n,
\end{equation*}
and set $C(n) = 0$ if $n \in \numQ$ is not a non-negative integer. A straightforward calculation shows that
\begin{multline*}
T_l \eta^{-3}(\tau)\eta^7(2\tau)\\
=q^{\frac{7}{24}}\sum_{n \in \numZ}\left(C\left(\frac{7l-11}{24}+ln\right)-lC\left(\frac{1}{l}\left(-\frac{11l-7}{24}+n\right)\right)\right)q^n.
\end{multline*}
By comparing first terms of Fourier series expansions of $T_l \eta^{-3}(\tau)\eta^7(2\tau)$ and $c_l\eta(\tau)\eta^3(2\tau)$, we find that $c_l=C\left(\frac{7l-11}{24}\right)$.
Some non-trivial congruence relations can be deduced from such identities. If $C\left(\frac{7l-11}{24}\right)\neq 0$, then for any integer $n$, we have
\begin{equation}
\label{eq:congEtaQuotExam}
C\left(\frac{7l-11}{24}+ln\right) \equiv lC\left(\frac{1}{l}\left(-\frac{11l-7}{24}+n\right)\right) \pmod{C\left(\frac{7l-11}{24}\right)}.
\end{equation}
For instance, take $l=13709 = 5+24\cdot 571$, which is a prime. Then \eqref{eq:congEtaQuotExam} becomes
\begin{equation*}
C(3998+13709n) \equiv 41C\left(\frac{n-6283}{13709}\right) \pmod{3\cdot 67}.
\end{equation*}
We present the first ten coefficients in Table \ref{table:C3998+13709n}.
\begin{table}[b]
\centering
\caption{The coefficients $C(3998+13709n)$ of $\eta^{-3}(\tau)\eta^7(2\tau)$$=q^{11/24}\sum_{n \in \numgeq{Z}{0}}C(n)q^n$ for $n=0,1,2,\dots,9$\label{table:C3998+13709n}}
\begin{tabular}{llll}
\toprule
$n$ & $C(3998+13709n)$ & $n$ & $C(3998+13709n)$ \\
\midrule
$0$ & $-3\cdot67$ & $5$ & $-3\cdot67$\\
$1$ & $3\cdot67$ & $6$ & $-3\cdot5\cdot67$\\
$2$ & $2^2\cdot3\cdot67$ & $7$ & $3\cdot7\cdot67$\\
$3$ & $-3^2\cdot67$ & $8$ & $3\cdot5\cdot67$\\
$4$ & $-3^2\cdot67$ & $9$ & $3^2\cdot67$\\
\bottomrule
\end{tabular}
\end{table}

In fact, for all function pairs listed in Table \ref{table:Lintegral}, we have identities like \eqref{eq:identityEtaQuotExam} and congruences like \eqref{eq:congEtaQuotExam} that we have presented for $(\eta^{-3}_1\eta^7_2, \eta^1_1\eta^3_2)$.
\end{examp}

\subsection{The case of half-integral weights}
In this subsection, let $k \in \frac{1}{2}+\numgeq{Z}{0}$. It is a bit harder to handle this case, since whether $\delta=l^{2k}\prod_{2 \nmid r_n'-r_n}n$ is a square depends on $l$. Nevertheless, let $l_0$ be the least positive integer such that $l_0\prod_{2 \nmid r_n'-r_n}n$ is a square. Then $l^{2k}\prod_{2 \nmid r_n'-r_n}n$ is a square if and only if $l=l_0m^2$ for some positive integer $m$. Hence, we solve the congruences \eqref{eq:charCom1} and \eqref{eq:charCom2}, with $l=l_0m^2$ for the variable $m\in \numgeq{Z}{1}$. All solutions $(\eta^{\mathbf{r}}, \eta^{\mathbf{r'}}, l_0, m)$ have been found, again by calling the SageMath function \lstinline{findl(N, k)}, which are listed in Table \ref{table:Lhalfintegral} in \ref{apx:tables}. For the SageMath code, see \ref{sec:Usage of code}.

We summarize as follows.
\begin{thm}
\label{thm:TlHalfintegralWeight}
Let $(N, k)$ be a pair listed in Table \ref{table:NkDimone} with $k$ a half-integer. Let $\eta^{\mathbf{r}}$ and $\eta^{\mathbf{r'}}$ be holomorphic eta-quotients on $\Gzt{N}$ of weight $k$. Then $\eta^{\mathbf{r}}$ and $\eta^{\mathbf{r'}}$ must be functions listed in Table \ref{table:EtaDimonePrime} or in Table \ref{table:EtaDimoneNonPrime}, of the corresponding $(N, k)$. Moreover, let $l$ be a positive integer. Then $\mathbf{r}$, $\mathbf{r'}$ and $l$ satisfy \eqref{eq:charCom1}, \eqref{eq:charCom2}, and that $l^{2k}\prod_{2 \nmid r_n'-r_n}n$ is a square, if and only if there is some $m\in \numgeq{Z}{1}$, such that $l=l_0m^2$, and that $(\eta^\mathbf{r}, \eta^\mathbf{r'}, l_0, m^2 \bmod{24})$ is an entry in Table \ref{table:Lhalfintegral}. For such a tuple $(\eta^\mathbf{r}, \eta^\mathbf{r'}, l)$, we have $T_l \eta^\mathbf{r}=c\eta^\mathbf{r'}$, with some $c \in \numC$.
\end{thm}

\begin{examp}
\label{examp:onedimExamp2}
Look at the entry $(\eta^{-6}_1\eta^{15}_2\eta^{-6}_4, \eta^{-3}_1\eta^6_2\eta^0_4)$ in Table \ref{table:Lhalfintegral}, for which $N=4$ and $k=3/2$. The corresponding characters are $\widetilde{\tbtmat{1}{0}{0}{l}}$-compatible, if and only if $l=2m^2$, where $m^2 \equiv 0,2,4,16 \bmod{24}$ and $m > 0$ by Theorem \ref{thm:charCompatibleEtaDirichlet}. For simplicity, we only consider the case $m=2^{\beta}$ with $\beta \in \numgeq{Z}{1}$. Thus $\rad{l} \mid \rad{N}$, and consequently Proposition \ref{prop:radlcaseFourier} applies. Hence we have
\begin{multline*}
T_l \eta^{-6}(\tau)\eta^{15}(2\tau)\eta^{-6}(4\tau) = 2^{-\frac{6\beta+3}{4}} \\
\cdot\sum_{0\leq \lambda < 2^{2\beta+1}}\etp{-\frac{3}{8}\lambda}\eta^{-6}\left(\frac{\tau+\lambda}{2^{2\beta+1}}\right)\eta^{15}\left(\frac{\tau+\lambda}{2^{2\beta}}\right)\eta^{-6}\left(\frac{\tau+\lambda}{2^{2\beta-1}}\right).
\end{multline*} 
Note that $\eta^{-2}(\tau)\eta^{5}(2\tau)\eta^{-2}(4\tau)$ is just the classical theta series $\theta(\tau)=\sum_{n \in \numZ}q^{n^2}$ and that $\eta^{-1}(\tau)\eta^2(2\tau)=\theta_1(\tau)$ where $\theta_1(\tau)=\frac{1}{2}\sum_{n \in \numZ}\legendre{4}{n}q^{n^2/8}=q^{1/8}\sum_{n \geq 0}q^{n(n+1)/2}$. To prove these identities compare Fourier expansions and use either the Sturm bound or Theorem \ref{thm:valence}. For $\theta_1(\tau)$, see \cite{ORW95} for details.
Thus, we have
\begin{gather*}
\eta^{-6}(\tau)\eta^{15}(2\tau)\eta^{-6}(4\tau)=\sum_{n \geq 0}r_3(n)q^n, \\
\eta^{-3}(\tau)\eta^6(2\tau)=q^{3/8}\sum_{n \geq 0}r_3'(n)q^n,
\end{gather*}
where $r_3(n)$ represents the number of solutions of $n=x^2+y^2+z^2$ for $x,\,y,\,z \in \numZ$, and $r_3'(n)$ represents the number of ways expressing $n$ as a sum of three triangular numbers (including $0(0+1)/2$ and taking into account the order of summands). Now by \eqref{eq:radlcaseFourier2} in Proposition \ref{prop:radlcaseFourier}, we have
\begin{equation*}
T_l \eta^{-6}(\tau)\eta^{15}(2\tau)\eta^{-6}(4\tau) =2^{\frac{2\beta+1}{4}}q^{3/8}\sum_{n \geq 0}r_3(3\cdot4^{\beta-1}+2^{2\beta+1}n)q^n
\end{equation*}
Apply Theorem \ref{thm:TlHalfintegralWeight}; we obtain
\begin{equation}
r_3(3\cdot4^{\beta-1}+2^{2\beta+1}n)=r_3(3\cdot4^{\beta-1})r_3'(n),
\end{equation}
for any $\beta \in \numgeq{Z}{1}$ and $n \in \numgeq{Z}{0}$.
\end{examp}

We can find more one-dimensional spaces, if we restrict the consideration to cusp forms. We list all possible pairs $(N, k)$ with the property $[\pslZ : \Gzp{N}]k=12$ in Table \ref{table:Nk12}, for which we have $\dim_\numC S_k(\Gzt{N}, v) \leq 1$ by Lemma \ref{lemm:dimone}.
\begin{table}[ht]
\centering
\caption{All pairs $(N,k)$ such that $[\pslZ : \Gzp{N}]k=12$ \label{table:Nk12}}
\begin{tabular}{lllll}
\toprule
$(1, 12)$ & $(2, 4)$ & $(3, 3)$ &  $(4, 2)$ & $(5, 2)$ \\
$(6, 1)$ & $(7, 3/2)$ & $(8, 1)$ &$(9, 1)$   & $(11, 1)$ \\
$(12, 1/2)$&$ (14, 1/2)$&$ (15, 1/2)$&$(16, 1/2)$&$ (23, 1/2)$ \\
\bottomrule
\end{tabular}
\end{table}
There are many more holomorphic eta-quotients that are also cusp forms for these $(N, k)$'s than those considered in Theorem \ref{thm:TlIntegralWeight} and Theorem \ref{thm:TlHalfintegralWeight}. Listing all of them and the admissible operators acting on them would lead to a too long paper. We only present the case $(N, k)=(4, 2)$. There are $85$ holomorphic eta-quotients that are also cusp forms, of weight $2$ on the group $\Gzt{4}$. All pairs of such forms and the admissible operators between them are showed in Table \ref{table:N4k2L} (There may be some pairs of such eta-quotients with no $T_l$ acting on them. Such pairs are not listed in Table \ref{table:N4k2L}). By Theorem \ref{thm:charCompatibleEtaDirichlet}, for any tuple $(\eta^\mathbf{r}, \eta^\mathbf{r'}, l)$ in Table \ref{table:N4k2L}, there is a constant $c \in \numC$ such that $T_l \eta^\mathbf{r}=c\eta^\mathbf{r'}$. The reader may use SageMath code to produce those pairs of eta-quotients and admissible operators not listed here. See \ref{sec:Usage of code} for details.

\section{Eta-powers and Operators Acting on Them}
\label{sec:eta-powers}
This section is devoted to eta-powers $\eta^r(\tau)$ with $r=0, 1,2,\dots,24$, and operators acting on them, as the second application. Although Theorem \ref{thm:TlIntegralWeight} and \ref{thm:TlHalfintegralWeight} include the following lemma (except for $r=0$, which is a trivial case), we prefer to state it independently.
\begin{lemm}
\label{lemm:etaPowerEigen}
Suppose that $r \in \{0,\,1,\,2,\dots,24\}$, and that $l \in \numgeq{Z}{1}$ satisfying $24 \mid r(l-1)$. If in addition, $r$ is odd, then we require that $l$ is a square. Then there is a constant $c$, depending on $r$ and $l$, such that
\begin{equation}
T_l \eta^r = c\cdot\eta^r.
\end{equation}
\end{lemm}
\begin{proof}
This follows immediately from Corollary \ref{coro:compatibleFullmodularGroup}, the dimension formulas $\dim_{\numC}M_{r/2}(\sltZ, v_\eta^r)=1$ for $r=0,\,1,\,2,\dots,23$ and $\dim_{\numC}S_{12}(\sltZ, v_\eta^{24})=1$. These dimension formulas have been discussed in Lemma \ref{lemm:dimone}.
\end{proof}
The main task of this section is to write $T_l \eta^r$ explicitly and to determine $c$ in Lemma \ref{lemm:etaPowerEigen}.
\begin{lemm}
\label{lemm:etaPowerTl}
Suppose that $r \in \numZ$, $l \in \numgeq{Z}{1}$ satisfying $24 \mid r(l-1)$. If in addition $r$ is odd, then $l$ is required to be a square. Then we have
\begin{multline}
\label{eq:etaPowerTl}
T_l \eta^r(\tau)= l^{-r/4}\sum_{a>0,\,a \mid l}a^{r/2} \\
\cdot\sum_{\twoscript{0 \leq b < d,\,d=l/a}{gcd(a, b, d)=1}}\psi^r_{a,b}\cdot\etp{-\frac{r}{24}(bd+3d-3)}\eta^r\left(\frac{a\tau+b}{d}\right),
\end{multline}
where $\psi^r_{a,b}=1$ if $2\mid r$, and $\psi^r_{a,b}=\legendre{-b}{\gcd(a, l/a)}$ if $2 \nmid r$.
\end{lemm}
\begin{proof}
We begin with \eqref{eq:Tlgeneral} where we set $f=\eta^r$ (and $k=r/2$, $N=1$). The first sum in the right-hand side of \eqref{eq:Tlgeneral} corresponds with the term with $a=1$ in \eqref{eq:etaPowerTl}. It remains to prove that, if $1<a \mid l$, $d=l/a$, $0 \leq b < d$ with $\gcd(a, b, d)=1$, then
\begin{equation}
\label{eq:requiredIdentity}
\left(v_\eta^r\vert_{\alpha}v_\eta^r\right)^{-1}\widetilde{\tbtMat{a}{b}{0}{d}}=\psi^r_{a,b}\cdot\etp{-\frac{r}{24}(bd+3d-3)}.
\end{equation}
Using \eqref{eq:v1v2charGeneral} with $N=1$ and $v_1=v_2=v_\eta^r$, and noting the cocycle involved is $1$, we obtain
\begin{equation}
\label{eq:v1v2etatemp}
\left(v_\eta^r\vert_{\alpha}v_\eta^r\right)^{-1}\widetilde{\tbtMat{a}{b}{0}{d}}=v_\eta^{-r}\widetilde{\tbtMat{a-(l-1)z}{b-(l-1)xz}{-(l-1)y}{d-(l-1)xy}}.
\end{equation}
We split the deduction into two cases: $2 \mid l$ and $2 \nmid l$. If $2 \mid l$, then $8 \mid r$, since we have $24 \mid r(l-1)$ by assumption. Notice that $(-b+ax)y + dz=1$ by our choice of $y$ and $z$ in \eqref{eq:matrixab0dDecom}. Hence we can always choose an odd $y$, so $2 \nmid (l-1)y$. Now we can use the first case of \eqref{eq:etaChar} to calculate the right-hand side of \eqref{eq:v1v2etatemp}. Using $\etp{r(l-1)/24}=1$, the result is simplified to $\etp{-\frac{r}{24}(bd+3d-3)}$, as required (In this case, $\psi^r_{a,b}=1$). The other case is $2 \nmid l$. This time $2 \mid (l-1)y$. Hence we use the second case of \eqref{eq:etaChar} instead. The right-hand side of \eqref{eq:v1v2etatemp} then equals $\legendre{-(l-1)y}{d-(l-1)xy}^{-r}\etp{-\frac{r}{24}(bd+3d-3)}$. If $2 \mid r$, then $\legendre{-(l-1)y}{d-(l-1)xy}^{-r}=1=\psi^r_{a,b}$, and hence \eqref{eq:requiredIdentity} holds. Otherwise if $2 \nmid r$, we shall prove that $\legendre{-(l-1)y}{d-(l-1)xy}=\legendre{-b}{\gcd(a, d)}$, from which the required identity \eqref{eq:requiredIdentity} follows. Note that since $2 \nmid r$, $l$ is required to be a square by assumption. Hence $4 \mid (l-1)y$. Thus the function $n \mapsto \legendre{-(l-1)y}{n}$ is periodic with respect to $n \in \numZ$, with a period $\abs{(l-1)y}$. Since $0 < d \mid l$, we have $\legendre{-(l-1)}{d}=\legendre{1}{d}=1$. Consequently
\begin{equation*}
\legendre{-(l-1)y}{d-(l-1)xy}=\legendre{-(l-1)y}{d}=\legendre{y}{d}=\legendre{-b+ax}{d}.
\end{equation*}
The last equality holds since $(-b+ax)y + dz=1$. Now we write $a=a_0^2a_1$ and $d=d_0^2d_1$ with $a_0,\,a_1,\,d_0,\,d_1$ being positive integers and $\gcd(a_0, a_1)=\gcd(d_0, d_1)=1$. Since $ad=l$ is a square, the integers $a_1$ and $d_1$ have the same prime divisors. Also, we have $\gcd(a,d)=\gcd(a_0, d_0)^2\gcd(a_1,d_1)$. It follows that
\begin{equation*}
\legendre{-b+ax}{d}=\legendre{-b+ax}{d_1}=\legendre{-b+ax}{\gcd(a_1, d_1)}=\legendre{-b}{\gcd(a, d)}.
\end{equation*}
Combining these calculations of Kronecker-Jacobi symbols gives $\legendre{-(l-1)y}{d-(l-1)xy}=\legendre{-b}{\gcd(a, d)}$, as required.
\end{proof}

From inserting \eqref{eq:Prn} into \eqref{eq:etaPowerTl} and rearranging the order of summation, we derive the following Fourier expansion (We emphasize that $r$ may be any integer, with $24 \mid r(l-1)$ and when $2 \nmid r$, $l$ is required to be a square.)
\begin{multline}
\label{eq:FourierTletar}
l^{r/4}T_l \eta^r(\tau)=\sum_{n \in \numZ}q^{\frac{n}{24l}} \\
\cdot\left(\sum_{0 < a \mid l}a^{r/2}\etp{-\frac{r(l/a-1)}{8}}\sum_{\twoscript{0 \leq b < l/a}{\gcd(b,a,l/a)=1}}\psi_{a,b}^r\cdot\etp{\frac{n-rl^2}{24al}b}P_r\left(\frac{n/a^2-r}{24}\right)\right).
\end{multline}
In fact, we can prove that many of the coefficients in the above Fourier series expansion must be zero.
\begin{lemm}
\label{lemm:FourierCoeffVanish}
The coefficient of the $q^{\frac{n}{24l}}$-term in the right-hand side of \eqref{eq:FourierTletar} vanishes, unless $l \mid n$ and $n/l \equiv r \bmod{24}$.
\end{lemm}
\begin{proof}
Note that by assumption, $v_\eta^r$ is $\alpha$-compatible to itself, where $\alpha=\widetilde{\tbtmat{1}{0}{0}{l}}$. Hence the function $l^{r/4}T_l \eta^r(\tau)$ is a meromorphic modular form on the group $\sltZ$, of weight $r/2$, and with the character $v_\eta^r$, using Part 4 of Proposition \ref{prop:basicPropTalpha}. Therefore
\begin{equation*}
l^{r/4}T_l\eta^r(\tau+1)=v_\eta^r\widetilde{\tbtmat{1}{1}{0}{1}}l^{r/4}T_l\eta^r(\tau)=\etp{r/24}l^{r/4}T_l\eta^r(\tau).
\end{equation*}
Using the uniqueness of Fourier or Laurent expansions, the Fourier expansion of $l^{r/4}T_l \eta^r(\tau)$ at the cusp $\rmi\infty$ must be an infinite sum of terms of the form $c_m\cdot q^{m/24}$ with $m\in r+24\numZ$ and $c_m \in \numC$. Again using the uniqueness property, the $q^{\frac{n}{24l}}$-term in the right-hand side of \eqref{eq:FourierTletar} must vanish, unless $l \mid n$ and $n/l \equiv r \bmod{24}$.
\end{proof}
To simplify the expression \eqref{eq:FourierTletar}, we need to evaluate the sum over $b$, for which the following simple relation plays a necessary role.
\begin{lemm}
\label{lemm:24divd2}
Let $r$ be an integer, and $l$ a positive integer, such that $24 \mid r(l-1)$. Then for any positive divisor $d$ of $l$, we have $24 \mid r(d^2-1)$.
\end{lemm}
\begin{proof}
Note that $(d_1d_2)^2-1=d_1^2(d_2^2-1)+d_1^2-1$. So $24 \mid r(d_1^2-1)$ and $24 \mid r(d_2^2-1)$ together would imply that $24 \mid r((d_1d_2)^2-1)$. This means we need only prove the desired relation for $d$ being a prime power $p^\beta$ ($\beta>0$). It is easy to see that if $p \geq 5$, then $24 \mid p^2-1$, and hence $24 \mid r(p^{2\beta}-1)$. If $p=2$, then $8 \mid r$, since in this case $2 \nmid (l-1)$. On the other hand, we have $3 \mid 2^{2\beta}-1$. Consequently $24 \mid r(2^{2\beta}-1)$. The remaining case is $p=3$, in which $3 \mid r$ and $8 \mid 3^{2\beta}-1$. Hence $24 \mid r(3^{2\beta}-1)$.
\end{proof}

Now we can obtain a computable formula for the Fourier coefficients of $T_l \eta^r(\tau)$. First is the case when $r$ is even, which is of course easier than the case when $r$ is odd. Recall that the \textit{M\"obius function} $\mu$ on $\numgeq{Z}{1}$ maps a product of $j$ distinct primes to $(-1)^j$, $\mu(1)=1$, and maps other positive integers to $0$.
\begin{lemm}
\label{lemm:FourierCoeffReven}
Let $r$ be an even integer, and $l$ a positive integer, such that $24 \mid r(l-1)$. Then the Fourier coefficient of the $q^{n/(24l)}$-term in the right-hand side of \eqref{eq:FourierTletar} is
\begin{equation}
\label{eq:fourierCoeffTletarEven}
\sum_{0 < a \mid l}a^{r/2}(-1)^{r(d-1)/4}\sum_{t}\mu(t)\frac{d}{t}P_r\left(\frac{n/a^2-r}{24}\right),
\end{equation}
where $d=l/a$, and $t$ ranges over positive integers divisible by $24al/\gcd(24al, n-rl^2)$ and dividing $\gcd(a,d)$ in the inner sum. Moreover, the quantity $r(d-1)/4$ occurring in \eqref{eq:fourierCoeffTletarEven} is always an integer. In particular, if $l$ is square-free, then \eqref{eq:fourierCoeffTletarEven} may be further simplified to
\begin{equation}
\label{eq:fourierCoeffTletarEvenSquareFree}
l\cdot\sum_{\twoscript{0 < a \mid l}{24al \mid n-rl^2}}a^{r/2-1}(-1)^{r(d-1)/4}P_r\left(\frac{n/a^2-r}{24}\right).
\end{equation}
\end{lemm}
\begin{proof}
We only prove \eqref{eq:fourierCoeffTletarEven}, since under the special conditions stated in this lemma the formula \eqref{eq:fourierCoeffTletarEvenSquareFree} follows immediately from \eqref{eq:fourierCoeffTletarEven}. First we prove that $4 \mid r(d-1)$. If $4 \mid r$, then definitely $4 \mid r(d-1)$. Otherwise $r \equiv 2\bmod{4}$. Hence $4 \mid l-1$, using the assumption $24 \mid r(l-1)$. Thus $2 \mid d-1$, and $4 \mid r(d-1)$, as desired. Next we shall prove that, for $a \mid l$ with $\frac{n/a^2-r}{24}$ being an integer, the following holds:
\begin{equation}
\label{eq:needToProveSumMu}
\sum_{\twoscript{0 \leq b < d}{\gcd(b,a,d)=1}}\etp{\frac{n-rl^2}{24al}b}=\sum_{t \in \mathfrak{A}}\mu(t)\frac{d}{t},
\end{equation}
where $\mathfrak{A}=\{t \in \numgeq{Z}{1} \colon t \mid \gcd(a,d),\, 24al/\gcd(24al, n-rl^2)\mid t\}$. From this and \eqref{eq:FourierTletar} the aimed assertion \eqref{eq:fourierCoeffTletarEven} follows. Using that $\sum_{t \mid n}\mu(t)$ is equal to $1$ if $n=1$ and to $0$ if $n>1$, we have
\begin{equation}
\label{eq:changeSumMu}
\sum_{\twoscript{0 \leq b < d}{\gcd(b,a,d)=1}}\etp{\frac{n-rl^2}{24al}b} = \sum_{t \mid \gcd(a, d)}\mu(t)\sum_{0 \leq b < d/t}\etp{\frac{n-rl^2}{24al}tb}.
\end{equation}
Note that $\frac{n-rl^2}{24a^2} = \frac{n/a^2-r}{24}-\frac{r(d^2-1)}{24}$ is an integer since $\frac{n/a^2-r}{24}$ is by assumption and $\frac{r(d^2-1)}{24}$ is by Lemma \ref{lemm:24divd2}. Now if $24al/\gcd(24al, n-rl^2)\nmid t$, then $\frac{d}{t} \nmid \frac{n-rl^2}{24a^2}$ and hence
\begin{equation*}
\sum_{0 \leq b < d/t}\etp{\frac{n-rl^2}{24al}tb} = \sum_{0 \leq b < d/t}\etp{\frac{n-rl^2}{24a^2}\frac{b}{d/t}} = 0;
\end{equation*}
while if $24al/\gcd(24al, n-rl^2)\mid t$ then $\sum_{0 \leq b < d/t}\etp{\frac{n-rl^2}{24al}tb}=d/t$. Inserting these into \eqref{eq:changeSumMu} gives \eqref{eq:needToProveSumMu}, as required.
\end{proof}

We now deal with the case $2 \nmid r$. The major difficulty is an evaluation of a Gauss character sum, for which we introduce some necessary notations.
\begin{deff}
\label{def:rad}
Let $m$ be a positive integer, $p$ be a prime, and $\beta$ be a non-negative integer. By $p^{\beta}\parallel m$, we mean that $p^{\beta}\mid m$, but $p^{\beta+1}\nmid m$. Moreover, we define $\rad(m)$, $\radE(m)$, $\radO(m)$, $\radp(m)$, $\irad(m)$ and $\iradp(m)$ by the following formulas:
\begin{align*}
\radE(m) &= \prod_{\twoscript{p^{\beta}\parallel m}{2 \mid \beta}}p  \qquad &\radO(m) &= \prod_{\twoscript{p^{\beta}\parallel m}{2 \nmid \beta}}p \\
\rad(m) &= \radE(m)\radO(m) \qquad &\radp(m) &= \radE(m)^2\radO(m) \\
\irad(m) &= m/\rad(m) \qquad &\iradp(m)&=m/\radp(m).
\end{align*}
The subscript $E$ stands for \textit{even}, $O$ stands for \textit{odd}, and irad stands for \textit{inverse radical}, which are the reasons of the names. We also define $f(a, b,\dots)$ as $f(\gcd(a,b,\dots))$ with $f$ being one of the above six functions.
\end{deff}
We remark that $\radE(m)$ contains all primes that only appear in the square-full part of $m$ and $\radO(m)$ contains primes that appear in the square-free part of $m$ (and possibly also the square-full part).
\begin{lemm}
\label{lemm:aGaussSum}
Let $m$ be a positive odd number, and $t$ be an integer. Let $u$ be the number of primes dividing $\radO(m)$ that are congruent to $3$ modulo $4$. If $\irad(m) \nmid t$, then $\sum_{b=0}^{m-1}\legendre{b}{m}\etp{\frac{tb}{m}}=0$.
Otherwise, if $\irad(m) \mid t$, we have
\begin{multline}
\label{eq:aGaussSum}
\sum_{b=0}^{m-1}\legendre{b}{m}\etp{\frac{tb}{m}} \\
=\rmi^{u^2}\frac{m}{\sqrt{\radp(m)}}\legendre{t/\iradp(m)}{\radO(m)}\prod_{p \mid \radE(m)}\left(p-1-p\legendre{t/\irad(m)}{p}^2\right).
\end{multline}
\end{lemm}
\begin{proof}
The case $m=1$ is trivial, so we assume that $m > 1$, and write $m=p_1^{\beta_1}\dots p_s^{\beta_s}$, where $s \geq 1$, $\beta_j \geq 1$ and $p_j$ are distinct odd primes for $1 \leq j \leq s$.
Put $m_j=m/p_j^{\beta_j}$, and let $M_j$ be any integer such that $m_jM_j \equiv 1 \bmod p_j^{\beta_j}$, for $1 \leq j \leq s$. Then the map
\begin{align*}
\numZ/p_1^{\beta_1}\numZ \times \numZ/p_2^{\beta_2}\numZ \times\dots \numZ/p_s^{\beta_s}\numZ &\rightarrow \numZ/m\numZ \\
(b_1+p_1^{\beta_1}\numZ, b_2+p_2^{\beta_2}\numZ, \dots, b_s+p_s^{\beta_s}\numZ) &\mapsto (\sum_{1 \leq j \leq s}b_jm_jM_j+m\numZ)
\end{align*}
is a ring isomorphism. Therefore, using a change of variables $b=\sum_{1 \leq j \leq s}b_jm_jM_j$, we have
\begin{equation}
\label{eq:GaussSumDecomp}
\sum_{b=0}^{m-1}\legendre{b}{m}\etp{\frac{tb}{m}}=\prod_{1 \leq j \leq s}\sum_{b_j=0}^{p_j^{\beta_j}-1}\legendre{b_j}{p_j^{\beta_j}}\etp{\frac{tM_jb_j}{p_j^{\beta_j}}}.
\end{equation}
We need to evaluate $\sum_{b_j}\legendre{b_j}{p_j^{\beta_j}}\etp{tM_jb_j/p_j^{\beta_j}}$. Apply a change of variables $b_j=b_{j,0}+p_j\cdot b_{j,1}$ with $0 \leq b_{j,0} < p_j$ and $0 \leq b_{j,1} < p_j^{\beta_j-1}$. There are two cases: $2 \nmid \beta_j$, and $2 \mid \beta_j$.

Consider the case $2 \nmid \beta_j$. Note that $\legendre{b_j}{p_j^{\beta_j}}=\legendre{b_j}{p_j}$ in this case; we obtain
\begin{equation}
\label{eq:b0b1decomp}
\sum_{b_j=0}^{p_j^{\beta_j}-1}\legendre{b_j}{p_j^{\beta_j}}\etp{\frac{tM_jb_j}{p_j^{\beta_j}}}=\sum_{b_{j,0}=0}^{p_j-1}\legendre{b_{j,0}}{p_j}\etp{\frac{tM_jb_{j,0}}{p_j^{\beta_j}}}\sum_{b_{j,1}=0}^{p_j^{\beta_j-1}-1}\etp{\frac{tM_jb_{j,1}}{p_j^{\beta_j-1}}}
\end{equation}
If $p_j^{\beta_j-1} \nmid t$, then $p_j^{\beta_j-1} \nmid M_jt$, since $\gcd(p_j, M_j)=1$. Hence the sum in the above equality vanishes. Otherwise, $p_j^{\beta_j-1} \mid t$, so
\begin{equation}
\label{eq:b0b1decompContinue}
\sum_{b_j=0}^{p_j^{\beta_j}-1}\legendre{b_j}{p_j^{\beta_j}}\etp{\frac{tM_jb_j}{p_j^{\beta_j}}}=p_j^{\beta_j-1}\legendre{(t/p_j^{\beta_j-1})M_j}{p_j}\sum_{b_{j,0}=0}^{p_j-1}\legendre{b_{j,0}}{p_j}\etp{\frac{b_{j,0}}{p_j}}.
\end{equation}
In the above equality, we use the fact that the classical Gauss sum $\sum_{b=0}^{p-1}\legendre{b}{p}\etp{\frac{xb}{p}}$ is separable for every $x \in \numZ$. Recall Gauss' classical result $\sum_{b=0}^{p-1}\legendre{b}{p}\etp{\frac{b}{p}}=\varepsilon_p\sqrt{p}$, where $\varepsilon_p=1$, if $p \equiv 1 \bmod 4$, and $\varepsilon_p=\rmi$, if $p \equiv 3 \bmod 4$ (consult \cite[Section 4.7]{Nat00} for a detailed proof). Thus
\begin{equation}
\label{eq:GaussSumPrimepowerOdd}
\sum_{b_j=0}^{p_j^{\beta_j}-1}\legendre{b_j}{p_j^{\beta_j}}\etp{\frac{tM_jb_j}{p_j^{\beta_j}}}=p_j^{\beta_j-1}\sqrt{p_j}\varepsilon_p\legendre{tM_j/p_j^{\beta_j-1}}{p_j}.
\end{equation}

Now consider the case $2 \mid \beta_j$. We have
\begin{equation*}
\sum_{b_j=0}^{p_j^{\beta_j}-1}\legendre{b_j}{p_j^{\beta_j}}\etp{\frac{tM_jb_j}{p_j^{\beta_j}}}=\sum_{0 < b_{j,0} < p_j}\etp{\frac{tM_jb_{j, 0}}{p_j^{\beta_j}}}\sum_{0 \leq b_{j, 1} < p_j^{\beta_j-1}}\etp{\frac{tM_jb_{j,1}}{p_j^{\beta_j-1}}}.
\end{equation*}
If $p_j^{\beta_j-1} \nmid t$, then the sum over $b_j$ vanishes, since the inner sum in the right-hand side vanishes, and if $p_j^{\beta_j-1} \mid t$ then
\begin{align}
\sum_{b_j=0}^{p_j^{\beta_j}-1}\legendre{b_j}{p_j^{\beta_j}}\etp{\frac{tM_jb_j}{p_j^{\beta_j}}}&=p_j^{\beta_j-1}\sum_{0 < b_{j,0} < p_j}\etp{\frac{tM_j/p_j^{\beta_j-1}}{p_j}b_{j, 0}}\notag\\
&=p_j^{\beta_j-1}\left(p_j-1-p_j\legendre{tM_j/p_j^{\beta_j-1}}{p_j}^2\right). \label{eq:GaussSumPrimepowerEven}
\end{align}
The last equality can be deduced by splitting the calculation into two cases: $p_j \mid tM_j/p_j^{\beta_j-1}$ and $p_j \nmid tM_j/p_j^{\beta_j-1}$.

We turn back to the evaluation of the quantity $\sum_{b=0}^{m-1}\legendre{b}{m}\etp{\frac{tb}{m}}$. If $\irad(m) \nmid t$, then $p_j^{\beta_j-1} \nmid t$ for some $j$. Therefore $\sum_{b=0}^{m-1}\legendre{b}{m}\etp{\frac{tb}{m}}=0$ by \eqref{eq:GaussSumDecomp} and what we have proved for the prime power cases. Otherwise, we should consider the case $\irad(m) \mid t$. Inserting \eqref{eq:GaussSumPrimepowerOdd} and \eqref{eq:GaussSumPrimepowerEven} into \eqref{eq:GaussSumDecomp} gives
\begin{multline}
\label{eq:GaussSumMidFormula}
\sum_{b=0}^{m-1}\legendre{b}{m}\etp{\frac{tb}{m}}\\
= \rmi^u\frac{m}{\sqrt{\radp(m)}}\prod_{p_j \mid \radO(m)}\legendre{tM_j/p_j^{\beta_j-1}}{p_j}\prod_{p_j \mid \radE(m)}\left(p_j-1-p_j\legendre{tM_j/p_j^{\beta_j-1}}{p_j}^2\right).
\end{multline}
Note that
\begin{align*}
\prod_{p_j \mid \radO(m)}\legendre{tM_j/p_j^{\beta_j-1}}{p_j} &=\legendre{t/\iradp(m)}{\radO(m)}\prod_{\twoscript{p_j \mid \radO(m)}{p_{j'} \mid \radO(m),\,j \neq j'}}\legendre{p_{j'}}{p_j} \\
&=\legendre{t/\iradp(m)}{\radO(m)}(-1)^{u(u-1)/2},
\end{align*}
\begin{align*}
\prod_{p_j \mid \radE(m)}\left(p_j-1-p_j\legendre{tM_j/p_j^{\beta_j-1}}{p_j}^2\right) &=\prod_{p_j \mid \radE(m)}\left(p_j-1-p_j\legendre{t/p_j^{\beta_j-1}}{p_j}^2\right) \\
&=\prod_{p_j \mid \radE(m)}\left(p_j-1-p_j\legendre{t/\irad(m)}{p_j}^2\right).
\end{align*}

Inserting these formulas into \eqref{eq:GaussSumMidFormula} gives \eqref{eq:aGaussSum}, as desired.
\end{proof}

We are now at a position to give the following formula, which is an analogue of Lemma \ref{lemm:FourierCoeffReven} for the case $2 \nmid r$.
\begin{lemm}
\label{lemm:FourierCoeffRodd}
Let $r$ be an odd integer, and $l$ a positive perfect square, such that $24 \mid r(l-1)$. Let $n$ be an integer. Then the Fourier coefficient of the $q^{n/24}$-term in the right-hand side of \eqref{eq:FourierTletar} is
\begin{multline}
\label{eq:FourierCoeffRodd}
\sum_{\twoscript{0 < a \mid l}{24a \mid \rad(a,d)(n-rl)}}\legendre{2\cdot(-1)^{(r-1)/2}}{d}\legendre{\radE(a,d)}{\radO(a,d)}\frac{d\cdot a^{r/2}}{\sqrt{\radp(a, d)}}P_r\left(\frac{nl/a^2-r}{24}\right) \\
\times \legendre{\rad(a,d)(n-rl)/(24a)}{\radO(a, d)}\prod_{p \mid \radE(a, d)}\left(p-1-p\legendre{\rad(a,d)(n-rl)/(24a)}{p}^2\right),
\end{multline}
where $d$ stands for $l/a$. In particular, if $l$ is the square of a square-free integer, then \eqref{eq:FourierCoeffRodd} becomes
\begin{multline}
\label{eq:FourierCoeffRoddSpecial}
\sum_{\twoscript{0 < a \mid l}{24a \mid \rad(a,d)(n-rl)}}\legendre{2\cdot(-1)^{(r-1)/2}}{d}\frac{d\cdot a^{r/2}}{\sqrt{\rad(a, d)}}P_r\left(\frac{nl/a^2-r}{24}\right) \\
\times \legendre{\rad(a,d)(n-rl)/(24a)}{\rad(a, d)}.
\end{multline}
\end{lemm}
\begin{proof}
Assume that $n \equiv r \bmod 24$. By \eqref{eq:FourierTletar}, the coefficient of $q^{n/24}$-term is
\begin{equation}
\label{eq:lemmFourierCoeffRoddeq1}
\sum_{0 < a \mid l}a^{r/2}\etp{-\frac{r(d-1)}{8}}\sum_{0 \leq b < d}\legendre{-b}{\gcd(a,d)}\cdot\etp{\frac{n-rl}{24a}b}P_r\left(\frac{nl/a^2-r}{24}\right).
\end{equation}
We shall evaluate $\sum_{b}\legendre{-b}{\gcd(a,d)}\cdot\etp{\frac{n-rl}{24a}b}$ occurring in the above sum. We may assume that $24 \mid nl/a^2-r$, since the factor $P_r\left(\frac{nl/a^2-r}{24}\right)$ otherwise vanishes. Since $l$ is a square and $24 \mid r(l-1)$, we can conclude that $gcd(a, d)$ is odd. Moreover, we have $\frac{n-rl}{24a}\cdot a \in \numZ$, since $\frac{n-rl}{24}=\frac{n-r}{24}-\frac{r(l-1)}{24}\in \numZ$. Also, $\frac{n-rl}{24a}\cdot d \in \numZ$, because of the assumption $24 \mid nl/a^2-r$, the fact $\frac{nl/a^2-rd^2}{24}=\frac{nl/a^2-r}{24}-\frac{r(d^2-1)}{24}$, and Lemma \ref{lemm:24divd2}. It follows that $\frac{n-rl}{24a} \in \gcd(a, d)^{-1}\numZ$. By setting $m=\gcd(a, d)$ and $t=\frac{\gcd(a,d)(n-rl)}{24a}$ in Lemma \ref{lemm:aGaussSum}, we obtain that, if $24a \mid \rad(a,d)(n-rl)$, then
\begin{multline}
\label{eq:lemmFourierCoeffRoddGaussSum}
 \sum_{0 \leq b < d}\legendre{-b}{\gcd(a,d)}\cdot\etp{\frac{n-rl}{24a}b}=\legendre{-1}{\gcd(a,d)}\frac{d}{\gcd(a,d)}\rmi^{u_{a,d}^2}\frac{\gcd(a,d)}{\sqrt{\radp(a, d)}} \\
 \times \legendre{\gcd(a,d)\frac{n-rl}{24a}/\iradp(a,d)}{\radO(a,d)}\prod_{p \mid \radE(a, d)}\left(p-1-p\legendre{\gcd(a,d)\frac{n-rl}{24a}/\irad(a,d)}{p}^2\right),
\end{multline}
where $u_{a, d}$ denotes the number of primes dividing $\radO(a, d)$ that are congruent to $3$ modulo $4$. On the other hand, if $24a \nmid \rad(a,d)(n-rl)$, then $\irad(m) \nmid t$, and hence by Lemma \ref{lemm:aGaussSum}
\begin{equation*}
\sum_{0 \leq b < d}\legendre{-b}{\gcd(a,d)}\cdot\etp{\frac{n-rl}{24a}b}=\sum_{0 \leq b < d}\legendre{-b}{m}\cdot\etp{\frac{tb}{m}}=0.
\end{equation*}
Taking into account this relation and \eqref{eq:lemmFourierCoeffRoddGaussSum}, the quantity \eqref{eq:lemmFourierCoeffRoddeq1} can be simplified to
\begin{multline}
\sum_{\twoscript{0 < a \mid l}{24a \mid \rad(a,d)(n-rl)}}\legendre{\radE(a,d)}{\radO(a,d)}\frac{d\cdot a^{r/2}}{\sqrt{\radp(a,d)}} P_r\left(\frac{nl/a^2-r}{24}\right) \\
\times \etp{-\frac{r(d-1)}{8}}\rmi^{u_{a, d}^2}\legendre{-1}{\gcd(a, d)} \\
\times  \legendre{\rad(a,d)(n-rl)/(24a)}{\radO(a, d)}\prod_{p \mid \radE(a, d)}\left(p-1-p\legendre{\rad(a,d)(n-rl)/(24a)}{p}^2\right).
\end{multline}
To acquire \eqref{eq:FourierCoeffRodd}, based on the above formula, it remains to prove that
\begin{equation}
\label{eq:lemmFourierCoeffRodduad}
\etp{-\frac{r(d-1)}{8}}\rmi^{u_{a, d}^2}\legendre{-1}{\gcd(a, d)}=\legendre{2\cdot(-1)^{(r-1)/2}}{d}.
\end{equation}
Since $l$ is an odd square, $2\mid u_{a,d}$ if $a\equiv d \equiv 1 \bmod 4$, and $2\nmid u_{a,d}$ if $a\equiv d \equiv 3 \bmod 4$. It follows that $\rmi^{u_{a,d}^2}=\etp{\frac{1-d}{8}}\legendre{2}{d}$, and $\legendre{-1}{\gcd(a, d)}=\legendre{-1}{d}$, from which \eqref{eq:lemmFourierCoeffRodduad} follows.

If $n \not\equiv r \bmod 24$, then the quantity \eqref{eq:FourierCoeffRodd} vanishes since there is no $a$ satisfying $24a \mid \rad(a,d)(n-rl)$, while Lemma \ref{lemm:FourierCoeffVanish} tells us the Fourier coefficient of the $q^{n/24}$-term in the right-hand side of \eqref{eq:FourierTletar} is really zero.
\end{proof}

The main theorem of this section now follows from all the above lemmas. We emphasize that $d$ always stands for $l/a$ in the following theorem, and that the function $P_r(n)$ is defined in \eqref{eq:Prn} and the sentence that follows it. We introduce a symbol $R(n;r,l)$ as follows. If $r$ is even, then
\begin{equation*}
R(n;r,l)=\sum_{0 < a \mid l}a^{r/2}(-1)^{r(d-1)/4}\sum_{\twoscript{0 < t \mid \gcd(a,d)}{24a \mid \gcd(24a,n-rl)t}}\mu(t)\frac{d}{t}P_r\left(\frac{nl/a^2-r}{24}\right),
\end{equation*}
and if $r$ is odd, then $R(n;r,l)$ denotes the quantity \eqref{eq:FourierCoeffRodd}.
\begin{thm}
\label{thm:mainThmEtar}
Suppose that $r \in \{0,\,1,\,2,\dots,24\}$, and that $l \in \numgeq{Z}{1}$ satisfying $24 \mid r(l-1)$. If $r$ is odd, we require that $l$ is a perfect square. Then for every $n \in \numZ$, we have
\begin{equation}
\label{eq:mainThmEtarEvenOdd}
R(n;r,l)=R(r;r,l)P_r\left(\frac{n-r}{24}\right).
\end{equation}
\end{thm}
\begin{proof}
Note that the left-hand side of the formula \eqref{eq:mainThmEtarEvenOdd} is the Fourier coefficient of the $q^{n/24}$-term of $l^{r/4}T_l \eta^r(\tau)$ by Lemma \ref{lemm:FourierCoeffReven} and Lemma \ref{lemm:FourierCoeffRodd}. On the other hand, the right-hand side of the formula \eqref{eq:mainThmEtarEvenOdd} is the Fourier coefficient of the $q^{n/24}$-term of $R(r;r,l)\eta^r(\tau)$. Thus, by Lemma \ref{lemm:etaPowerEigen}, the desired formula follows, up to a constant factor. By inserting $r$ for $n$ in the formula, we find that the constant factor is $1$.
\end{proof}

\begin{examp}
Some special cases of Theorem \ref{thm:mainThmEtar} have been previously discovered by Newman. The case when $r$ is even and $l$ is a prime power was proved in \cite{New55} and \cite{New56}. The case when $r$ is odd and $l$ is the square of a prime was proved in \cite{New58}. Our method is quite different from Newman's, and it seems that the general formula \eqref{eq:mainThmEtarEvenOdd} for all admissible $l$ can not follow from Newman's results in which $l$ are prime powers directly. Newman's results were also proved in \cite[Section 4]{vL57} using general Hecke operators as introduced by Wohlfahrt, while van Lint did not describe $T_l$ when $l$ is not a prime power explicitly.
\end{examp}

\section{Express $F_{r, p^{\beta}}(\tau)$ as Linear Combinations of Eta-quotients}
\label{sec:Express}
We now turn to the third application of the theory developed in Section \ref{sec:Operators}. We will study functions of the form $T_{p^{\beta}} \eta^r$ but on smaller groups $\Gzt{p}$, instead of $\sltZ$. These functions can be explicitly written as
\begin{equation}
\label{eq:expressSecMainFun}
F_{r, p^{\beta}}(\tau)=\sum_{m \in \frac{p^{\ob}r}{24}+\numZ}P_r\left(p^{\beta}m-\frac{r}{24}\right)q^m=q^{p^{\ob}r/24}\sum_{n \in\numZ}P_r\left(p^{\beta}n+\frac{r(p^{\beta+\ob}-1)}{24}\right)q^n,
\end{equation}
where $r$ is an integer, $p$ a prime satisfying $24 \mid r(p^2-1)$, $\beta$ is a positive integer and $\ob$ is equal to $1$ ($0$ respectively) if $\beta$ is odd (even respectively). The aim is to give a criterion on whether we can express these functions as linear combinations of eta-quotients of level $p$, and an algorithm for calculating the coefficients. First, we shall prove that these functions are weakly holomorphic modular forms and are actually of the form $T_{p^{\beta}} \eta^r$ on the group $\Gzt{p}$. We remind the reader that the operator $T_l$ depends not only on $l$, but also on the group, the characters and the weight. See the paragraph after Definition \ref{deff:doubleCosetOperators}. Therefore the operator $T_l$ may differ from that one with the same name used in Section \ref{sec:eta-powers}.

For the property of being weakly holomorphic modular forms, we consider a wider class of functions:
\begin{lemm}
\label{lemm:TletarlevelN}
Suppose $l$ and $N$ are positive integers such that $\rad(l)\mid\rad(N)$, and $\mathbf{r'}=(r_1', r_2', \dots)$ is an integral sequence such that $r_n' \neq 0$ implies $n \mid N$ for $n \in \numgeq{Z}{1}$. Let $r=\sum_{n \mid N}r_n'$. If the following conditions hold:
\begin{align}
\sum_{n \mid N}(nl-1)r_n' &\equiv 0 \pmod{24}, \label{eq:assume1TletarlevelN}\\
\sum_{n \mid N}\frac{N}{n}(nl-1)r_n' &\equiv 0 \pmod{24}, \label{eq:assume2TletarlevelN}\\
l^{\abs r}\prod_{\twoscript{n \mid N}{2 \nmid r_n'}}n &\text{ is a square, } \label{eq:assume3TletarlevelN}
\end{align}
then the function
\begin{equation}
\label{eq:TletarlevelNFunction}
\sum_{m \in \frac{1}{24}\sum_{n \mid N}nr_n'+\numZ}P_r\left(lm-\frac{r}{24}\right)q^m
\end{equation}
is in the space $\MFormW{r/2}{N}{v_{\mathbf{r'}}}$ (see \eqref{eq:etaQuoChar} for the definition of $v_{\mathbf{r'}}$). If in addition, $r>0$, then it belongs to $\CForm{r/2}{N}{v_{\mathbf{r'}}}$.
\end{lemm}
\begin{proof}
Set $\mathbf{r}=(r,0, \dots)$, that is, $r_1=r$ but $r_n = 0$ for any $n>1$. Then Theorem \ref{thm:charCompatibleEtaDirichlet} and the three assumptions \eqref{eq:assume1TletarlevelN}, \eqref{eq:assume2TletarlevelN} and \eqref{eq:assume3TletarlevelN} together imply that $v_{\mathbf{r}}$ and $v_{\mathbf{r'}}$ are $\widetilde{\tbtmat{1}{0}{0}{l}}$-compatible. Therefore, regarding $\eta^r(\tau)$ as a weakly holomorphic modular form on the group $\Gzt{N}$, we have $T_l \eta^r(\tau) \in \MFormW{r/2}{N}{v_{\mathbf{r'}}}$ by Part 5 of Proposition \ref{prop:basicPropTalpha}. Since $\rad(l)\mid\rad(N)$, the function \eqref{eq:TletarlevelNFunction} is a constant times  $T_l \eta^r(\tau)$, by \eqref{eq:radlcaseFourier2} and the fact that $v_{\mathbf{r'}}(\widetilde{T})=\etp{\frac{1}{24}\sum_{n \mid N}nr_n'}$. Thus \eqref{eq:TletarlevelNFunction} belongs to $\MFormW{r/2}{N}{v_{\mathbf{r'}}}$. The last assertion follows from which we have proved, formula \eqref{eq:etaOrder} and the definition of cusp forms.
\end{proof}

Since the main purpose of this section is to investigate functions $F_{r, p^{\beta}}$, we shall only consider the case that $l$ is a prime power. Although it is possible to also study the case when $N$ is a prime power we will only consider the case when $N$ is a prime in this paper. We will return to study the more general case in a forthcoming paper.
\begin{lemm}
\label{lemm:Tpbetarlevelp}
Suppose $p$ is a prime, $\beta$ is a positive integer, and $r$ is an integer. Set $N=p$ and $l=p^{\beta}$ in Lemma \ref{lemm:TletarlevelN}. Then there exists some $\mathbf{r'}$ satisfying \eqref{eq:assume1TletarlevelN}, \eqref{eq:assume2TletarlevelN} and \eqref{eq:assume3TletarlevelN}, such that $r_j'=0$ unless $j=1$ or $p$, and $r_1'+r_p'=r$, if and only if $24 \mid r(p^2-1)$. In this case, all $\mathbf{r'}$ satisfying these conditions, which are determined by $r_1'$ and $r_p'$, are given exactly by the following equations. If $\beta$ is odd, then
\begin{equation}
\label{eq:TpbetarlevelpSolutionbodd}
r_1'=-\frac{24}{\gcd(12, p-1)}y \qquad r_p'=r-r_1',
\end{equation}
with $y \in \numZ$. If $\beta$ is even, then
\begin{equation}
\label{eq:TpbetarlevelpSolutionbeven}
r_1'=r - \frac{24}{\gcd(12, p-1)}y \qquad r_p'=r-r_1',
\end{equation}
with $y \in \numZ$.
\end{lemm}
\begin{proof}
Note that in the case $N=p$ and $l=p^{\beta}$, \eqref{eq:assume1TletarlevelN}, \eqref{eq:assume2TletarlevelN} and \eqref{eq:assume3TletarlevelN} are equivalent to the following system of congruences
\begin{align}
p^{\beta}(p-1)r_1'&\equiv r(p^{\beta+1}-1) \pmod{24} \label{eq:systemEquiva}\\
(p-1)r_1' &\equiv r(p^{\beta+1}-1) \pmod{24} \notag\\
r_1' &\equiv (\beta-1)r \pmod{2} \notag\\
r_p' &= r-r_1'.\notag
\end{align}
We will show that this system is solvable for $r_1',\,r_p'\in\numZ$ if and only if $24 \mid r(p^2-1)$, and that when solvable, the solutions of this system are exactly those given by \eqref{eq:TpbetarlevelpSolutionbodd} when $\beta$ is odd, or by \eqref{eq:TpbetarlevelpSolutionbeven} when $\beta$ is even. We break the proof into cases.

\textit{Case 1.} $\beta$ is odd and $p\geq5$. Then $24 \mid p^{\beta+1}-1$, so the system \eqref{eq:systemEquiva} is equivalent to $\frac{24}{\gcd(12, p-1)} \mid r_1'$ and $r_p'=r-r_1'$ since the least common multiple of $\frac{24}{\gcd(24, p-1)}$ and $2$ is $\frac{24}{\gcd(12, p-1)}$. Also in this case the statement $24 \mid r(p^2-1)$ is true. Therefore, the desired statements follow.

\textit{Case 2.} $\beta$ is odd and $p=2$. The system \eqref{eq:systemEquiva} becomes $2^{\beta}r_1'\equiv r(2^{\beta+1}-1) \bmod{24}$, $r_1'\equiv r(2^{\beta+1}-1) \bmod{24}$ and $r_1' \equiv 0\bmod{2}$. If it has a solution ($r_1'$, $r_p'$), then by subtracting the second congruence from the first, we obtain $24 \mid (2^{\beta}-1)r_1'$, so $24 \mid r_1'$. Hence $24 \mid r(2^2-1)$, and the solutions are exactly those given by \eqref{eq:TpbetarlevelpSolutionbodd}. Conversely, suppose $24 \mid r(2^2-1)$. Then the system is solvable.

We leave the remaining cases to the reader.
%
%
%
%

\end{proof}
We introduce some notations. Let $p,\,\beta,\,r$ be as in Lemma \ref{lemm:Tpbetarlevelp}. For any $n\in \numZ$ and $y \in \numZ$, we define the numbers $a_n^{r,p^\beta, y}$ by the following formula:
\begin{equation}
\label{eq:etaQuotientsUsedToExpressTletar}
\eta^r(p^{\ob}\tau)\left(\frac{\eta(p\tau)}{\eta(\tau)}\right)^{\frac{24}{\gcd(12,p-1)}y}=q^{\frac{p^{\ob}r}{24}+\frac{p-1}{\gcd(12,p-1)}y}\sum_{n \in \numgeq{Z}{0}}a_n^{r,p^\beta, y}q^n.
\end{equation}
The function that occurs in the above formula is an eta-quotient, whose character is $v_{\mathbf{r'}}$, where the only nonzero components $r_1'$ and $r_p'$ are just \eqref{eq:TpbetarlevelpSolutionbodd} or \eqref{eq:TpbetarlevelpSolutionbeven}, according to whether $\beta$ is odd or even. Note that the characters $v_\mathbf{r'}$ are in fact independent of $y$. In the case when $F_{r, p^{\beta}}$ is not identically zero this follows immediately by observing that it must have a unique character. The explicit formulas \eqref{eq:etaChar} and \eqref{eq:etaQuoChar} can be used to show this independence in general and we leave the details to the reader.


The sequences $a_n^{r,p^\beta, y}$ can also be given by recursive relations, as the following lemma shows.
\begin{lemm}
\label{lemm:crpynRecur}
Denote the sum of positive divisors of a positive integer $k$ by $\sigma(k)$, and set $\sigma(k)=0$ if k is not a positive integer. Let $p$ be a prime, $\beta$ be a positive integer, and $r$, $y$ be integers. Then we have $a_0^{r,p^\beta, y}=1$ and
\begin{multline}
a_n^{r,p^\beta, y}=-\frac{1}{n}\\
\cdot\sum_{1\leq k \leq n}\left(rp^{\ob}\sigma\left(\frac{k}{p^{\ob}}\right)+\frac{24}{\gcd(12,p-1)}y\left(p\sigma\left(\frac{k}{p}\right)-\sigma(k)\right)\right)a_{n-k}^{r,p^\beta, y},
\end{multline}
for any $n\in\numgeq{Z}{1}$.
\end{lemm}
\begin{proof}
Insert \eqref{eq:defEta} into the left-hand side of \eqref{eq:etaQuotientsUsedToExpressTletar}, and take the logarithmic derivative of both sides with respect to $q$. The desired identity follows using the well-known identity
\begin{equation*}
\sum_{n \geq 1}\frac{nq^n}{1-q^n}=\sum_{n \geq 0}\sigma(n)q^n.
\end{equation*}
\end{proof}
\begin{examp}
Put $Y=24y/\gcd(12,p-1)$, if $\beta$ is odd, or $Y=24y/\gcd(12,p-1)-r$, if $\beta$ is even. For $p \geq 5$, we have
\begin{align}
a_1^{r,p^\beta, y} &= Y, \\
a_2^{r,p^\beta, y} &= \frac{3}{2}Y+\frac{1}{2}Y^2, \\
a_3^{r,p^\beta, y} &= \frac{4}{3}Y+\frac{3}{2}Y^2+\frac{1}{6}Y^3.
\end{align}
\end{examp}

We use these eta-quotients (i.e. \eqref{eq:etaQuotientsUsedToExpressTletar}) as building blocks of linear combinations. Put in other words, we shall give a criterion on whether $F_{r, p^{\beta}}(\tau)$ can be expressed as a finite linear combination of these eta-quotients, and give an algorithm for calculating the coefficients.

Before stating the main theorem of this section, it is natural to ask whether these eta-quotients are linear independent. The answer is yes. Of course, this is well-known. In fact, the weight $0$ weakly holomorphic modular forms $\left(\frac{\eta(p\tau)}{\eta(\tau)}\right)^{\frac{24}{\gcd(12,p-1)}y}$ are widely used in the literature. For example, \cite[Eq.~12]{New55} is essentially the same as the one used here. But we also give a proof here, for the reader's convenience.
\begin{lemm}
\label{lemm:etaQuotLinearIndep}
The functions given by the left-hand side of \eqref{eq:etaQuotientsUsedToExpressTletar} with $y$ ranging over the integers are $\numC$-linear independent.
\end{lemm}
\begin{proof}
It suffices to prove that, the order of these functions at the cusp $\rmi\infty$ (or $1/p$ equivalently) are different for different $y$. Thus using \eqref{eq:etaOrder}, we see that
\begin{equation}
\ord\nolimits_{1/p}\left(\eta^r(p^{\ob}\tau)\left(\frac{\eta(p\tau)}{\eta(\tau)}\right)^{\frac{24}{\gcd(12,p-1)}y}\right)=\frac{p^{\ob}r}{24}+\frac{p-1}{\gcd(12,p-1)}y,
\end{equation}
which is indeed injective as a function of $y$.

\end{proof}

Now comes the main theorem of this section. For a prime $p$, set $\pi_p=\frac{p-1}{\gcd(12,p-1)}$. If there exist integers $y$ such that $P_r\left(\frac{p^{\beta+\ob}-1}{24}r+p^\beta \pi_p y\right) \neq 0$, then let $y_0$ be the smallest integer with this property, and define a sequence $c_y$ for integral $y \geq y_0$ recursively by
\begin{align}
c_{y_0} &= P_r\left(\frac{p^{\beta+\ob}-1}{24}r+p^\beta \pi_p y_0\right),\notag\\
c_y &= P_r\left(\frac{p^{\beta+\ob}-1}{24}r+p^\beta \pi_p y\right)-\sum_{y_0\leq y'< y}c_{y'}a_{\pi_p\cdot(y-y')}^{r,p^\beta,y'},\label{eq:defOfCyLinear}
\end{align}
where $a_{\pi_p(y-y')}^{r,p^\beta,y'}$ was defined in \eqref{eq:etaQuotientsUsedToExpressTletar}. Define an integer $y_1$ as follows. If $r \geq 0$ (or $r<0$ respectively), then $y_1$ is the smallest integer greater than or equal to $y_0$ and $\frac{\gcd(12,p-1)}{24}\frac{p^{\beta-\ob}-1}{p^{\beta-1}(p-1)}r$ (or $\frac{\gcd(12,p-1)}{24}\frac{p^{\beta+1}-p^{1-\ob}}{p-1}(-r)$ respectively).
\begin{thm}
\label{thm:mainThmLinearCombina}
Let $p$ be a prime, $r$ an integer such that $24 \mid r(p^2-1)$. Let $\beta$ be a positive integer, and define $\ob=1$ if $\beta$ is odd, and $\ob=0$ otherwise.
Suppose the smallest integer $y$ such that $P_r\left(\frac{p^{\beta+\ob}-1}{24}r+p^\beta y\right) \neq 0$ is a multiple of $\pi_p$.
Then the following two formulas are equivalent:
\begin{equation}
\label{eq:mainThmLinearCombina}
\sum_{m \in \frac{p^{\ob}r}{24}+\numZ}P_r\left(p^\beta m-\frac{r}{24}\right)q^m=\sum_{y \in \numZ}c_y\eta^r(p^{\ob}\tau)\left(\frac{\eta(p\tau)}{\eta(\tau)}\right)^{\frac{24}{\gcd(12,p-1)}y},
\end{equation}
where the sum in the right-hand side is a finite sum, i.e. only finitely many of $c_y$ are nonzero, and
\begin{multline}
\label{eq:mainThmLinearCombinaEqui}
P_r\left(\frac{p^{\beta+\ob}-1}{24}r+p^\beta \pi_p y_0+p^\beta n\right) \\
=\sum_{\twoscript{y_0\leq y\leq y_0+n/\pi_p}{y \in \numZ}}c_y a_{n-\pi_p(y-y_0)}^{r,p^\beta,y},\qquad n=0,\,1,\,2,\dots,\pi_p(y_1-y_0).
\end{multline}
If \eqref{eq:mainThmLinearCombina} holds, then in fact $c_y=0$ unless $y_0\leq y\leq y_1$. In particular, if $p \in \{2,\,3,\,5,\,7,\,13\}$, then \eqref{eq:mainThmLinearCombina} holds for any $r\in\numZ$ satisfying $24 \mid r(p^2-1)$.
\end{thm}
\begin{proof}
\eqref{eq:mainThmLinearCombina} $\Rightarrow$ \eqref{eq:mainThmLinearCombinaEqui}. Note that if $c_y \neq 0$ in the right-hand side of \eqref{eq:mainThmLinearCombina}, then $y \geq y_0$ by \eqref{eq:defOfCyLinear}. Then \eqref{eq:mainThmLinearCombinaEqui} follows from comparing the coefficients of $q^{p^{\ob}r/24+\pi_py_0+n}$-term of the two sides of \eqref{eq:mainThmLinearCombina}, for $n=0,\,1,\,2,\dots,\pi_p(y_1-y_0)$.

\eqref{eq:mainThmLinearCombinaEqui} $\Rightarrow$ \eqref{eq:mainThmLinearCombina}. We prove more precisely \eqref{eq:mainThmLinearCombina} with the summation range $y\in\numZ$ in the right-hand side replaced by $y_0\leq y\leq y_1$, which simultaneously implies the assertion after \eqref{eq:mainThmLinearCombinaEqui}. Writing down the Fourier expansion, we have
\begin{multline}
\sum_{y_0\leq y\leq y_1}c_y\eta^r(p^{\ob}\tau)\left(\frac{\eta(p\tau)}{\eta(\tau)}\right)^{\frac{24}{\gcd(12,p-1)}y} \\
=\sum_{m\in \frac{p^{\ob}r}{24} + \pi_py_0+\numgeq{Z}{0}}\left(\sum_{\twoscript{n' \geq 0,\,y_0\leq y\leq y_1}{\frac{p^{\ob}r}{24}+\pi_py+n'=m}} c_ya_{n'}^{r,p^\beta,y}\right)q^m.
\end{multline}
Set
\begin{align*}
G(\tau)&=\sum_{m \in \frac{p^{\ob}r}{24}+\numZ}P_r\left(p^\beta m-\frac{r}{24}\right)q^m-\sum_{y_0 \leq y\leq y_1}c_y\eta^r(p^{\ob}\tau)\left(\frac{\eta(p\tau)}{\eta(\tau)}\right)^{\frac{24}{\gcd(12,p-1)}y} \\
&=q^{\frac{p^{\ob}r}{24}+\pi_py_0}\sum_{n \in \numgeq{Z}{0}}\bigg(P_r\left(\frac{p^{\beta+\ob}-1}{24}r+p^\beta \pi_p y_0+p^\beta n\right) \\
&-\sum_{\twoscript{y_0\leq y\leq y_1}{y \leq y_0+n/\pi_p}}c_y a_{n-\pi_p(y-y_0)}^{r,p^\beta,y}\bigg)q^n.
\end{align*}
We shall prove that $G(\tau)$ is identically zero, from which \eqref{eq:mainThmLinearCombina} follows. We prove by contradiction, so assume that $G(\tau)$ is not identically zero. From the discussion after \eqref{eq:etaQuotientsUsedToExpressTletar}, we know that $G(\tau)\in\MFormW{r/2}{p}{v_{\mathbf{r'}}}$, where $v_{\mathbf{r'}}$ is the character of \eqref{eq:etaQuotientsUsedToExpressTletar} for any $y\in\numZ$. By \eqref{eq:mainThmLinearCombinaEqui}, $\ord_{1/p}(G)>\frac{p^{\ob}r}{24}+\pi_py_1$. Next we estimate $\ord_{1/1}(G)$. Break the consideration into two cases: $r\geq 0$ and $r<0$. If $r \geq 0$, we have
\begin{gather}
\ord\nolimits_{1/1}\left(\sum_{m \in \frac{p^{\ob}r}{24}+\numZ}P_r\left(p^\beta m-\frac{r}{24}\right)q^m\right) \geq \frac{r}{24p^\beta} \label{eq:estOrder11}\\
\ord\nolimits_{1/1}\left(\sum_{y_0 \leq y\leq y_1}c_y\eta^r(p^{\ob}\tau)\left(\frac{\eta(p\tau)}{\eta(\tau)}\right)^{\frac{24}{\gcd(12,p-1)}y}\right) \geq \frac{r}{24p^{\ob}}-\frac{\pi_p}{p}y_1,
\end{gather}
by Proposition \ref{prop:orderTl} and \eqref{eq:etaOrder} (both orders may be $+\infty$, which means the corresponding function is identically zero). By the choice of $y_1$ we have $\frac{r}{24p^{\ob}}-\frac{\pi_p}{p}y_1 \leq \frac{r}{24p^\beta}$. Thus $\ord_{1/1}(G) \geq \frac{r}{24p^{\ob}}-\frac{\pi_p}{p}y_1$. We can obtain the same inequality in the case $r<0$, with \eqref{eq:estOrder11} replaced by
\begin{equation}
\ord\nolimits_{1/1}\left(\sum_{m \in \frac{p^{\ob}r}{24}+\numZ}P_r\left(p^\beta m-\frac{r}{24}\right)q^m\right) \geq \frac{rp^\beta}{24}.
\end{equation}
Applying Theorem \ref{thm:valence} (or \eqref{eq:valenceCusp}) to $G(\tau)$, using the estimates of $\ord_{1/p}(G)$ and $\ord_{1/1}(G)$, and using the fact that $\Gzt{p}$ has exactly two cusps, $1/1$ and $1/p$, whose widths are $p$ and $1$ respectively, we obtain
\begin{equation*}
p\left(\frac{r}{24p^{\ob}}-\frac{\pi_p}{p}y_1\right)+\left(\frac{p^{\ob}r}{24}+\pi_py_1\right)<\frac{1}{12}(p+1)\frac{r}{2}.
\end{equation*}
However, after a simplification, the left-hand side is actually $\frac{1}{12}(p+1)\frac{r}{2}$, which contradicts the inequality. Therefore $G(\tau)$ is identically zero.

Finally, if $p \in \{2,\,3,\,5,\,7,\,13\}$, then $\pi_p=1$. Hence if $F_{r,p^\beta}(\tau)$ is not identically zero, then the smallest integer $y$ such that $P_r\left(\frac{p^{\beta+\ob}-1}{24}r+p^\beta y\right) \neq 0$ is really a multiple of $\pi_p$ and \eqref{eq:mainThmLinearCombinaEqui} holds since in this case it is just a reformulation of the definition of $c_y$, i.e. \eqref{eq:defOfCyLinear}. Thus \eqref{eq:mainThmLinearCombina} holds in this case. On the other hand, if $F_{r,p^\beta}(\tau)$ is identically zero, then by definition, $c_y$ are always zero. Therefore \eqref{eq:mainThmLinearCombina} still holds.
\end{proof}

\begin{rema}
Theorem 1.1 of \cite{DLZ19} is essentially the case that $p \in \{2,\,3,\,5,\,7,\,13\}$ of our theorem, while the authors of \cite{DLZ19} prove their Theorem 1.1 using modular equations and Atkin U-operators. They obtain infinite families of congruences of $P_r(n)$ using such identities after analysing linear recursive relations of the coefficients.
\end{rema}

The task to justify whether the smallest integer $y$ such that $P_r\left(\frac{p^{\beta+\ob}-1}{24}r+p^\beta y\right) \neq 0$ is a multiple of $\pi_p$ and whether \eqref{eq:mainThmLinearCombinaEqui} holds, and to compute the coefficients $c_y$ in \eqref{eq:mainThmLinearCombina} can be settled by our SageMath code. See \ref{sec:Usage of code} for the usage. From some experiments we find that when $\pi_p \neq 1$, the assumption on the smallest integer $y$ holds in many cases but \eqref{eq:mainThmLinearCombinaEqui} only holds among a few of these cases. It seems that these cases can be exhausted, which together with the problem of expressing $F_{r,p^\beta}$ as linear combinations of eta-quotients of prime power levels, will be studied in a forthcoming paper. In the remaining of this section, we present some examples.
\begin{examp}
\label{examp:r2Ident}
We present some examples for $r=2$. For $p=5$, we have
\begin{equation}
\sum_{m \in 5^{\ob}/12 + \numgeq{Z}{0}} P_{2}\left(5^\beta m - \frac{1}{12}\right)q^m = (-1)^{(\beta+\ob)/2}\eta^{2}(5^{\ob}\tau)
\end{equation}
for any positive integer $\beta$, since in this case, $y_0=0$, $y_1=1$, and 
\begin{gather*}
c_0=P_2\left(\frac{5^{\beta+\ob}-1}{12}\right)=(-1)^{(\beta+\ob)/2} \\
c_1=P_2\left(\frac{5^{\beta+\ob}-1}{12}+5^\beta\right)-P_2\left(\frac{5^{\beta+\ob}-1}{12}\right)a_1^{2,5^\beta, 0}=0.
\end{gather*}
The quantity $c_0$ was evaluated recursively by setting $l=5^2$, $n=2\cdot 5^{\beta+\ob}$ and $r=2$  in \eqref{eq:mainThmEtarEvenOdd}. To verify that $c_1=0$, note that $a_1^{2,5^\beta, 0}$ is equal to $0$ if $\beta$ is odd, and to $-2$ if $\beta$ is even. Therefore we need to prove that $P_2\left(\frac{5^{\beta+\ob}-1}{12}+5^\beta\right)$ is equal to $0$ if $\beta$ is odd, and to $-2\cdot(-1)^{(\beta+\ob)/2}$ if $\beta$ is even. This can be established inductively by setting $l=5^2$, $n=2\cdot 5^{\beta+\ob}+24\cdot5^\beta$ and $r=2$ in \eqref{eq:mainThmEtarEvenOdd}.

Similarly, for $p=7,\,11$ and any $\beta\in\numgeq{Z}{1}$ we have
\begin{equation}
\sum_{m \in p^{\ob}/12 + \numgeq{Z}{0}} P_{2}\left(p^\beta m - \frac{1}{12}\right)q^m = \eta^{2}(p^{\ob}\tau).
\end{equation}

For $p$ not less than $13$, more terms in linear combinations occur:
\begin{align*}
q^{1/12}\sum_{n \in \numgeq{Z}{0}} P_{2}\left(13^1n + 1\right)q^n &= -2\eta^{2}(\tau) - \eta^{2}(13\tau) \\
q^{1/12}\sum_{n \in \numgeq{Z}{0}} P_{2}\left(13^2n + 14\right)q^n &= 3\eta^{2}(\tau) + 2\eta^{2}(13\tau) \\
q^{1/12}\sum_{n \in \numgeq{Z}{0}} P_{2}\left(13^3n + 183\right)q^n &= -4\eta^{2}(\tau) - 3\eta^{2}(13\tau) \\
q^{17/12}\sum_{n \in \numgeq{Z}{0}} P_{2}\left(17^1n + 24\right)q^n &= -\eta^{2}(17\tau) \\
q^{1/12}\sum_{n \in \numgeq{Z}{0}} P_{2}\left(17^2n + 24\right)q^n &= -\eta^{2}(\tau) \\
q^{1/12}\sum_{n \in \numgeq{Z}{0}} P_{2}\left(37^1n + 3\right)q^n &= 2\eta^{2}(\tau)- \eta^{2}(37\tau).
\end{align*}

\end{examp}

\begin{examp}
\label{examp:rn1Ident}
We present some examples with $r=-1$, where the coefficients are given by the partition function, $P_{-1}$. The quantity $P_{-1}(n)$ equals the number of ways $n$ can be written as a sum of positive integers where repetition is allowed and the order of summands is not taken into account, while $P_{-1}(0)=1$.
\begin{multline*}
q^{23/24}\sum_{n \in \numgeq{Z}{0}} P_{-1}\left(5^2n + 24\right)q^n \\
= 3^2\cdot 5^2\cdot 7\eta^{-7}(\tau)\eta^{6}(5\tau) + 2^2\cdot 5^5\cdot 13\eta^{-13}(\tau)\eta^{12}(5\tau) + 3^2 \cdot 5^7 \cdot 7\eta^{-19}(\tau)\eta^{18}(5\tau) \\
\shoveright{+ 2 \cdot 3 \cdot 5^{10}\eta^{-25}(\tau)\eta^{24}(5\tau) + 5^{12}\eta^{-31}(\tau)\eta^{30}(5\tau).}\\
\shoveleft{q^{19/24}\sum_{n \in \numgeq{Z}{0}} P_{-1}\left(5^3n + 99\right)q^n} \\
\shoveright{= 5^{3} \cdot 1353839\eta^{-6}(\tau)\eta^{5}(5\tau) + 2^{2} \cdot 5^{6} \cdot 471256553\eta^{-12}(\tau)\eta^{11}(5\tau)} \\
\shoveright{+ 2^{2} \cdot 5^{9} \cdot 63499 \cdot 526573\eta^{-18}(\tau)\eta^{17}(5\tau)} \\
\shoveright{+ 3^{2} \cdot 5^{10} \cdot 13 \cdot 179 \cdot 5569 \cdot 465649\eta^{-24}(\tau)\eta^{23}(5\tau)} \\
\shoveright{+ 5^{13} \cdot 17 \cdot 349 \cdot 58346710427\eta^{-30}(\tau)\eta^{29}(5\tau)} \\
\shoveright{+ 2 \cdot 5^{15} \cdot 17 \cdot 563 \cdot 2423 \cdot 121683643\eta^{-36}(\tau)\eta^{35}(5\tau)} \\
\shoveright{+ 3 \cdot 5^{20} \cdot 73 \cdot 1427 \cdot 1403411069\eta^{-42}(\tau)\eta^{41}(5\tau)} \\
\shoveright{+ 3 \cdot 5^{23} \cdot 7^{2} \cdot 38047 \cdot 99960857\eta^{-48}(\tau)\eta^{47}(5\tau)} \\
\shoveright{+ 3 \cdot 5^{25} \cdot 2154851 \cdot 385954601\eta^{-54}(\tau)\eta^{53}(5\tau)} \\
\shoveright{+ 3 \cdot 5^{28} \cdot 79 \cdot 6887225324173\eta^{-60}(\tau)\eta^{59}(5\tau)} \\
\shoveright{+ 3 \cdot 5^{30} \cdot 11 \cdot 283 \cdot 2179 \cdot 198780817\eta^{-66}(\tau)\eta^{65}(5\tau)} \\
\shoveright{+ 2^{2} \cdot 5^{34} \cdot 19 \cdot 37 \cdot 110694925669\eta^{-72}(\tau)\eta^{71}(5\tau)} \\
\shoveright{+ 5^{37} \cdot 19 \cdot 241 \cdot 20669855203\eta^{-78}(\tau)\eta^{77}(5\tau)} \\
\shoveright{+ 2^{2} \cdot 5^{39} \cdot 11 \cdot 43 \cdot 101 \cdot 1151 \cdot 523763\eta^{-84}(\tau)\eta^{83}(5\tau)} \\
\shoveright{+ 2^{2} \cdot 5^{42} \cdot 43 \cdot 79 \cdot 107 \cdot 15580267\eta^{-90}(\tau)\eta^{89}(5\tau)} \\
\shoveright{+ 2^{2} \cdot 3 \cdot 5^{44} \cdot 2789 \cdot 14083 \cdot 38377\eta^{-96}(\tau)\eta^{95}(5\tau)} \\
\shoveright{+ 5^{51} \cdot 13 \cdot 17 \cdot 23 \cdot 223 \cdot 3319\eta^{-102}(\tau)\eta^{101}(5\tau)} \\
\shoveright{+ 3^{2} \cdot 5^{54} \cdot 43 \cdot 1027487\eta^{-108}(\tau)\eta^{107}(5\tau)} \\
\shoveright{+ 5^{56} \cdot 29 \cdot 5864437\eta^{-114}(\tau)\eta^{113}(5\tau)} \\
\shoveright{+ 2^{2} \cdot 5^{58} \cdot 14585089\eta^{-120}(\tau)\eta^{119}(5\tau)} \\
\shoveright{+ 2^{2} \cdot 5^{59} \cdot 97 \cdot 103 \cdot 1979\eta^{-126}(\tau)\eta^{125}(5\tau)} \\
\shoveright{+ 5^{63} \cdot 11 \cdot 29 \cdot 31 \cdot 67\eta^{-132}(\tau)\eta^{131}(5\tau)} \\
\shoveright{+ 5^{66} \cdot 23 \cdot 29 \cdot 31\eta^{-138}(\tau)\eta^{137}(5\tau)} \\
\shoveright{+ 5^{68} \cdot 31 \cdot 73\eta^{-144}(\tau)\eta^{143}(5\tau)} \\
\shoveright{+ 5^{71} \cdot 31\eta^{-150}(\tau)\eta^{149}(5\tau) + 5^{73}\eta^{-156}(\tau)\eta^{155}(5\tau).}\\
\shoveleft{q^{23/24}\sum_{n \in \numgeq{Z}{0}} P_{-1}\left(7^2n + 47\right)q^n} \\
\shoveright{= 2 \cdot 7^{2} \cdot 19 \cdot 67\eta^{-5}(\tau)\eta^{4}(7\tau) + 2 \cdot 7^{4} \cdot 43 \cdot 569\eta^{-9}(\tau)\eta^{8}(7\tau)} \\
\shoveright{+ 7^{5} \cdot 11 \cdot 13 \cdot 9923\eta^{-13}(\tau)\eta^{12}(7\tau) + 2^{5} \cdot 3 \cdot 5^{2} \cdot 7^{7} \cdot 17 \cdot 61\eta^{-17}(\tau)\eta^{16}(7\tau)} \\
\shoveright{+ 2 \cdot 3 \cdot 7^{9} \cdot 23 \cdot 17351\eta^{-21}(\tau)\eta^{20}(7\tau) + 7^{11} \cdot 1437047\eta^{-25}(\tau)\eta^{24}(7\tau)} \\
\shoveright{+ 3^{2} \cdot 7^{12} \cdot 523 \cdot 859\eta^{-29}(\tau)\eta^{28}(7\tau) + 2^{4} \cdot 7^{15} \cdot 11 \cdot 919\eta^{-33}(\tau)\eta^{32}(7\tau)} \\
\shoveright{+ 2^{3} \cdot 3 \cdot 7^{17} \cdot 13 \cdot 103\eta^{-37}(\tau)\eta^{36}(7\tau) + 2 \cdot 3^{2} \cdot 7^{18} \cdot 41 \cdot 43\eta^{-41}(\tau)\eta^{40}(7\tau)} \\
\shoveright{+ 2^{4} \cdot 3 \cdot 5 \cdot 7^{20} \cdot 13\eta^{-45}(\tau)\eta^{44}(7\tau) + 2^{2} \cdot 3 \cdot 7^{22} \cdot 17\eta^{-49}(\tau)\eta^{48}(7\tau)} \\
\shoveright{+ 2^{3} \cdot 7^{24}\eta^{-53}(\tau)\eta^{52}(7\tau) + 7^{25}\eta^{-57}(\tau)\eta^{56}(7\tau).} \\
\shoveleft{q^{11/24}\sum_{n \in \numgeq{Z}{0}} P_{-1}\left(13^1n + 6\right)q^n} \\
\shoveright{= 11\eta^{-2}(\tau)\eta^{1}(13\tau) + 2^{2} \cdot 3^{2} \cdot 13\eta^{-4}(\tau)\eta^{3}(13\tau)} \\
\shoveright{+ 2 \cdot 13^{2} \cdot 19\eta^{-6}(\tau)\eta^{5}(13\tau) + 2^{2} \cdot 5 \cdot 13^{3}\eta^{-8}(\tau)\eta^{7}(13\tau)} \\
\shoveright{+ 2 \cdot 3 \cdot 13^{4}\eta^{-10}(\tau)\eta^{9}(13\tau) + 13^{5}\eta^{-12}(\tau)\eta^{11}(13\tau)} \\
+ 13^{5}\eta^{-14}(\tau)\eta^{13}(13\tau).
\end{multline*}

\end{examp}

\begin{examp}
\label{examp:p19Ident}
We give one more sort of examples with $p=19$ and $\beta=1$.
\begin{align*}
q^{19/24}\sum_{n \in \numgeq{Z}{0}} P_{1}\left(19^1n + 15\right)q^n &= -\eta^{1}(19\tau),\\
q^{19/12}\sum_{n \in \numgeq{Z}{0}} P_{2}\left(19^1n + 30\right)q^n &= \eta^{2}(19\tau),\\
q^{19/8}\sum_{n \in \numgeq{Z}{0}} P_{3}\left(19^1n + 45\right)q^n &= -19\eta^{3}(19\tau),\\
q^{1/6}\sum_{n \in \numgeq{Z}{0}} P_{4}\left(19^1n + 3\right)q^n &= 2^{3}\eta^{4}(\tau) - 19\eta^{4}(19\tau),\\
q^{19/4}\sum_{n \in \numgeq{Z}{0}} P_{6}\left(19^1n + 90\right)q^n &= 19^{2}\eta^{6}(19\tau),\\
q^{1/3}\sum_{n \in \numgeq{Z}{0}} P_{8}\left(19^1n + 6\right)q^n &= 2^{3} \cdot 7\eta^{8}(\tau) - 19^{3}\eta^{8}(19\tau),\\
q^{95/12}\sum_{n \in \numgeq{Z}{0}} P_{10}\left(19^1n + 150\right)q^n &= 19^{4}\eta^{10}(19\tau),\\
q^{1/2}\sum_{n \in \numgeq{Z}{0}} P_{12}\left(19^1n + 9\right)q^n &= 2^{2} \cdot 11 \cdot 19\eta^{12}(\tau) - 19^{5}\eta^{12}(19\tau),\\
q^{2/3}\sum_{n \in \numgeq{Z}{0}} P_{16}\left(19^1n + 12\right)q^n &= 2^{3} \cdot 11 \cdot 13 \cdot 29\eta^{16}(\tau) - 19^{7}\eta^{16}(19\tau),\\
q^{5/6}\sum_{n \in \numgeq{Z}{0}} P_{20}\left(19^1n + 15\right)q^n &= 2^{3} \cdot 11 \cdot 23 \cdot 397\eta^{20}(\tau) - 19^{9}\eta^{20}(19\tau),\\
q\sum_{n \in \numgeq{Z}{0}} P_{24}\left(19^1n + 18\right)q^n &= 2^{2} \cdot 5 \cdot 7^{2} \cdot 11 \cdot 23 \cdot 43\eta^{24}(\tau) - 19^{11}\eta^{24}(19\tau).
\end{align*}
\end{examp}

\section{Miscellaneous Observations}
\label{sec:Miscellaneous}
\subsection{More Identities Involving Eta-quotients}
\label{subsec:MoreIdent}
The right-hand sides of all identities occurring in Section \ref{sec:Express} (i.e. Example \ref{examp:r2Ident}, Example \ref{examp:rn1Ident} and Example \ref{examp:p19Ident}) are some linear combinations of eta-quotients of prime levels. Here we present some identities of non-prime levels which can also be proved by Theorem \ref{thm:charCompatibleEtaDirichlet}. We have
\begin{align*}
\eta(\tau)\eta^{-2}(2\tau)\eta^5(4\tau)\eta^{-2}(8\tau)&=q^{11/24}\sum_{n\in\numgeq{Z}{0}}P_2(2n)q^n,\\
\eta(\tau)\eta^{-1}(4\tau)\eta^{2}(8\tau)&=-\frac{1}{2}q^{13/24}\sum_{n\in\numgeq{Z}{0}}P_2(2n+1)q^n,\\
\eta^{-5}(\tau)\eta^{5}(4\tau)\eta^{-2}(8\tau)&=q^{-1/24}\sum_{n\in\numgeq{Z}{0}}P_{-2}(2n)q^n,\\
\eta^{-5}(\tau)\eta^2(2\tau)\eta^{-1}(4\tau)\eta^{2}(8\tau)&=\frac{1}{2}q^{11/24}\sum_{n\in\numgeq{Z}{0}}P_{-2}(2n+1)q^n,\\
\eta^{-2}(\tau)\eta^{10}(2\tau)\eta^{-4}(4\tau)&=q^{1/12}\sum_{n\in\numgeq{Z}{0}}P_{4}(2n)q^n,\\
\eta^{2}(\tau)\eta^{-2}(2\tau)\eta^{4}(4\tau)&=-\frac{1}{4}q^{7/12}\sum_{n\in\numgeq{Z}{0}}P_{4}(2n+1)q^n,\\
\eta^{-14}(\tau)\eta^{14}(2\tau)\eta^{-4}(4\tau)&=q^{-1/12}\sum_{n\in\numgeq{Z}{0}}P_{-4}(2n)q^n,\\
\eta^{-10}(\tau)\eta^{2}(2\tau)\eta^{4}(4\tau)&=\frac{1}{4}q^{5/12}\sum_{n\in\numgeq{Z}{0}}P_{-4}(2n+1)q^n.
\end{align*}
All the above can be proved by first observing that on $\Gzt{8}$, the character of the eta-power $\eta^r$ related to the right-hand side of each identity and that of the eta-quotient in the left-hand side are $\widetilde{\tbtmat{1}{0}{0}{2}}$-compatible (using Theorem \ref{thm:charCompatibleEtaDirichlet}), and then using the Sturm bound, or Theorem \ref{thm:valence} for the difference between the left-hand side and a suitable multiple of $T_2\eta^r$. Also, related to suitable operators $T_4$ on $\Gzt{8}$, we have following identities:
\begin{align}
\eta^{5}(\tau)\eta^{-4}(2\tau)\eta^{5}(4\tau)\eta^{-2}(8\tau)&=q^{1/24}\sum_{n\in\numgeq{Z}{0}}P_{4}(4n)q^n, \label{eq:N8P40}\\
\eta^{-1}(\tau)\eta^{2}(2\tau)\eta^{5}(4\tau)\eta^{-2}(8\tau)&=-\frac{1}{4}q^{7/24}\sum_{n\in\numgeq{Z}{0}}P_{4}(4n+1)q^n, \label{eq:N8P41}\\
\eta^{5}(\tau)\eta^{-2}(2\tau)\eta^{-1}(4\tau)\eta^{2}(8\tau)&=\frac{1}{2}q^{13/24}\sum_{n\in\numgeq{Z}{0}}P_{4}(4n+2)q^n, \label{eq:N8P42}\\
\eta^{-1}(\tau)\eta^{4}(2\tau)\eta^{-1}(4\tau)\eta^{2}(8\tau)&=\frac{1}{8}q^{19/24}\sum_{n\in\numgeq{Z}{0}}P_{4}(4n+3)q^n. \label{eq:N8P43}
\end{align}
From these four identities, some linear relations on eta-quotients and some arithmetical identities on $P_4(n)$ follow. For example, adding suitable multiples of the four identities together, we obtain
\begin{multline*}
\eta^4(\frac{\tau}{4})=\eta^{5}(\tau)\eta^{-4}(2\tau)\eta^{5}(4\tau)\eta^{-2}(8\tau)-4\eta^{-1}(\tau)\eta^{2}(2\tau)\eta^{5}(4\tau)\eta^{-2}(8\tau)\\
+2\eta^{5}(\tau)\eta^{-2}(2\tau)\eta^{-1}(4\tau)\eta^{2}(8\tau)+8\eta^{-1}(\tau)\eta^{4}(2\tau)\eta^{-1}(4\tau)\eta^{2}(8\tau).
\end{multline*}
For another example, by multiplying \eqref{eq:N8P40} and \eqref{eq:N8P43}, multiplying \eqref{eq:N8P41} and \eqref{eq:N8P42}, and comparing two results, we obtain
\begin{equation*}
\sum_{\twoscript{n_1,n_2\in\numZ}{n_1+n_2=n}}P_4(4n_1)P_4(4n_2+3)=-\sum_{\twoscript{n_1,n_2\in\numZ}{n_1+n_2=n}}P_4(4n_1+1)P_4(4n_2+2),
\end{equation*}
for any integer $n$.

It should be mentioned that some authors have recently obtained equalities similar to \eqref{eq:N8P40} -- \eqref{eq:N8P43} and others using modular equations and Atkin $U_l$ operators studied in \cite{AL70}. See, for instance, \cite{CT10}, \cite{BO12} and \cite{XY13}.

\subsection{A Further Application of Theorem \ref{thm:charCompatibleEtaDirichlet} Involving Theta Functions}
\label{subsec:appMainThmTheta}
Let $L$ be a positive definite lattice of even rank $r$, that is to say, $L$ is a free $\numZ$-module of rank $r$, and there is a positive definite quadratic form $Q$ defined on $L$. We require that $Q$ is even integral, i.e., $Q(L)\subseteq\numZ$. Set $D=\det(Q)=\det\left(B(e_i,e_j)\right)_{1 \leq i,j \leq r}$ and $\Delta=(-1)^{r/2}D$ where $\{e_i\}$ is a $\numZ$-basis of $L$ and $B(x,y)=Q(x+y)-Q(x)-Q(y)$. The theta function associated with $(L, Q)$ is
\begin{equation}
\vartheta(\tau)=\sum_{X \in L}\etp{\tau Q(X)}.
\end{equation}
A well-known result of Schoeneberg \cite{Sch39} says that $\vartheta\in\MForm{r/2}{D}{\chi_\Delta}$, where $\chi_\Delta(d)=\legendre{\Delta}{d}$. It can be shown that $\chi_\Delta$ is actually a Dirichlet character modulo $D$, since $\Delta \equiv 0,1 \bmod 4$. (For a proof, see \cite[Lemma 14.3.20]{CoS17}.) We want to seek characters $v$ and operators $T_l$ from $\MForm{r/2}{D}{v}$ to $\MForm{r/2}{D}{\chi_\Delta}$.
\begin{lemm}
\label{lemm:vrchiCompa}
Let $r$ be an even positive integer, and $L$ be a lattice of rank $r$, with positive definite quadratic form $Q$. Suppose $Q$ is even integral. Let $D$ and $\Delta$ be as above. Let $\mathbf{r}=(r_1, r_2,\dots)$ be an integral sequence such that $r_n=0$ unless $n \mid D$, and the sum of all components $r_n$ is zero. Let $l$ be a positive integer. Then the characters $v_{\mathbf{r}}$ and $\chi_\Delta$ are $\widetilde{\tbtmat{1}{0}{0}{l}}$-compatible if and only if the following conditions hold:
\begin{gather}
l\cdot\sum_{n\mid D}\frac{D}{n}r_n \equiv 0 \pmod{24}, \label{eq:vrchi1}\\
\sum_{n\mid D}nr_n \equiv 0 \pmod{24}, \label{eq:vrchi2}\\
2 \nmid Dl \Rightarrow \prod_{2 \nmid r_n}n \equiv 1 \bmod 4, \label{eq:vrchi3}\\
\radO(D) = \radO\left(\prod\nolimits_{2\nmid r_n}n\right) \text{ and } 4\mid r. \label{eq:vrchi4}
\end{gather}
In this case, $T_l$ maps $\MForm{r/2}{D}{v_{\mathbf{r}}}$ to $\MForm{r/2}{D}{\chi_\Delta}$.
\end{lemm}
\begin{proof}
An immediate application of Theorem \ref{thm:charCompatibleEtaDirichlet}, by inserting the zero sequence for $\mathbf{r'}$, $D$ for $N$, the trivial character modulo $D$ for $\chi$ and $\chi_\Delta$ for $\chi'$ in that theorem.
\end{proof}
For the simplest case $4 \mid r$ and $D=p$, where $p$ is a prime congruent to $1$ modulo $8$, all possible pairs $(r_1,r_p)$ subject to the conditions in Lemma \ref{lemm:vrchiCompa} are exactly the following
\begin{equation*}
r_1=2t+1,\qquad r_p=-(2t+1),
\end{equation*}
where $t\in\numZ$ satisfies $3\mid \frac{p-1}{8}(2t+1)$ and $l$ can be any positive integer. Nevertheless, since $v_\mathbf{r}$ and $\chi_\Delta$ are $\widetilde{\tbtmat{1}{0}{0}{1}}$-compatible, they are actually the same character on $\Gzt{D}$. Thus in this case we have nothing new. (That is to say, in this case, the operators $T_l$ are essentially the same as the usual Hecke operators.) So is the case $D=p^2$ with $p>3$.

Next consider the case $D=3^2$. It can be shown that all tuples $(r_1,r_3,r_9)$ subject to the conditions in Lemma \ref{lemm:vrchiCompa} are exactly the following
\begin{align*}
r_1&=3t, \\
r_3&=4s, \\
r_{9}&=-3t-4s,
\end{align*}
with $t,s\in\numZ$. In this case $l$ should be a positive integer satisfying $3 \mid ls$. Now the character $\chi_9$ is trivial on $\Gzt{9}$, and\
\begin{equation}
\label{eq:someGamma09char}
v_{\mathbf{r}}\widetilde{\tbtmat{a}{b}{c}{d}}=\etp{\frac{2cds}{27}},
\end{equation}
which can be proved by \eqref{eq:etaChar} and \eqref{eq:etaQuoChar}. Note that when $s \equiv 1, 2 \bmod 3$, the character $v_{\mathbf{r}}$ can not be induced by any Dirichlet character, since its values depend on $c$. We summarize:
\begin{prop}
Let $v$ be the character given by \eqref{eq:someGamma09char} (which depends on $s\in\numZ$). Let $r$ be a positive integer divisible by $4$. If $s \equiv 1, 2 \bmod 3$, then for any positive integer $l$ divisible by $3$, the operator $T_l$ maps $\MForm{r/2}{9}{v}$ to $\MFormt{r/2}{9}$. On the other hand, if $s \equiv 0 \bmod 3$, then $T_l$ maps $\MFormt{r/2}{9}$ to itself, for any positive integer $l$.
\end{prop}

We emphasize that $T_l$ depends not only on $l$, but also on the source character, the target character, and the weight. See the paragraph after Definition \ref{deff:doubleCosetOperators}. It is an interesting question that whether we can find some functions $f_1,\,f_2,\dots$ in $\MForm{r/2}{9}{v}$ of simple Fourier coefficients and some operators $T_{l_1},\,T_{l_2},\dots$ such that $\vartheta=T_{l_1}f_1+T_{l_2}f_2+\dots$. This question and questions in more general cases (including theta functions of odd rank) will be treated in a later paper, as further applications of generalized double coset operators.

\appendix

\section{Tables}
\label{apx:tables}

\begin{longtable}{llllll}
\caption{For each $\eta^{\mathbf{r}}$ and $\eta^{\mathbf{r'}}$ in Table \ref{table:EtaDimonePrime} or Table \ref{table:EtaDimoneNonPrime} of integral weights, the set of positive integers $l$ such that $v_{\mathbf{r}}$ and $v_{\mathbf{r'}}$ are $\widetilde{\tbtmat{1}{0}{0}{l}}$-compatible (the variable $j$ ranges over non-negative integers) \label{table:Lintegral}} \\
\toprule
$\eta^{\mathbf{r}}$ & $\eta^{\mathbf{r'}}$ & $l$ & $\eta^{\mathbf{r}}$ & $\eta^{\mathbf{r'}}$ & $l$ \\
\midrule
\endfirsthead
\multicolumn{6}{l}{Table \ref{table:Lintegral}. (continued)} \\
\toprule
$\eta^{\mathbf{r}}$ & $\eta^{\mathbf{r'}}$ & $l$ & $\eta^{\mathbf{r}}$ & $\eta^{\mathbf{r'}}$ & $l$ \\
\midrule
\endhead
\bottomrule
\multicolumn{6}{r}{To be continued} \\[2ex]
\endfoot
\bottomrule
\endlastfoot
\multicolumn{6}{l}{$N=1$, $k=1,\,2,\dots,11$} \\
$\eta^{2}_{1}$&$\eta^{2}_{1}$&$1+12j$&$\eta^{4}_{1}$&$\eta^{4}_{1}$&$1+6j$\\ $\eta^{6}_{1}$&$\eta^{6}_{1}$&$1+4j$&$\eta^{8}_{1}$&$\eta^{8}_{1}$&$1+3j$\\ $\eta^{10}_{1}$&$\eta^{10}_{1}$&$1+12j$&$\eta^{12}_{1}$&$\eta^{12}_{1}$&$1+2j$\\ $\eta^{14}_{1}$&$\eta^{14}_{1}$&$1+12j$&$\eta^{16}_{1}$&$\eta^{16}_{1}$&$1+3j$\\ $\eta^{18}_{1}$&$\eta^{18}_{1}$&$1+4j$&$\eta^{20}_{1}$&$\eta^{20}_{1}$&$1+6j$\\ $\eta^{22}_{1}$&$\eta^{22}_{1}$&$1+12j$& & & \\
\\
\end{longtable}

{\footnotesize
\begin{longtable}{llllllll}
\caption{For each $\eta^{\mathbf{r}}$ and $\eta^{\mathbf{r'}}$ in Table \ref{table:EtaDimonePrime} or Table \ref{table:EtaDimoneNonPrime} of half-integral weights, the positive integer $l_0$ and and the set of $m \in \numgeq{Z}{1}$ such that $v_{\mathbf{r}}$ and $v_{\mathbf{r'}}$ are $\widetilde{\tbtmat{1}{0}{0}{l_0m^2}}$-compatible \label{table:Lhalfintegral}} \\
\toprule
$\eta^{\mathbf{r}}$ & $\eta^{\mathbf{r'}}$ & $l_0$ & $m^2 \bmod{24}$ & $\eta^{\mathbf{r}}$ & $\eta^{\mathbf{r'}}$ &  $l_0$ & $m^2 \bmod{24}$ \\
\midrule
\endfirsthead
\multicolumn{8}{l}{Table \ref{table:Lhalfintegral}. (continued)} \\
\toprule
$\eta^{\mathbf{r}}$ & $\eta^{\mathbf{r'}}$ & $l_0$ & $m^2 \bmod{24}$ & $\eta^{\mathbf{r}}$ & $\eta^{\mathbf{r'}}$ &  $l_0$ & $m^2 \bmod{24}$ \\
\midrule
\endhead
\bottomrule
\multicolumn{8}{r}{To be continued} \\[2ex]
\endfoot
\bottomrule
\endlastfoot
\multicolumn{8}{l}{$N=1$, $k=1/2,\,3/2,\dots,23/2$}\\
$\eta^{1}_{1}$&$\eta^{1}_{1}$&$1$&$1$&$\eta^{3}_{1}$&$\eta^{3}_{1}$&$1$&$1,9$\\ $\eta^{5}_{1}$&$\eta^{5}_{1}$&$1$&$1$&$\eta^{7}_{1}$&$\eta^{7}_{1}$&$1$&$1$\\ $\eta^{9}_{1}$&$\eta^{9}_{1}$&$1$&$1,9$&$\eta^{11}_{1}$&$\eta^{11}_{1}$&$1$&$1$\\ $\eta^{13}_{1}$&$\eta^{13}_{1}$&$1$&$1$&$\eta^{15}_{1}$&$\eta^{15}_{1}$&$1$&$1,9$\\ $\eta^{17}_{1}$&$\eta^{17}_{1}$&$1$&$1$&$\eta^{19}_{1}$&$\eta^{19}_{1}$&$1$&$1$\\ $\eta^{21}_{1}$&$\eta^{21}_{1}$&$1$&$1,9$&$\eta^{23}_{1}$&$\eta^{23}_{1}$&$1$&$1$\\
\\
\end{longtable}
}

\begin{longtable}{llllll}
\caption{For each holomorphic $\eta^{\mathbf{r}}$ and $\eta^{\mathbf{r'}}$ on $\Gzt{4}$ and of weight $2$ that are also cusp forms, the set of positive integers $l$ such that $v_{\mathbf{r}}$ and $v_{\mathbf{r'}}$ are $\widetilde{\tbtmat{1}{0}{0}{l}}$-compatible (the variable $j$ ranges over non-negative integers) \label{table:N4k2L}} \\
\toprule
$\eta^{\mathbf{r}}$ & $\eta^{\mathbf{r'}}$ & $l$ & $\eta^{\mathbf{r}}$ & $\eta^{\mathbf{r'}}$ & $l$ \\
\midrule
\endfirsthead
\multicolumn{6}{l}{Table \ref{table:N4k2L}. (continued)} \\
\toprule
$\eta^{\mathbf{r}}$ & $\eta^{\mathbf{r'}}$ & $l$ & $\eta^{\mathbf{r}}$ & $\eta^{\mathbf{r'}}$ & $l$ \\
\midrule
\endhead
\bottomrule
\multicolumn{6}{r}{To be continued} \\[2ex]
\endfoot
\bottomrule
\endlastfoot
$\eta^{-7}_{1}$$\eta^{18}_{2}$$\eta^{-7}_{4}$&$\eta^{-7}_{1}$$\eta^{18}_{2}$$\eta^{-7}_{4}$&$1+24j$&$\eta^{-7}_{1}$$\eta^{18}_{2}$$\eta^{-7}_{4}$&$\eta^{-3}_{1}$$\eta^{10}_{2}$$\eta^{-3}_{4}$&$5+24j$\\ $\eta^{-7}_{1}$$\eta^{18}_{2}$$\eta^{-7}_{4}$&$\eta^{-1}_{1}$$\eta^{6}_{2}$$\eta^{-1}_{4}$&$7+24j$&$\eta^{-7}_{1}$$\eta^{18}_{2}$$\eta^{-7}_{4}$&$\eta^{3}_{1}$$\eta^{-2}_{2}$$\eta^{3}_{4}$&$11+24j$\\ $\eta^{-6}_{1}$$\eta^{15}_{2}$$\eta^{-5}_{4}$&$\eta^{-6}_{1}$$\eta^{15}_{2}$$\eta^{-5}_{4}$&$1+24j$&$\eta^{-6}_{1}$$\eta^{15}_{2}$$\eta^{-5}_{4}$&$\eta^{-5}_{1}$$\eta^{15}_{2}$$\eta^{-6}_{4}$&$4+24j$\\ $\eta^{-6}_{1}$$\eta^{15}_{2}$$\eta^{-5}_{4}$&$\eta^{-3}_{1}$$\eta^{9}_{2}$$\eta^{-2}_{4}$&$4+24j$&$\eta^{-6}_{1}$$\eta^{15}_{2}$$\eta^{-5}_{4}$&$\eta^{-2}_{1}$$\eta^{9}_{2}$$\eta^{-3}_{4}$&$7+24j$\\ $\eta^{-6}_{1}$$\eta^{15}_{2}$$\eta^{-5}_{4}$&$\eta^{-1}_{1}$$\eta^{3}_{2}$$\eta^{2}_{4}$&$4+24j$&$\eta^{-6}_{1}$$\eta^{15}_{2}$$\eta^{-5}_{4}$&$\eta^{1}_{1}$$\eta^{-3}_{2}$$\eta^{6}_{4}$&$4+24j$\\ $\eta^{-6}_{1}$$\eta^{15}_{2}$$\eta^{-5}_{4}$&$\eta^{2}_{1}$$\eta^{3}_{2}$$\eta^{-1}_{4}$&$13+24j$&$\eta^{-6}_{1}$$\eta^{15}_{2}$$\eta^{-5}_{4}$&$\eta^{6}_{1}$$\eta^{-3}_{2}$$\eta^{1}_{4}$&$19+24j$\\ $\eta^{-6}_{1}$$\eta^{16}_{2}$$\eta^{-6}_{4}$&$\eta^{-6}_{1}$$\eta^{16}_{2}$$\eta^{-6}_{4}$&$1+12j$&$\eta^{-6}_{1}$$\eta^{16}_{2}$$\eta^{-6}_{4}$&$\eta^{2}_{1}$$\eta^{0}_{2}$$\eta^{2}_{4}$&$5+12j$\\ $\eta^{-5}_{1}$$\eta^{12}_{2}$$\eta^{-3}_{4}$&$\eta^{-5}_{1}$$\eta^{12}_{2}$$\eta^{-3}_{4}$&$1+24j$&$\eta^{-5}_{1}$$\eta^{12}_{2}$$\eta^{-3}_{4}$&$\eta^{-3}_{1}$$\eta^{12}_{2}$$\eta^{-5}_{4}$&$7+24j$\\
\end{longtable}

\section{Usage of Sage Code}
\label{sec:Usage of code}
The Sage code can be obtained from the repository \cite{Zhu23}. It is an ipynb file which should be opened in a Sage Jupyter notebook. The code in the first cell should be run when the file is opened, which contains all Python/Sage functions dealing with eta-quotients and compatibility. We list the usage of some of them below.
\begin{itemize}
    \item To calculate the character of an eta-quotient \eqref{eq:etaQuoChar}, use \lstinline{eta_char_gamma0N_exp(N, es, *A)}. For example, a call \lstinline!eta_char_gamma0N_exp(4, {1: 2, 2: 3, 4: -5}, 1, 0, 4, 1)! returns $v_{\eta^2(\tau)\eta^3(2\tau)\eta^{-5}(4\tau)}\widetilde{\tbtmat{1}{0}{4}{1}}$. This function uses the character formula of eta function involving Dedekind sums, not \eqref{eq:etaChar}, to allow fractional powers.
    \item To calculate the order of an eta-quotient, use \lstinline{eta_order_gamma0N(N, es, cusp)}. For example, a call \lstinline!eta_order_gamma0N(4, {1: 2, 2: 3, 4: -5}, 1/4)! returns the order of $\eta^2(\tau)\eta^3(2\tau)\eta^{-5}(4\tau)$ at the cusp $1/4=\rmi\infty$.
    \item To calculate the quantity $P_r(n)$ in \eqref{eq:Prn}, use \lstinline{p_num(r, n)}. If a sequence $P_r(n)$ for $n \leq n_1$ is needed, one should call \lstinline{p_num_init(r, n1)} first to increase the speed.
    \item To calculate the coefficient $a_n$ in the expression $q^{-\sum_{n}nr_n/24}\eta^{\mathbf{r}}=\sum_{n \in \numgeq{Z}{0}}a_nq^n$, use \lstinline{coeff_eta_quotient(N, es, n)}. For example, a call \lstinline!coeff_eta_quotient(4, {1: 2, 2: 3, 4: -5}, 3)! returns $a_3$ in $q^{1/2}\eta^2(\tau)\eta^3(2\tau)\eta^{-5}(4\tau)=\sum_{n \in \numgeq{Z}{0}}a_nq^n$. The coefficients in Table \ref{table:C3998+13709n} are obtained using this function. To acquire the whole series including the first fractional exponent, call \lstinline!eta_gamma_N(4, {1: 2, 2: 3, 4: -5}, 10)!, where the argument $10$ means we need a power series of $O(q^{10})$.
    \item To check Theorem \ref{thm:charCompatibleEtaDirichlet}, call \lstinline!check_compatible(N, dchar_left, echar_left, dchar_right, echar_right, l)!, where the arguments represent $N$, $\chi$, $\mathbf{r}$, $\chi'$, $\mathbf{r'}$, $l$ respectively. This function uses the method described in Remark \ref{rema:charCompatible}, not the conditions \eqref{eq:charCom1}, \eqref{eq:charCom2}, \eqref{eq:charCom3}, \eqref{eq:charCom4}. The ipynb file includes an example listing all $1 \leq l \leq 48$ such that $\eta_1^2\eta_2^{-3}\eta_4^5$ and $\eta_1^3\eta_2^3\eta_4^{-2}$ are $\widetilde{\tbtmat{1}{0}{0}{l}}$-compatible.
    \item To generate data in Table \ref{table:NkDimone} and Table \ref{table:Nk12}, use \lstinline!gen_etaquot_dimone(b_cusp=False)! or \lstinline!gen_etaquot_dimone(b_cusp=True)!.
    \item To generate data in Table \ref{table:EtaDimonePrime} and Table \ref{table:EtaDimoneNonPrime}, use \lstinline{gen_holo_etaquot(N, k)}.
    \item To generate data in Table \ref{table:Lintegral}, Table \ref{table:Lhalfintegral} and Table \ref{table:N4k2L}, use \lstinline!findl(N, k)! or \lstinline!findl(N, k, b_cusp=True)!. This function can also produce those compatible operators not given in the paper. For example, The ipynb file includes a cell showing the result of \lstinline!findl(6, 1, omitempty=True, b_cusp=True)!, which contains many pairs of eta-quotients of level $6$ and weight $1$ and compatible operators.
    \item There are two functions calculating the Gauss sum in Lemma \ref{lemm:aGaussSum}. The first is \lstinline!some_gauss_sum(m, t)! which calculates using definition, and the second is \lstinline!some_gauss_sum_formula_all(m, t)!, which calculates using the right-hand side of \eqref{eq:aGaussSum}. The reader may check Lemma \ref{lemm:aGaussSum} by comparing the results returned by these two functions.
    \item There are two functions calculating the coefficient of the $q^{n/(24l)}$-term in \eqref{eq:FourierTletar}. One is \lstinline!coef_Tletar_new(l, r, n)! which calculates by definition; the other is \lstinline!coef_Tletar_good_formula(l, r, n)! which calculates using Lemma \ref{lemm:FourierCoeffReven} or Lemma \ref{lemm:FourierCoeffRodd}. By comparing the results returned by the two functions, the reader may check the validity of these two lemmas.
    \item To check Theorem \ref{thm:mainThmEtar} numerically, use \lstinline!check_main_theorem_onTletar(l, r, n_min, n_max)!. This would compare the two sides of \eqref{eq:mainThmEtarEvenOdd} for $n$ in \lstinline!range(n_min,n_max+1)!.
    \item To evaluate the coefficient $a_n^{r,p^\beta, y}$ in \eqref{eq:etaQuotientsUsedToExpressTletar}, call \lstinline!some_eta_quotient_coef(r, p, y, n, beta)!. To check Lemma \ref{lemm:crpynRecur}, call \lstinline!check_some_eta_quot_coef_recur(r, p, y, n, beta)!.
    \item To check whether the condition \eqref{eq:mainThmLinearCombinaEqui} holds, call \lstinline!check_cy(r, p, beta)!. When it holds, call \lstinline{show_identity_latex(r, p, beta)} to print the resulted identity \eqref{eq:mainThmLinearCombina}.
\end{itemize}

\section*{Acknowledgements}
Both authors are supported by the National Natural Science Foundation of China, Grant no. 12271405, to which we express our thanks. We appreciate the referee's careful, patient reading and professional, detailed suggestions. The `only if' part of Theorem \ref{thm:charCompatibleEtaDirichlet} is encouraged by the referee.


\bibliographystyle{elsarticle-num} 
\bibliography{main}


\end{document}